\input epsf

\documentclass[12pt]{article}

\usepackage{latexsym,amsmath,amssymb,amsfonts,amscd,multirow,bm}
\usepackage{epsfig,subfigure,verbatim,epstopdf,graphics}
\usepackage{color}
\usepackage{slashbox}
\newtheorem{thm}{Theorem}[section]

\newtheorem{lem}[thm]{Lemma}

\usepackage{bibentry}
\nobibliography*

\usepackage{scalerel}

\newcommand\reallywidehat[1]{%
\begin{array}{c}
\stretchto{
 \scaleto{
   \scalerel*[\widthof{$#1$}]{\bigwedge}
   {\rule[-\textheight/2]{1ex}{\textheight}} 
 }{1.0\textheight} 
}{0.75ex}\\           
#1                 
\end{array}
}

\newcommand{\bx}{\mathbf{x}}

\newcommand{\bE}{\mathbf{E}}

\newfont{\iams}{msbm9}

\newcommand{\commentbis}[1]{}
\newcommand{\be}{\begin{eqnarray}}
\newcommand{\ee}{\end{eqnarray}}
\newcommand{\beno}{\begin{eqnarray*}}
\newcommand{\eeno}{\end{eqnarray*}}
\newcommand{\barr}[1]{\begin{array}{#1}}
\newcommand{\earr}{\end{array}}

\newcommand{\Dt}{\Delta t}
\newcommand{\df}{\partial}

\newcommand{\mE}{{\mathcal E}}
\def\Real{\mathbb R}
\newcommand{\beq}{\begin{equation}}
\newcommand{\eeq}{\end{equation}}
\newcommand{\beqa}{\begin{eqnarray}}
\newcommand{\eeqa}{\end{eqnarray}}

\newcommand{\Kx}{{K_x}}
\newcommand{\Kv}{{K_v}}
\newcommand{\bv}{{\bf v}}
\newcommand{\bB}{{\bf B}}
\newcommand{\bJ}{{\bf J}}

\newcommand{\bU}{{\bf U}}
\newcommand{\bV}{{\bf V}}
\newcommand{\bn}{{\bf n}}
\newcommand{\Ox}{{\Omega_x}}
\newcommand{\Ov}{{\Omega_v}}

\newcommand{\nh}{{n+1/2}}
\newcommand{\na}{{n+1}}
\newcommand{\mT}{{\mathcal T}}
\newcommand{\mG}{{\mathcal G}}
\newcommand{\mS}{{\mathcal S}}
\newcommand{\mU}{{\mathcal U}}
\newcommand{\mW}{{\mathcal W}}
\newcommand{\mZ}{{\mathcal Z}}
\newcounter{nsez}

\title
{  Energy-conserving Discontinuous Galerkin Methods for the Vlasov-Maxwell System}

\author{  Yingda Cheng
\thanks{Department of Mathematics, Michigan State University,
East Lansing, MI 48824 U.S.A.
 {\tt ycheng@math.msu.edu}}
  \and
Andrew J. Christlieb
\thanks{Department of Mathematics, Michigan State University,
East Lansing, MI 48824 U.S.A.
 {\tt christlieb@math.msu.edu}}
   \and
Xinghui Zhong
\thanks{Corresponding author. Department of Mathematics, Michigan State University,
East Lansing, MI 48824 U.S.A.
 {\tt zhongxh@math.msu.edu}}
}

\date{\today}

\begin{document}


\maketitle

\begin{abstract}
In this paper, we  generalize the idea in our previous work for the Vlasov-Amp\`{e}re (VA) system [\bibentry{cheng_va}] and develop energy-conserving discontinuous Galerkin (DG) methods for the Vlasov-Maxwell (VM) system. The VM system is a fundamental model in the simulation of collisionless magnetized plasmas. Compared to [\bibentry{cheng_va}], additional care needs to be taken for both the temporal and spatial discretizations to achieve similar type of conservation when the magnetic field is no longer negligible. Our proposed schemes conserve the total particle number and the total energy at the same time,  therefore can obtain accurate and physically relevant numerical solutions. The main components of our methods include second order and above, explicit or implicit energy-conserving temporal discretizations, and  DG methods for Vlasov and Maxwell's equations with carefully chosen numerical fluxes.   Benchmark numerical tests such as the streaming Weibel instability are provided to validate the accuracy and conservation of the schemes.

\end{abstract}



\section{Introduction}
\label{sec:intro}
\setcounter{equation}{0}
\setcounter{figure}{0}
\setcounter{table}{0}
In this paper, we develop energy-conserving numerical schemes for Vlasov-Maxwel (VM) systems. The VM system is an important  equation for the  modeling of collisionless magnetized plasmas. In this model, the Vlasov equation  describes the time  evolution of the probability distribution function of collisionless charged particles with long-range interactions. 
The evolution of the electromagnetic field
is modeled by the Maxwell's equation.
Here we restrict our attention to the VM equation for a single species of nonrelativistic electrons   while the ions are treated as uniform fixed background. Under the scaling of the characteristic time by the inverse of the plasma frequency $\omega_p^{-1}$, length scaled by the Debye length $\lambda_D$, and characteristic electric and magnetic field by $\bar{E}=\bar{B}=-m c \omega_p/e$, the dimensionless  VM equations become
 \begin{eqnarray}
&&\partial_t f  +  \bv \cdot \nabla_{\bx} f  + (\bE + \bv \times \bB) \cdot \nabla_\bv f = 0~, \qquad (\bx, \bv) \in \Omega= \Omega_x \times \Ov  \notag\\
&&\frac{\df \bE}{\df t} = \nabla_\bx \times \bB - \bJ, \quad
\frac{\df \bB}{\df t} = -  \nabla_\bx \times \bE~, \qquad \bx \in \Ox \quad\label{modelvm} \\
&&\nabla_\bx \cdot \bE =  \rho-\rho_i, \qquad
\nabla_\bx \cdot \bB = 0~, \quad\notag
\end{eqnarray}
with the density and current density defined by
\begin{equation*}
\rho(\bx, t)= \int_\Ov f(\bx, \bv, t)d\bv,\qquad \bJ(\bx, t)= \int_\Ov f(\bx, \bv, t)\,\bv d\bv,
\end{equation*}
and $\rho_i$ being the ion density.
In this model, $f=f(\bx,\bv,t)$ is  the probability distribution function ($pdf$) for  finding an electron at position $\bx$ with velocity $\bv$ at time $t$. $\Ox$ denotes the physical domain, while $\Ov=\Real^n$ represents the velocity domain.
It is well-known that the VM system conserves the total particle number $\int_\Ox \int_\Ov f \, d\bv d\bx$,     and  the total energy $$TE=\frac{1}{2} \int_\Ox \int_\Ov f |\bv|^2 d\bv d\bx + \frac{1}{2}\int_\Ox |\bE|^2+|\bB|^2 d\bx, $$
which is composed of the kinetic and electromagnetic energy. Moreover, any functional of the form $\int_\Ox \int_\Ov G(f) \, d\bv d\bx$ is a constant of motion.

Various types of numerical methods have been developed to compute the VM system.
 This includes the  popular Particle-in-cell (PIC) methods \cite{brackbill1982implicit, Jacobs_Hesthaven_2006, Hesthaven_09_divclean, Markidis20117037}.  In PIC methods, the macro-particles are advanced in a Lagrangian framework, while the field equations are solved on a mesh.  On the other hand, in recent years, there has been growing interest in computing kinetic equations in a deterministic fashion (i.e. the direct computation for the solutions to the Vlasov equations under Eulerian or semi-Lagrangian framework).  Deterministic solvers enjoy the advantage of producing highly accurate results without having any statistical noise. In the literature, semi-Lagrangian methods \cite{califano1998ksw, mangeney2002nsi, califano1965ikp, califano2001ffm}, spectral methods \cite{eliasson2003numerical}, finite difference method \cite{umeda2009two}, and Runge-Kutta DG (RKDG) methods \cite{cheng_vm} have been developed for VM systems.
 From computational point of the view, the main challenges for the deterministic simulations of the VM systems include: high dimensionality of the Vlasov equation, conservation of macroscopic quantities,  multiple temporal and spatial scales encountered in applications, and the desire to be able to work on unstructured meshes for real applications on complicated geometry in $\Ox$.

For most methods in the VM literature, the conservation of the total particle number is achieved, but  the conservation of total energy is not addressed, rather it was left to the accuracy of the scheme. For simulations in longer time ranges, the spurious energy created or annihilated by numerical methods could build up and lead to unphysical results, such as plasma self heating or cooling \cite{cohen1989performance}. This issue will be more prominent if we use under-resolved mesh or large time steps.
The total energy is a quantity that depends nonlinearly on the probability distribution function and the electromagnetic field, and it could  serve as a type of nonlinear stability bound for the scheme. Recently, several PIC methods have been proposed to conserve the total energy for VA, VM or Vlasov-Poisson (VP) system. In \cite{chen2011energy}, PIC for VA equations is developed; it is fully implicit, energy and charge conserving. In \cite{Markidis20117037}, PIC for VM system is developed, in which Maxwell's equation is solved on Yee's lattice \cite{yee1966numerical} and implicit midpoint method  is used as the time integrator. In \cite{Filbet_compare_2003, ayuso2012high},  finite difference and DG methods were proposed to conserve the total energy of VP systems. There is also abundant literature on energy-conserving Maxwell solvers. 
They include but are not limited to, the finite-difference time-domain (FDTD) method \cite{yee1966numerical, taflove2000computational}, finite-element time-domain \cite{rodrigue2001vector}, finite-volume time-domain methods \cite{piperno2002nondiffusive}, and discontinuous Galerkin time-domain (DGTD) methods \cite{fezoui2005convergence, piperno_2006, Chung2012}.

In this paper, we  generalize our energy-conserving methods for the VA system \cite{cheng_va} and develop energy-conserving DG methods for VM systems.   The schemes  in \cite{cheng_va} are proven to conserve the total particle number and the total energy on the fully discrete level for the VA system.  For the VM system, additional care needs to be taken to ensure such conservation properties, since the magnetic field is no longer negligible.
We aim to   address some of the common challenges for deterministic solvers. The issue of high dimensionality is treated by a new  splitting for the VM system so that the resulting equations are in reduced dimensions and still preserve energy conservation. We design   spatial discretizations by DG methods  with appropriate flux  to maintain energy conservation and still being able  to deal with filamentation. The symplectic integrators for Maxwell's equations  are carefully coupled with suitable time integrators for Vlasov equations to achieve fully discrete energy conservation. Implicit and explicit methods are designed under the same framework to deal with application problems with different stiffness.  The schemes designed have potential to be implemented on general unstructured mesh in $\Ox$.

Before we proceed, we would like to remark on a few assumptions and limitations  for our computation. As usual, we 
assume that $f(\bx, \bv, t)$ remains  compactly  supported in $\bv$,  given that it is initially so.    Whether or not the three-dimensional VM system is globally well-posed as a Cauchy problem is a major open problem.  The limited results of global existence without uniqueness of weak solutions and  well-posedness and regularity of solutions assuming either some symmetry or near  neutrality constitute the present  extent of  knowledge  \cite{GlasseyARMA86,GlasseyCMP87, GlasseyCMP88, DipernaCPAM89,GlasseyCMP97, GlasseyARMA98II,GlasseyARMA98I}.   In this paper, we will always take $\Ov$ to be finite and assume that $\Ov$ is taken large enough, so that the numerical solution $f_h \approx 0$ at $\partial \,\Ov$.  This can be achieved by  enlarging the velocity domain, and some related discussions can be found in \cite{cheng_vm}.
Another issue is related to the Gauss's law, i.e. the last two equations in \eqref{modelvm}. On the PDE level, those relations can be derived from the remaining part of the VM system;  therefore,  the numerical methods proposed in this paper are formulated for the VM system without those parts. We want to stress that even though in principle the initial satisfaction of these constraints is sufficient for their satisfaction for all time to certain accuracy, in certain circumstance one may need to consider explicitly  such divergence conditions  in order to produce  physically relevant numerical simulations \cite{Munz:2000, barth2006role}. In this paper, we do not attempt to address such issues. In particular, we will  present our numerical scheme in the general setting, and then discuss the details in 1D2V case by streaming Weibel instability.

The remaining part of the paper is  organized as follows:
in Sections \ref{sec:temporal} and \ref{sec:full}, the numerical schemes and their   properties are discussed.  In particular, Section \ref{sec:temporal} is devoted to the temporal discretizations, while in Section \ref{sec:full} the fully discrete methods are outlined. Section \ref{sec:numerical} includes the simulation results, and we conclude with a few remarks in Section \ref{sec:conclusion}.

\section{Numerical methods: temporal discretizations}
\label{sec:temporal}
\setcounter{equation}{0}
\setcounter{figure}{0}
\setcounter{table}{0}

In this section, we will describe the first main component of our schemes: energy-conserving temporal discretizations.
We leave the variables $(\bx, \bv)$ continuous in the discussions, and therefore, the time integrators introduced in this section can potentially be coupled with other spatial discretizations than those considered in Section \ref{sec:full}.

Before we discuss the details of our methods, we want to emphasize the relations of our methods with the symplectic integrators.
Symplectic integrators \cite{leimkuhler2005simulating, hairer2006geometric}
 for the   Hamiltonian systems  are known to possess as a conserved quantity, which is a Hamiltonian that is slightly perturbed from the original one.
Those methods are widely used   for the Maxwell's equation to preserve the electromagnetic energy. Some of the methods we proposed below are of this nature, while some others, e.g. those in Section \ref{sec:split} are motivated by and tailored to the specific structure of the VM system.

The outline of this section is as follows: we will first establish second-order explicit and implicit energy-conserving temporal discretizations in Section \ref{sec:secondorder}. Then to treat the fully implicit method more efficiently without inverting in the $(\bx, \bv)$ space, we propose an operator splitting in Section \ref{sec:split}. Finally, we will discuss how to improve the method beyond second order in Section \ref{sec:high}.

\subsection{Second order schemes}
\label{sec:secondorder}
In this subsection, we introduce  four types of methods for the coupled VM system, namely (1) explicit for Vlasov and Maxwell, (2) explicit for Vlasov and implicit for Maxwell, (3) implicit for Vlasov and explicit for Maxwell, (4) fully implicit schemes. Those four methods can potentially work for  VM equations in  various regimes when different types of stiffness occur. We will first define the methods and defer the rigorous   proof for   energy conservation to Theorem \ref{thm:energy}.

A prototype (1) scheme can be constructed by the leapfrog method for the Maxwell's equation and second order explicit Runge-Kutta method for Vlasov equation. To advance from $\{f^n, \bE^n, \bB^n\}$ to $\{f^{n+1}, \bE^{n+1}, \bB^{n+1}\}$, we use the scheme \eqref{scheme1}  and denote it as $\textnormal{\bf Scheme-1}(\Dt)$, i.e. 
$$(f^{n+1}, \bE^{n+1}, \bB^{n+1}) = \textnormal{\bf Scheme-1}(\Dt)(f^{n},\bE^{n},\bB^{n}).$$
\underline{$\textnormal{\bf Scheme-1}(\Dt)$}
\begin{subequations}
\label{scheme1}
\begin{align}
&\frac{f^{n+1/2}-f^n}{\Delta t/2} +  \bv \cdot \nabla_{\bx} f^n  +  (\bE^n + \bv \times \bB^n) \cdot \nabla_\bv f^n = 0~,\\
&\frac{\bB^{n+1/2}-\bB^n}{\Delta t/2}=-  \nabla_\bx \times \bE^n,\\
&\frac{\bE^{n+1}-\bE^n}{\Delta t}=  \nabla_\bx \times \bB^{n+1/2}-\bJ^{n+1/2}, \quad \textrm{where}\, \,\bJ^{n+1/2}= -\int f^{n+1/2} \bv d\bv\\
&\frac{\bB^{n+1}-\bB^\nh}{\Delta t/2}=-  \nabla_\bx \times \bE^{n+1},\\
&\frac{f^{n+1}-f^n}{\Delta t} +  \bv \cdot \nabla_{\bx} f^{n+1/2}  +  \left(\frac{1}{2}\left(\bE^n+\bE^{n+1}\right) + \bv \times \bB^{n+1/2}\right) \cdot \nabla_\bv f^{n+1/2} = 0.\
\end{align}
\end{subequations}

On the other hand, scheme of type (2) can be designed based on an implicit midpoint method for the Maxwell's equation. We denote the scheme 
\eqref{scheme2} by $\textnormal{\bf Scheme-2}(\Dt)$.
\\
\underline{$\textnormal{\bf Scheme-2}(\Dt)$}
\begin{subequations}
\label{scheme2}
\begin{align}
&\frac{f^{n+1/2}-f^n}{\Delta t/2} +  \bv \cdot \nabla_{\bx} f^n  +  (\bE^n + \bv \times \bB^n) \cdot \nabla_\bv f^n = 0~,\\
&\frac{\bB^{n+1}-\bB^n}{\Delta t}=-  \nabla_\bx \times \left(\frac{\bE^n+\bE^\na}{2} \right),\\
&\frac{\bE^{n+1}-\bE^n}{\Delta t}=  \nabla_\bx \times \left(\frac{\bB^n+\bB^\na}{2} \right)-\bJ^{n+1/2}, \quad \textrm{where}\, \,\bJ^{n+1/2}= -\int f^{n+1/2} \bv d\bv\\
&\frac{f^{n+1}-f^n}{\Delta t} +  \bv \cdot \nabla_{\bx} f^{n+1/2}  +  \left(\frac{1}{2}\left(\bE^n+\bE^{n+1}\right) + \frac{1}{2} \bv \times (\bB^n+\bB^{n+1})\right) \cdot \nabla_\bv f^{n+1/2} = 0.\
\end{align}
\end{subequations}
This method would work well for low frequency plasmas as the normalized  speed of light $c^\nu \rightarrow \infty$, where $\nu$ is  the characteristic speed of the electrons. This type of semi-implicit schemes are used quite often in PIC methods, where particles are evolved explicitly and field equations are solved implicitly, see for example \cite{brackbill1982implicit}. 

Scheme of type (3) can be formulated by using an implicit midpoint method for the Vlasov equation, and leap frog method for the Maxwell's equation. We denote the method \eqref{scheme3}
to be $\textnormal{\bf Scheme-3}(\Dt)$.
\\
{\bf \underline{$\textnormal{\bf Scheme-3}(\Dt)$}}
\begin{subequations}
\label{scheme3}
\begin{align}
&\frac{\bE^{n+1/2}-\bE^n}{\Delta t/2}= \nabla_\bx \times \bB^{n}-\bJ^{n}, \quad \textrm{where}\, \,\bJ^{n}= -\int f^{n} \bv d\bv\\
&\frac{\bB^{n+1}-\bB^n}{\Delta t}=-  \nabla_\bx \times \bE^\nh,\\
&\frac{f^{n+1}-f^n}{\Delta t} +  \bv \cdot \nabla_{\bx} \frac{ f^n+f^\na}{2}  +  \left(\bE^\nh + \bv \times \frac{\bB^n+\bB^{n+1}}{2}\right) \cdot \nabla_\bv \frac{ f^n+f^\na}{2} = 0~,\\
&\frac{\bE^{n+1}-\bE^\nh}{\Delta t/2}=  \nabla_\bx \times \bB^{n+1}-\bJ^{n+1}, \quad \textrm{where}\, \,\bJ^{n+1}= -\int f^{n+1} \bv d\bv.
\end{align}
\end{subequations}
This scheme should apply to the case when the Vlasov equation is stiff, while Maxwell's equation is not stiff. For simplicity, here we still consider our model equation \eqref{modelvm} in the nonrelativistic setting. We remark that implicit solves for the high dimensional Vlasov equation has to be implemented efficiently to make this scheme competitive. 

For plasma simulations, some applications incur stiffness in both Vlasov and Maxwell's equations. This includes the case of multi-species simulations, where the electron time scale is much faster than the ion time scale.  In those cases,   fully implicit methods are desirable. For our model of single species of nonrelativistic electrons \eqref{modelvm}, schemes of type (4) can be directly formulated by using implicit midpoint methods on the whole VM system, and we denote it to be $\textnormal{\bf Scheme-4}(\Dt)$.
\\
\underline{$\textnormal{\bf Scheme-4}(\Dt)$}
 \begin{subequations}
\label{scheme4}
\begin{align}
&\frac{\bB^{n+1}-\bB^n}{\Delta t}=-  \nabla_\bx \times \left(\frac{\bE^n+\bE^\na}{2} \right),\\
&\frac{\bE^{n+1}-\bE^n}{\Delta t}=  \nabla_\bx \times \left(\frac{\bB^n+\bB^\na}{2} \right)-\frac{\bJ^n+\bJ^\na}{2},  \\
&\frac{f^{n+1}-f^n}{\Delta t} +  \bv \cdot \nabla_{\bx} \frac{ f^n+f^\na}{2}   +  \left(\frac{1}{2}\left(\bE^n+\bE^{n+1}\right) + \frac{1}{2} \bv \times \left(\bB^n+\bB^{n+1}\right)\right) \cdot \nabla_\bv \frac{ f^n+f^\na}{2}  = 0.\
\end{align}
\end{subequations}

However the computation of this method is very demanding as it requires  inversion of a nonlinear high-dimensional coupled system. We will address this issue in detail under the splitting framework in the next subsection.

Through simple Taylor expansions, we can verify that the schemes above are all second order accurate in time.  In the next theorem, we will establish energy conservation for those methods. To simplify the discussion, we always assume periodic boundary conditions at $\Ox$ boundaries in this paper. For other boundary conditions, additional contributions from $\partial \Ox$ has to be considered, as is the case for the PDE itself.
\begin{thm}
\label{thm:energy}
The schemes introduced in this subsection preserve the discrete total energy $TE_n=TE_{n+1}$, where
$$2\ (TE_n)=\int_\Ox \int_\Ov f^n |\bv|^2 d\bv d\bx +\int_\Ox |\bE^n|^2 + |\bB^n|^2d\bx$$
in  $\textnormal{\bf Scheme-2}(\Dt)$ and $\textnormal{\bf Scheme-4}(\Dt)$, and
\begin{eqnarray*}
2\ (TE_n)&=&  \int_\Ox \int_\Ov f^n |\bv|^2 d\bv d\bx + \int_\Ox (|\bE^n|^2+\bB^{n-1/2}\cdot\bB^\nh) d\bx \\
&=&\int_\Ox \int_\Ov f^n |\bv|^2 d\bv d\bx +\int_\Ox |\bE^n|^2 + |\bB^n|^2 d\bx-\frac{\Delta t^2}{4} \int_\Ox |\nabla_\bx \times \bE^n|^2 d\bx
\end{eqnarray*}
in  $\textnormal{\bf Scheme-1}(\Dt)$, and
\begin{eqnarray*}
2\ (TE_n)&=& \int_\Ox \int_\Ov f^n |\bv|^2 d\bv d\bx + \int_\Ox (|\bB^n|^2+\bE^{n-1/2}\cdot\bE^\nh) d\bx. \\
&=&\int_\Ox \int_\Ov f^n |\bv|^2 d\bv d\bx +\int_\Ox |\bE^n|^2 + |\bB^n|^2d\bx-\frac{\Delta t^2}{4} \int_\Ox |\nabla_\bx \times \bB^n-\bJ^n|^2 d\bx
\end{eqnarray*}
in  $\textnormal{\bf Scheme-3}(\Dt)$.
\end{thm}
\emph{Proof.} We will prove the theorem for  $\textnormal{\bf Scheme-1}(\Dt)$. The other three proofs are similar and are omitted. Using the definition of  $\textnormal{\bf Scheme-1}(\Dt)$ in (\ref{scheme1}), 
\begin{eqnarray*}
&&\hspace{-6mm}\int_\Ox \int_\Ov \frac{f^{n+1}-f^n}{\Delta t} |\bv|^2 d\bv d\bx \\
&&\hspace{-6mm}=-\int_\Ox \int_\Ov \bv \cdot \nabla_{\bx} f^{n+1/2}  |\bv|^2 d\bv d\bx   -  \int_\Ox \int_\Ov  \left( \frac{1}{2}(\bE^n+\bE^{n+1})+ \bv \times \bB^{n+1/2}\right)  \cdot \nabla_\bv f^{n+1/2}  |\bv|^2 d\bv d\bx  \\
&&\hspace{-6mm} =-\int_\Ov  |\bv|^2 \bv \cdot (\int_\Ox  \nabla_{\bx} f^{n+1/2}  d\bx)\, d\bv   +  \int_\Ox \int_\Ov  (\bE^n+\bE^{n+1}+2 \, \bv \times \bB^{n+1/2}) \cdot f^{n+1/2}  \bv d\bv d\bx \\
&&\hspace{-6mm}=\int_\Ox   (\bE^n+\bE^{n+1}) \cdot\bJ^{n+1/2}   d\bx
\end{eqnarray*}
where in the second equality, we used the integration by parts, and in the  third equality we   employed the periodic boundary conditions in $\Ox$ domain. On the other hand, 
\begin{eqnarray*}
&&\hspace{-6mm}\int_\Ox  \frac{\bE^{n+1}-\bE^n}{\Delta t} \cdot (\bE^{n+1}+\bE^n)d\bx=-\int_\Ox  (\nabla_\bx \times \bB^{n+1/2}-\bJ^{n+1/2})\cdot (\bE^n+\bE^{n+1})    d\bx,
\end{eqnarray*}
and
\begin{eqnarray*}
&&\hspace{-6mm}\int_\Ox  \frac{\bB^{n+3/2}-\bB^{n-1/2}}{\Delta t} \cdot \bB^{n+1/2}d\bx=-\int_\Ox  \nabla_\bx \times (\bE^n+\bE^{n+1})\cdot \bB^{n+1/2}   d\bx.
\end{eqnarray*}
Therefore,
\begin{eqnarray*}
&&\hspace{-6mm}\int_\Ox \int_\Ov f^n |\bv|^2 d\bv d\bx \,+ \int_\Ox (|\bE^n|^2+\bB^{n-1/2}\cdot\bB^\nh) d\bx \\
&&\hspace{-6mm}=\int_\Ox \int_\Ov f^\na |\bv|^2 d\bv d\bx + \int_\Ox (|\bE^\na|^2+\bB^{n+1/2}\cdot\bB^{n+3/2}) d\bx. \hspace{45mm}  \Box
\end{eqnarray*}
  
 From this theorem, we can see that  {\bf Scheme-2} and {\bf Scheme-4} exactly preserve the total energy, while {\bf Scheme-1} and  {\bf Scheme-3} achieve near conservation of the total energy. The numerical energies from {\bf Scheme-1} and  {\bf Scheme-3} are second order modified version of the original total energy. This ensures that over the long run, the numerical energy will not deviate much from its actual value.
\subsection{An energy-conserving operator splitting for the VM system}
\label{sec:split}

Following the lines of our previous work \cite{cheng_va}, here we propose an operator splitting of the VM system to efficiently compute for the fully implicit   {\bf Scheme-4} . This splitting operator is specifically tailored to the VM system such that each split equation can still preserve the total energy.
In particular, for the model VM equation,   the operator splitting is done as follows:
\[
  \textrm{ (a)} \left\{
  \begin{array}{l}
\partial_t f  +  \bv \cdot \nabla_{\bx} f   = 0~, \\
\partial_t \bE = 0, \\
\partial_t \bB = 0,
  \end{array} \right.
  \textrm{(b)} \left\{
  \begin{array}{l}
\partial_t f  +   \bE  \cdot \nabla_\bv f   = 0~, \\
\partial_t \bE = -\bJ, \\
\partial_t \bB = 0,
  \end{array} \right.
  \textrm{(c)}\left\{
  \begin{array}{l}
\partial_t f  +  ( \bv \times \bB)  \cdot \nabla_\bv f   = 0~, \\
\displaystyle\partial_t \bE = \nabla_\bx \times \bB, \\
\partial_t \bB =  -  \nabla_\bx \times \bE.
  \end{array} \right.
\]
We can verify that each of the three equations is energy-conserving,
$$\frac{d}{dt} (\int_\Ox \int_\Ov f |\bv|^2 d\bv d\bx+ \int_\Ox (|\bE|^2+|\bB|^2) d\bx ) = 0.$$ In particular,
\[
  \textrm{ (a)} \left\{
  \begin{array}{l}
\displaystyle\frac{d}{dt} \int_\Ox \int_\Ov f |\bv|^2 d\bv d\bx = 0~, \\[4mm]
\displaystyle\frac{d}{dt}  \int_\Ox |\bE|^2 d\bx = 0, \\[4mm]
\displaystyle\frac{d}{dt}  \int_\Ox |\bB|^2 d\bx= 0,
  \end{array} \right.
  \textrm{(b)} \left\{
  \begin{array}{l}
\displaystyle\frac{d}{dt} (\int_\Ox \int_\Ov f |\bv|^2 d\bv d\bx+ \int |\bE|^2 d\bx ) = 0~, \\[4mm]
\displaystyle\frac{d}{dt}  \int_\Ox |\bB|^2 d\bx= 0,
  \end{array} \right.\]
  \[
    \textrm{(c)} \left\{
  \begin{array}{l}
\displaystyle\frac{d}{dt} \int_\Ox \int_\Ov f |\bv|^2 d\bv d\bx = 0~, \\[4mm]
\displaystyle\frac{d}{dt}  \int_\Ox(|\bE|^2+|\bB|^2) d\bx = 0.
  \end{array} \right.
\]
We can see that equation (a) contains the free streaming operator. In this equation, the electromagnetic fields are unchanged, and the kinetic energy is conserved. Equation (b) contains the interchange of kinetic and electric energy, while the magnetic field is unchanged. Equation (c)  is comprised of the Maxwell's equation and the rotation of $f$ under the magnetic field. In this process, the kinetic energy and the electromagnetic energy are conserved respectively. We also notice that combing equation (a) and (b) yields the VA system.

Using this splitting, we only have to solve each individual equation in an energy-conserving manner, and then carefully combine them using a splitting method of desired order \cite{mclachlan2002splitting, strang1968construction, yoshida1990construction, forest1990fourth} that can keep the conservation of energy for the split equations. Moreover, each of the equations is now essentially decoupled  and in lower dimensions,  therefore we can solve it more efficiently. Now let's discuss the details of the scheme for each split equation.

As for equation (a), we can use any implicit or explicit Runge-Kutta methods to solve it, and they all conserve the kinetic energy. To see this, consider the forward Euler
$$\frac{f^{n+1}-f^n}{\Delta t} +  \bv \cdot \nabla_{\bx} f^n=0,$$
or  backward Euler method
$$\frac{f^{n+1}-f^n}{\Delta t} +  \bv \cdot \nabla_{\bx} f^\na=0,$$
A simple check yields
$\int_\Ox \int_\Ov  f^n |\bv|^2 \,d\bv d\bx=\int_\Ox \int_\Ov  f^\na |\bv|^2 \,d\bv d\bx$. (Note that here we have abused the notation, and use superscript $n$, $n+1$ to denote the sub steps in computing equation (a), not the whole time step to compute the VM system). Therefore, we can pick a suitable Runge-Kutta method with desired order and property for this step. To be second order, one can for example use the implicit midpoint method,
\beq
\label{schemea}
\frac{f^{n+1}-f^n}{\Delta t} +  \bv \cdot \nabla_{\bx} \frac{f^n+f^\na}{2}=0.
\eeq

Equation (b) contains the main coupling effect of the Vlasov and Maxwell's equation, and has to be computed carefully to balance the kinetic and electric energies. We can use the methods studied in  Section \ref{sec:secondorder} to compute this equation. (We only need to include the corresponding terms as those appeared in equation (b)).  The resulting scheme will naturally preserve a discrete form of  the sum of kinetic and electric energies.
In particular, we will use
\begin{subequations}
\label{scheme4s}
\begin{align}
&\frac{f^{n+1}-f^n}{\Delta t}    +  \frac{1}{2}(\bE^n+\bE^{n+1})  \cdot \nabla_\bv \frac{f^n+f^{n+1}}{2} = 0~,\label{scheme2s:1}\\
&\frac{\bE^{n+1}-\bE^n}{\Delta t}=- \frac{1}{2} (\bJ^{n}+\bJ^{n+1}) \label{scheme2s:2}
\end{align}
\end{subequations}

Similarly, for equation (c), we will use
\begin{subequations}
\label{schemec}
\begin{align}
&\label{schemec:e1}\frac{\bB^{n+1}-\bB^n}{\Delta t}=-  \nabla_\bx \times \left(\frac{\bE^n+\bE^\na}{2} \right),\\
&\label{schemec:e2}\frac{\bE^{n+1}-\bE^n}{\Delta t}=  \nabla_\bx \times \left(\frac{\bB^n+\bB^\na}{2} \right),  \\
&\label{schemec:e3}\frac{f^{n+1}-f^n}{\Delta t}  + \frac{1}{2} \bv \times (\bB^n+\bB^{n+1}) \cdot \nabla_\bv \frac{ f^n+f^\na}{2}  = 0~,\
\end{align}
\end{subequations}
Notice that $\bB^{n+1}$ can  be computed by solving \eqref{schemec:e1} and \eqref{schemec:e2}.  Then we can plug $\bB^{n+1}$ into \eqref{schemec:e3} to solve  $f^{n+1}$.

Now let  $\textnormal{\bf Scheme-a}(\Dt)$  denote second order  schemes for equation (a) ,  $\textnormal{\bf Scheme-b}(\Dt)$  denote second order schemes for equation (b), 
and $\textnormal{\bf Scheme-c}(\Dt)$ denote second order schemes for equation (c), then   by Strang splitting 
\begin{eqnarray*}
\textnormal{\bf Scheme-a}(\Dt/2)\textnormal{\bf Scheme-b}(\Dt/2)\textnormal{\bf Scheme-c}(\Dt)\textnormal{\bf Scheme-b}(\Dt/2)\textnormal{\bf Scheme-a}(\Dt/2),
\end{eqnarray*}
is a second order scheme for the original VM system. 

\begin{thm}
\label{lem:energy_split}
If $\textnormal{\bf Scheme-a}$, $\textnormal{\bf Scheme-b}$, $\textnormal{\bf Scheme-c}$ are defined by \eqref{schemea},    \eqref{scheme4s},   \eqref{schemec}, respectively, then  
\begin{eqnarray*}
&&\textnormal{\bf Scheme-5}(\Dt):=\\
&&\textnormal{\bf Scheme-a}(\Dt/2)\textnormal{\bf Scheme-b}(\Dt/2)\textnormal{\bf Scheme-c}(\Dt)\textnormal{\bf Scheme-b}(\Dt/2)\textnormal{\bf Scheme-a}(\Dt/2),
\end{eqnarray*}
 preserves the discrete total energy $TE_n=TE_{n+1}$, where
$$2\ (TE_n)=\int_\Ox \int_\Ov f^n |\bv|^2 d\bv d\bx +\int_\Ox |\bE^n|^2 + |\bB^n|^2d\bx.$$\end{thm}
The proof for this theorem is straightforward by  discussion in this subsection and is omitted.

\subsection{Generalizations to higher order}
\label{sec:high}
Similar to the discussion in \cite{cheng_va}, we can generalize the  symmetric-in-time second order schemes to higher order based on previous works \cite{yoshida1990construction, forest1990fourth, de1992easily, sanz1991order}. In particular, we discuss the fourth order methods in details below.  The generalization to even higher order follows the same idea.

Let $\beta_1, \beta_2, \beta_3$ satisfy
$$\beta_1+\beta_2+\beta_3=1,  \qquad \beta^3_1+\beta^3_2+\beta^3_3=1, \qquad \beta_1=\beta_3,$$
 and we get $\beta_1=\beta_3=(2+2^{1/3}+2^{-1/3})/3\approx1.3512$, $\beta_2=1-2\beta_1\approx-1.7024$.

We define
$$\textnormal{\bf Scheme-3F}(\Delta t)=\textnormal{\bf Scheme-3}(\beta_1 \Delta t) \textnormal{\bf Scheme-3}(\beta_2 \Delta t) \textnormal{\bf Scheme-3}(\beta_3 \Delta t);$$
$$\textnormal{\bf Scheme-4F}(\Delta t)=\textnormal{\bf Scheme-4}(\beta_1 \Delta t) \textnormal{\bf Scheme-4}(\beta_2 \Delta t) \textnormal{\bf Scheme-4}(\beta_3 \Delta t);$$
$$\textnormal{\bf Scheme-5F}(\Delta t)=\textnormal{\bf Scheme-5}(\beta_1 \Delta t) \textnormal{\bf Scheme-5}(\beta_2 \Delta t) \textnormal{\bf Scheme-5}(\beta_3 \Delta t).$$
Then  $\textnormal{\bf Scheme-3F}(\Delta t), \textnormal{\bf Scheme-4F}(\Delta t),
 \textnormal{\bf Scheme-5F}(\Delta t)$ are all fourth order. On the other hand, this procedure won't work for $\textnormal{\bf Scheme-1}$ and $\textnormal{\bf Scheme-2}$, because they are not symmetric in time.
 
\begin{thm}
\label{thm:energy2}
{\bf Scheme-4F}  and {\bf Scheme-5F} 
 preserve the discrete total energy $TE_n=TE_{n+1}$, where
$$2\,  (TE_n)=\int_\Ox \int_\Ov f^n |\bv|^2 d\bv d\bx +\int_\Ox |\bE^n|^2 +\bB^n|^2 d\bx.$$
\end{thm}
The proof is straightforward by the properties of the second order methods {\bf Scheme-4} and {\bf Scheme-5} and is omitted.

However, it is  challenging to obtain the explicit form of the modified total energy for $\textnormal{\bf Scheme-3F}(\Delta t)$. 
Finally, we remark that for those negative time steps caused by $\beta_2<0$,  special cares need to be taken. For example, if numerical dissipation is added for the positive time steps, we need to make sure to add anti-dissipation for the negative time steps. This means that for the schemes described in Section \ref{sec:full}, the upwind fluxes for Vlasov and Maxwell's equations need to be changed to downwind fluxes for the negative time steps.
%
%

\section{Numerical methods: fully discrete schemes}
\setcounter{equation}{0}
\setcounter{figure}{0}
\setcounter{table}{0}
\label{sec:full}

In this section, we will discuss the spatial discretizations and formulate the fully discrete schemes. In particular, we consider two approaches: one being the unsplit schemes, the other being the split implicit schemes.

In this paper, 
we choose a discontinuous Galerkin (DG) discretization of the (x,v) variables since such methods can be made
      arbitrarily high-order while retaining exact mass conservation and, in conjunction
      with the previously described time discretizations, total energy conservation.
The DG method \cite{Cockburn_2000_history, Cockburn_2001_RK_DG} is a class
of finite element methods using discontinuous piecewise
polynomial space for the numerical solution and the test functions, and they are originally designed to solve conservation laws.  In recent years, high order DG schemes have shown their attractive properties in accuracy and conservation to simulate the VP system  \cite{rossmanith2011positivity, qiu2011positivity, heath2012discontinuous, Heath_thesis, Ayuso2009, Ayuso2010, cheng_vp}.
 RKDG schemes have been designed to solve  VM \cite{cheng_vm} systems, and semi-discrete total energy conservation have been established  \cite{cheng_vm}. In the discussions below, we will prove fully discrete conservation properties for our proposed methods. The DG methods with unsplit schemes are related to the methods in \cite{cheng_vm}. However,  the standard   Runge-Kutta methods are replaced by temporal schemes discussed in the previous section. For the split schemes, we consider the methods that are implemented based on the Gauss quadrature points \cite{cheng_va}.

\subsection{Notations}
In this subsection, we will introduce the mesh and polynomial space under consideration.
Let $\mT_h^x=\{\Kx\}$ and $\mT_h^v=\{\Kv\}$ be  partitions of $\Ox$ and $\Ov$, respectively,  with $\Kx$ and $\Kv$ being (rotated) Cartesian elements or 
simplices;  then $\mT_h=\{K: K=\Kx\times\Kv, \forall \Kx\in\mT_h^x, \forall \Kv\in\mT_h^v\}$ defines a partition of $\Omega$. Let $\mE_x$ be the set of the edges of $\mT_h^x$ and  $\mE_v$ be the set of the edges of $\mT_h^v$;  then the edges of $\mT_h$ will be $\mE=\{\Kx\times e_v: \forall\Kx\in\mT_h^x, \forall e_v\in\mE_v\}\cup \{e_x\times\Kv: \forall e_x\in\mE_x, \forall \Kv\in\mT_h^v\}$. Here we take into account the periodic boundary condition in the $\bx$-direction when defining $\mE_x$ and $\mE$. Furthermore, $\mE_v=\mE_v^i\cup\mE_v^b$ with $\mE_v^i$ and $\mE_v^b$ being the set of interior and boundary edges of $\mT_h^v$, respectively.  

We will make use of  the following discrete spaces
\begin{subequations}
\begin{align}
\mG_h^{k}&=\left\{g\in L^2(\Omega): g|_{K=\Kx\times\Kv}\in P^k(\Kx\times\Kv), \forall \Kx\in\mT_h^x, \forall \Kv\in\mT_h^v \right\}~,\label{eq:sp:f}\\
\mS_h^{k}&=\left\{g\in L^2(\Omega): g|_{K=\Kx\times\Kv}\in Q^k(\Kx)\times Q^k(\Kv), \forall \Kx\in\mT_h^x, \forall \Kv\in\mT_h^v \right\}~,\\
\mU_h^r&=\left\{\bU\in [L^2(\Omega_x)]^{d_x}: \bU|_\Kx\in [P^r(\Kx)]^{d_x}, \forall \Kx\in\mT_h^x \right\}~,\\
\mW_h^r&=\left\{w\in L^2(\Omega_x): w|_\Kx\in Q^r(\Kx), \forall \Kx\in\mT_h^x \right\}~,\\
\mZ_h^r&=\left\{z\in L^2(\Omega_v): w|_\Kv\in Q^r(\Kv), \forall \Kv\in\mT_h^v \right\}~,
\end{align}
\end{subequations}
where $P^k(D)$ denotes the set of polynomials of  total degree at most $k$ on $D$,  and $Q^k(D)$ denotes the set of polynomials of degree at most $k$ in each variable on $D$. Here  $k$ and $r$ are non-negative integers.
The discussion about those spaces for Vlasov equations can be found in \cite{cheng_vp, cheng_vm}.

For piecewise  functions defined with respect to $\mT_h^x$ or $\mT_h^v$, we further introduce the jumps and averages as follows. 
For any edge $e=\{K_.^+\cap K_.^-\}\in\mE_.$, with $\bn_.^\pm$ as the outward unit normal to $\partial K_.^\pm$,
$g^\pm=g|_{K_.^\pm}$, and $\bU^\pm=\bU|_{K_.^\pm}$, the jumps across $e$ are defined  as
\begin{equation*}
[g]_.={g^+}{\bn_.^+}+{g^-}{\bn_.^-},\qquad [\bU]_.={\bU^+}\cdot{\bn_.^+}+{\bU^-}\cdot{\bn_.^-},\quad [\bU]_\tau={\bU^+}\times{\bn_x^+}+{\bU^-}\times{\bn_x^-}
\end{equation*}
and the averages are
\begin{equation*}
\{g\}_{.}=\frac{1}{2}({g^+}+{g^-}),\qquad \{\bU\}_.=\frac{1}{2}({\bU^+}+{\bU^-}),
\end{equation*}
where $.$ are used to denote $x$ or $v$.

\subsection{Unsplit schemes and their properties}
\label{sec:unsplit}
In this subsection, we will describe the DG methods for the unsplit schemes {\bf Scheme-1}, {\bf Scheme-2}, {\bf Scheme-3}, {\bf Scheme-4}, {\bf Scheme-3F}, {\bf Scheme-4F} and discuss their properties.  For example, the scheme with {\bf Scheme-2}$(\Delta t)$  is formulated as follows: we look for $ f_h^{n+1/2},f_h^{n+1}\in\mG_h^k$, $\bB_h^{n+1/2}$, $\bB_h^{n+1}$, $\bE_h^{n+1} \in \mU_h^r$, 
such that for any $\psi_1, \psi_2 \in\mG_h^k$, $\bU, \bV\in \mU_h^r$,
\begin{subequations}
\label{dscheme1}
\begin{align}
&\int_K \frac{f_h^{n+1/2}-f_h^n}{\Delta t/2}  \psi_1 \,d\bv d\bx
- \int_K f_h^n\bv\cdot\nabla_\bx \psi_1 \,d\bv d\bx
- \int_K f_h^n (\bE_h^n +\bv \times \bB_h^n)\cdot\nabla_\bv \psi_1 \;\,d\bv d\bx\notag\\[2mm]
&\quad+ \int_{\Kv}\int_{\df\Kx} \widehat{f_h^n \bv\cdot \bn_x} \psi_1 ds_x d\bv 
+ \int_{\Kx} \int_{\df\Kv} \reallywidehat{(f_h^n\left( \bE_h^n  +\bv \times \bB_h^n\right)\cdot \bn_\bv)}{\psi_1} ds_v d\bx=0~,\label{dscheme1:1}\\[2mm]
&\int_{\Kx}\frac{\bB_h^{n+1}-\bB_h^{n}}{\Delta t} \cdot \bU\, d\bx=-\int_{\Kx}\frac{1}{2}\left(\bE_h^{n+1}+\bE_h^n\right)\cdot \nabla_\bx\times\bU d\bx 
-\int_{\df\Kx}\reallywidehat{\bn_x\times\frac{1}{2}\left(\bE_h^{n+1}+\bE_h^n\right)}\cdot \bU d\bx\label{dscheme1:2}\\[2mm]
&\int_\Kx\frac{\bE_h^{n+1}-\bE_h^n}{\Delta t}\cdot \bV d \bx=\int_{\Kx}\frac{1}{2}\left(\bB_h^{n+1}+\bB_h^n\right)\cdot \nabla_\bx\times\bV d\bx 
 -\int_\Kx \bJ_h^{n+1/2}\cdot \bV d \bx\notag\\[2mm]
 &\quad + \int_{\df\Kx}\reallywidehat{\bn_x\times\frac{1}{2}\left(\bB_h^{n+1}+\bB_h^n\right)}\cdot \bV d\bx, \label{dscheme1:3}\\[2mm]
& \int_K \frac{f_h^{n+1}-f_h^n}{\Delta t}  \psi_2 \,d\bv d\bx
- \int_K f_h^{n+1/2}\bv\cdot\nabla_\bx \psi_2 \,d\bv d\bx
+ \int_{\Kv}\int_{\df\Kx} \reallywidehat{f_h^{n+1/2} \bv\cdot \bn_x} \psi_2 ds_x d\bv \notag\\[2mm]
&\quad -  \int_K f_h^{n+1/2} \left(\frac{1}{2}(\bE_h^n+\bE_h^{n+1}) +\bv\times \frac{1}{2}\left(\bB_h^{n+1}+\bB_h^n\right)\right)\cdot\nabla_\bv \psi_2 \,d\bv d\bx\notag\\[2mm]
&\quad
+ \int_{\Kx} \int_{\df\Kv} \reallywidehat{\left(f_h^{n+1/2} \left(\frac{1}{2}(\bE_h^n+\bE_h^{n+1}) +\bv\times \frac{1}{2}\left(\bB_h^{n+1}+\bB_h^n\right)\right) \cdot \bn_\bv\right)} \psi_2 ds_v d\bx=0,\label{dscheme1:4}
\end{align}
\end{subequations}
with
$$ \bJ_h^{n+1/2}= \int_\Ov f_h^{n+1/2} \bv d\bv~.$$
Here $\bn_x$ and $\bn_v$ are outward unit normals of $\df\Kx$ and $\df\Kv$, respectively. All  `hat'  functions are numerical fluxes. For the Vlasov equation, the  fluxes in \eqref{dscheme1:1} are taken to be the central flux
\begin{align}
&\reallywidehat{f_h^n \bv\cdot \bn_x}= \{f_h^n\bv\}_x \cdot\bn_x,
\quad \reallywidehat{f_h^n ( \bE_h^n  +\bv \times \bB_h^n)\cdot \bn_\bv}= \{f_h^n ( \bE_h^n  +\bv \times \bB_h^n)\}_v \cdot\bn_x, \label{central:1}
\end{align}
or the upwind flux
\begin{subequations}
\begin{align}
\reallywidehat{f_h^n \bv\cdot \bn_x}& =\left(\{f_h^n\bv\}_x+\frac{|\bv\cdot\bn_x|}{2}[f_h^n]_x\right)\cdot\bn_x~,\label{upwind:1}\\
\reallywidehat{f_h^n ( \bE_h^n  +\bv \times \bB_h^n)\cdot \bn_\bv}&
=\left(\{f_h^n ( \bE_h^n  +\bv \times \bB_h^n)\}_v+\frac{|( \bE_h^n  +\bv \times \bB_h^n)\cdot\bn_\bv|}{2}[f_h^n]_\bv\right)\cdot\bn_v~.\label{upwind:2}
\end{align}
\end{subequations}
The flux terms in \eqref{dscheme1:4} are defined similarly. For the Maxwell's equation, we consider the central flux  
\begin{equation}
\label{central:2}
\begin{split}
&\reallywidehat{\bn_x\times \frac{1}{2}(\bE_h^{n+1}+\bE_h^n)}  =\bn_x\times \{ \frac{1}{2}(\bE_h^{n+1}+\bE_h^n)\}_x,\\
& \reallywidehat{\bn_x\times \frac{1}{2}(\bB_h^{n+1}+ \bB_h^n)} =\bn_x\times \{ \frac{1}{2}(\bB_h^{n+1}+ \bB_h^n)\}_x;
\end{split}
\end{equation}
 and the alternating flux
\begin{equation}
\label{alternating:1}
 \begin{split}
&\reallywidehat{\bn_x\times \frac{1}{2}(\bE_h^{n+1}+\bE_h^n)} =\bn_x\times \frac{1}{2}(\bE_h^{n+1,+}+\bE_h^{n,+}),\;\;\\
 & \reallywidehat{ \bn_x\times \frac{1}{2}(\bB_h^{n+1}+ \bB_h^n)} = \bn_x\times\frac{1}{2}(\bB_h^{n+1,-}+ \bB_h^{n,-});
 \end{split}
 \end{equation}
 or
\begin{equation}
\label{alternating:2}
\begin{split}
&\reallywidehat{\bn_x\times \frac{1}{2}(\bE_h^{n+1}+\bE_h^n)} =  \bn_x\times \frac{1}{2}(\bE_h^{n+1,-}+\bE_h^{n,-}),\;\;\\
&\reallywidehat{\bn_x\times \frac{1}{2}(\bB_h^{n+1}+ \bB_h^n)} = \bn_x\times \frac{1}{2}(\bB_h^{n+1,+}+ \bB_h^{n,+}).
\end{split}
\end{equation}

The fully discrete schemes with  {\bf Scheme-1}, {\bf Scheme-3}, {\bf Scheme-4},  {\bf Scheme-3F}, {\bf Scheme-4F} as time discretizations can be defined similarly, i.e. to use
DG discretization to approximate the derivatives of the $f$ in $\bx,\bv$, and the derivatives of $\bE, \bB$ in $\bx$. In particular, for negative time steps in the fourth order schemes {\bf Scheme-3F}, {\bf Scheme-4F}, we will use the downwind flux
\begin{subequations}
\begin{align}
\reallywidehat{f_h^n \bv\cdot \bn_x}& =\left(\{f_h^n\bv\}_x-\frac{|\bv\cdot\bn_x|}{2}[f_h^n]_x\right)\cdot\bn_x~,\label{downwind:1}\\
\reallywidehat{f_h^n ( \bE_h^n  +\bv \times \bB_h^n)\cdot \bn_\bv}&
=\left(\{f_h^n ( \bE_h^n  +\bv \times \bB_h^n)\}_v-\frac{|( \bE_h^n  +\bv \times \bB_h^n)\cdot\bn_\bv|}{2}[f_h^n]_\bv\right)\cdot\bn_v~,\label{downwind:2}
\end{align}
\end{subequations}
in the corresponding schemes.
To save space, we do not include the detailed descriptions of those methods here. 

The flux choices are crucial for the accuracy, stability and conservation properties of the methods. As shown in \cite{cheng_va}, the central flux for the Vlasov equations causes lack of numerical dissipation. When filamentation occurs, the numerical schemes will produce spurious oscillation, jeopardizing the quality of the solution. On the other hand, the flux choices for Maxwell's equation are especially important for energy conservation. In particular, for the semi-discrete schemes with $t$ continuous, the central and alternating fluxes preserve the the total energy, while the upwind flux causes energy dissipation \cite{cheng_vm}. Therefore, in this paper, we do not consider the upwind flux for the Maxwell's equation.
\medskip

Next, we will establish conservation properties of the fully discrete methods.

\begin{thm}[{Total particle number conservation}]
\label{thm:mass}
The scheme \eqref{dscheme1} preserves the total particle number of the system, i.e.
$$\int_\Ox \int_\Ov f_h^\na d\bv d\bx =\int_\Ox \int_\Ov f_h^n  d\bv d\bx .$$
This also holds for DG methods with time integrators {\bf Scheme-1}$(\Delta t)$, {\bf Scheme-3}$(\Delta t)$ and {\bf Scheme-4}$(\Delta t)$, {\bf Scheme-3F}$(\Delta t)$ and {\bf Scheme-4F}$(\Delta t)$.
\end{thm}
\emph{Proof.}  Let $\psi_2=1$ in \eqref{dscheme1:4}, and sum over all elements $K$, and we obtain the conservation property for  {\bf Scheme-2}$(\Delta t)$.  The proof for {\bf Scheme-1}$(\Delta t)$, {\bf Scheme-3}$(\Delta t)$, {\bf Scheme-4}$(\Delta t)$,  {\bf Scheme-3F}$(\Delta t)$ and {\bf Scheme-4F}$(\Delta t)$ is similar and is thus omitted. $\Box$
\begin{thm}[{Total energy conservation}]
\label{thm:fullyenergy2}
If $k \geq 2, r\geq 0$, the scheme \eqref{dscheme1} with either the upwind numerical flux \eqref{upwind:1}-\eqref{upwind:2} or the central
numerical flux \eqref{central:1} for the Vlasov equation, and either the central numerical flux of \eqref{central:2}  or the alternating numerical flux 
of \eqref{alternating:1} or \eqref{alternating:2} for the Maxwell's equation preserves the discrete total energy $TE_n=TE_{n+1}$, where
$$2(TE_n)=\int_\Ox \int_\Ov f_h^n |\bv|^2 d\bv d\bx +\int_\Ox |\bE_h^n|^2 + |\bB_h^n|^2d\bx.$$
This also holds for DG methods with time integrator $\textnormal{\bf Scheme-4}(\Delta t)$ and $\textnormal{\bf Scheme-4F}(\Delta t)$.
For DG methods with time integrator {\bf Scheme-1}$(\Delta t)$, the numerical energy  defined as
\begin{eqnarray*}
2(TE_n)&=&  \int_\Ox \int_\Ov f_h^n |\bv|^2 d\bv d\bx + \int_\Ox (|\bE_h^n|^2+\bB_h^{n-1/2}\cdot\bB_h^\nh) d\bx 
\end{eqnarray*}
is also preserved. The same holds   for DG methods with time integrator $\textnormal{\bf Scheme-3}(\Delta t)$ with the numerical  energy defined by
\begin{eqnarray*}
2(TE_n)&=& \int_\Ox \int_\Ov f_h^n |\bv|^2 d\bv d\bx +\int_\Ox (|\bB_h^n|^2+\bE_h^{n-1/2}\cdot\bE_h^\nh) d\bx. \\
\end{eqnarray*}
\end{thm}
\emph{Proof.}  We will only show the proof for $\textnormal{\bf Scheme-2}(\Delta t)$.
The proof for $\textnormal{\bf Scheme-1}(\Delta t)$, $\textnormal{\bf Scheme-3}(\Delta t)$,  {\bf Scheme-4}$(\Delta t)$ and $\textnormal{\bf Scheme-4F}(\Delta t)$ is similar. 
Let $\psi_2=|\bv|^2$ in  \eqref{dscheme1:4}. Note that  $|\bv|^2 \in\mG_h^k$ if $k\geq 2$ and it is continuous. Moreover, $\nabla_\bx \psi_2=0$, $\nabla_\bv \psi_2= 2 \bv$. Sum up over all elements $K$, we get 
\begin{eqnarray*}
\hspace{-3mm}\int_\Ox \int_\Ov \frac{f_h^{n+1}-f_h^n}{\Delta t} |\bv|^2 d\bv d\bx
 =  \int_\Ox \int_\Ov  (\bE_h^n+\bE_h^{n+1}) \cdot\left(f_h^{n+1/2}  \bv \right)d\bv d\bx
=\int_\Ox   (\bE_h^n+\bE_h^{n+1}) \cdot\bJ_h^{n+1/2}   d\bx.
\end{eqnarray*}
Denote $$\overline{\bB}=\displaystyle\frac{1}{2} (\bB_h^{n+1}+\bB_h^{n}),\quad \overline{\bE}= \displaystyle\frac{1}{2} (\bE_h^{n+1}+\bE_h^{n})$$and
$$\widetilde{\overline{\bB}}=\{\overline{\bB}\}_x,\; \widetilde{\overline{\bE}}=\{\overline{\bE}\}_x;\quad\quad
\widetilde{\overline{\bB}}=\overline{\bB}^-,\; \widetilde{\overline{\bE}}= \overline{\bE}^+;\quad\quad
\widetilde{\overline{\bB}}=\overline{\bB}^+,\; \widetilde{\overline{\bE}}= \overline{\bE}^-;$$ to be associated with the central flux \eqref{central:2}, and the alternating fluxes  \eqref{alternating:1}, \eqref{alternating:2}, respectively.
 Take $\bU = \bB_h^{n+1}+\bB_h^{n}\in\mU_h^r$ in \eqref{dscheme1:2} and $\bV= \bE_h^{n+1}+\bE_h^{n}\in\mU_h^r$ in \eqref{dscheme1:3} and
sum up over all elements $K_x$, we get
\begin{align*}
&\displaystyle\int_\Ox  \frac{\bE_h^{n+1}-\bE_h^n}{\Delta t} \cdot (\bE_h^{n+1}+\bE_h^n)d\bx+
\int_\Ox  \frac{\bB_h^{n+1}-\bB_h^n}{\Delta t} \cdot (\bB_h^{n+1}+\bB_h^n)d\bx\\
&\displaystyle=2\int_\Ox \overline{\bB}\cdot\nabla_x\times \overline{\bE}d\bx
                          -2\int_\Ox \overline{\bE}\cdot\nabla_x\times \overline{\bB}d\bx
                          +2\int_{\mE_x}\widetilde{\overline{\bB}}\cdot[\overline{\bE}]_\tau ds_x-2\int_{\mE_x}\widetilde{\overline{\bE}}\cdot[\overline{\bB}]_\tau ds_x
                          -2\int_\Ox  \overline{\bE}\cdot\bJ_h^{n+1/2}   d\bx\\
                         &=2 \int_{\mE_x}\left([\overline{\bE}\times\overline{\bB}]_x+\widetilde{\overline{\bB}}\cdot[\overline{\bE}]_\tau-\widetilde{\overline{\bE}}\cdot[\overline{\bB}_\tau]\right)ds_x-2\int_\Ox  \overline{\bE}\cdot\bJ_h^{n+1/2}   d\bx.
\end{align*}
Using the definitions of averages and  jumps, clearly we have the following identities
\begin{subequations}
\begin{equation}
 [\bU\times\bV]_x+\{\bV\}_x\cdot[\bU]_\tau-\{\bU\}_x\cdot[\bV]_\tau=0,
 \end{equation}
 \begin{equation}
 [\bU\times\bV]_x+\bV^+\cdot[\bU]_\tau-\bU^-\cdot[\bV]_\tau=0,
 \end{equation}
 \begin{equation}
 [\bU\times\bV]_x+\bV^-\cdot[\bU]_\tau-\bU^+\cdot[\bV]_\tau=0.
 \end{equation}
 \end{subequations}
 Thus $$[\overline{\bE}\times\overline{\bB}]_x+\widetilde{\overline{\bB}}\cdot[\overline{\bE}]_\tau-\widetilde{\overline{\bE}}\cdot[\overline{\bB}_\tau]=0$$
 holds for   both the central \eqref{central:2} and the alternating fluxes \eqref{alternating:1}, \eqref{alternating:2} in the Maxwell solver.
Therefore, 
\begin{align*}
\displaystyle\int_\Ox \int_\Ov f_h^{n+1} |\bv|^2 d\bv d\bx +\int_\Ox |\bE_h^{n+1}|^2 d\bx+\int_\Ox |\bB_h^{n+1}|^2 d\bx\\[2mm]
  \displaystyle=\int_\Ox \int_\Ov f_h^n |\bv|^2 d\bv d\bx +\int_\Ox |\bE_h^n|^2 d\bx+\int_\Ox |\bB_h^n|^2 d\bx
\end{align*}
and we are done.
 \hfill$\Box$

Likewise for the discussion in the previous section, It is challenging to obtain the explicit form of the numerical energy for $\textnormal{\bf Scheme-3F}(\Delta t)$.

On the other hand, for implicit schemes for Vlasov equation, i.e. $\textnormal{\bf Scheme-3}(\Delta t)$, 
$\textnormal{\bf Scheme-4}(\Delta t)$, $\textnormal{\bf Scheme-3F}(\Delta t)$ and $\textnormal{\bf Scheme-4F}(\Delta t)$, fully discrete $L^2$ stability can be established. Clearly this property is independent of choice of numerical fluxes in the Maxwell solver.
\begin{thm}[{$L^2$ stability}]
\label{thm:l2stability}
The  DG methods with time integrators $\textnormal{\bf Scheme-3}(\Delta t)$, 
$\textnormal{\bf Scheme-4}(\Delta t)$, $\textnormal{\bf Scheme-3F}(\Delta t)$ and $\textnormal{\bf Scheme-4F}(\Delta t)$ satisfy
$$\int_\Ox \int_\Ov |f_h^\na|^2 d\bv d\bx =\int_\Ox \int_\Ov |f_h^n|^2  d\bv d\bx $$
for central flux, and
$$\int_\Ox \int_\Ov |f_h^\na|^2 d\bv d\bx \leq \int_\Ox \int_\Ov |f_h^n|^2  d\bv d\bx $$
for upwind flux (Again for the negative time steps appearing in $\textnormal{\bf Scheme-3F}(\Delta t)$ and $\textnormal{\bf Scheme-4F}(\Delta t)$, we require the flux to be the downwind flux instead.).
\end{thm}
\emph{Proof.}  The proof is straightforward by taking the test function to be $\frac{1}{2}(f_h^n+f_h^{n+1})$ and is omitted.  \hfill$\Box$ 

\subsection{Split schemes and their properties}
In this subsection, we describe  fully discrete implicit schemes with operator splitting {\bf Scheme-5}$(\Delta t)$, {\bf Scheme-5F}$(\Delta t)$  and discuss their properties. The key idea is to solve each split equation in their respective reduced dimensions.

Below let's introduce some notations first. We look for
$f_h \in \mS_h^{k}$, for which
 we can pick a few nodal points to represent the degree of freedom for that element \cite{hesthaven2007nodal}. Suppose the nodes in $\Kx$ and $\Kv$ are  $\bx_\Kx^{(l)}$, $\bv_\Kv^{(m)}$, $l=1, \ldots, dof(k1)$, $m=1, \ldots, dof(k2)$,  respectively, then any $g \in \mS_h^{k}$ can be uniquely represented as $g=\sum_{l,m} g(\bx_\Kx^{(l)}, \bv_\Kv^{(m)}) L_x^{(l)}(\bx) L_v^{(m)}(\bv)$ on $K$, where $L_x^{(l)}(\bx), L_v^{(m)}(\bv)$ denote the $l$-th and $m$-th Lagrangian interpolating polynomials in $\Kx$ and $\Kv$, respectively.

Under this setting, the equations for $f$ in the split equations (a), (b), (c) can be solved in reduced dimensions. For example,  equation (a), we can fix a nodal point in $\bv$, say $\bv_\Kv^{(m)}$, then solve $\partial_t f (\bv_\Kv^{(m)})  +  \bv_\Kv^{(m)}\cdot \nabla_{\bx} f  (\bv_\Kv^{(m)})= 0$ by a DG method in the $\bx$ direction. We can use the time integrator discussed in the previous subsection,  and get an update of point values at $f(\bx_\Kx^{(l)}, \bv_\Kv^{(m)})$ for all $\Kx, l$.

The idea is similar for equation (b). We can fix a nodal point in $\bx$, say $\bx_\Kx^{(l)}$, then solve
\[ \left\{
  \begin{array}{l}
\displaystyle\partial_t f (\bx_\Kx^{(l)})+   \bE(\bx_\Kx^{(l)})  \cdot \nabla_\bv f  (\bx_\Kx^{(l)}) = 0~, \\[3mm]
\partial_t \bE (\bx_\Kx^{(l)}) = -\bJ (\bx_\Kx^{(l)}),
  \end{array} \right.
  \]
 in $\bv$ direction, and get an update of point values at $f(\bx_\Kx^{(l)}, \bv_\Kv^{(m)})$ for all $\Kv, m$.
 
Similarly, for  equation (c), we first  solve the following system
\[
\left\{
  \begin{array}{l}\displaystyle\partial_t \bE = \nabla_\bx \times \bB, \\
\partial_t \bB =  -  \nabla_\bx \times \bE,
  \end{array} \right.
  \]
  on $\Ox$.
Then we fix a nodal point in $\bx$, say $\bx_\Kx^{(l)}$, and use the computed magnetic field $\bB(\bx_\Kx^{(l)})$ to solve
\[ 
\partial_t f (\bx_\Kx^{(l)})+ \left( \bv\times\bB(\bx_\Kx^{(l)})\right)  \cdot \nabla_\bv f  (\bx_\Kx^{(l)}) = 0~, 
  \]
 in $\bv$ direction, and get an update of point values at $f(\bx_\Kx^{(l)}, \bv_\Kv^{(m)})$ for all $\Kv, m$. This procedure is quite general and can be implemented on  unstructured meshes on $\Ox$, $\Ov$, if the nodal points are defined to guarantee the accuracy of the methods.

For simplicity of discussion, for the remaining of this section we will only consider the VM system in a simple 1D2V setting on a Cartesian mesh.
The VM system now becomes
 \begin{subequations}
 \label{example}
\begin{align}
&  f_t + v_2 f_{x_2} + (E_1 + v_2 B_3)f_{v_1} + (E_2 -
v_1 B_3 )f_{v_2} = 0 \\
&\frac{\df B_3}{\df t} =  \frac{\df E_1}{\df x_2},\\
& \frac{\df E_1}{\df t} =  \frac{\df B_3}{\df x_2} -  j_1, \\
& \frac{\df E_2}{\df t} =  -  j_2,
\end{align}
\end{subequations}
where \beq j_1=\int_{-V_{2,c}}^{V_{2,c}} \int_{-V_{1,c}}^{V_{1,c}}
f(x_2, v_1, v_2, t) v_1 \,dv_1 dv_2,\quad j_2=\int_{-V_{2,c}}^{V_{2,c}} \int_{-V_{1,c}}^{V_{1,c}} f(x_2, v_1,
v_2, t) v_2 \,dv_1 dv_2.\eeq
Here,  $f=f(x_2, v_1, v_2,t)$, $\bE(x_2,t) = (E_1(x_2,t),E_2(x_2,t),0)$ and  $\bB(x_2,t) = (0,0, B_3(x_2,t))$. The computational domain is $\Omega=\Omega_{x_2}\times\Omega_{v_1}
\times\Omega_{v_2}=[0,L]\times [-V_{1,c}, V_{1,c}]\times [-V_{2,c}, V_{2,c}]$, where $V_{1,c},\, V_{2,c}$  are chosen appropriately large to guarantee $f$ vanishes at
 $\partial \Ov$. The mesh is  partitioned as follows:
\begin{align*}
 0&=x_{2,\frac{1}{2}}<x_{2,\frac{3}{2}}< \ldots
<x_{2,N_x+\frac{1}{2}}=L ,
 \qquad \\-V_{1,c}&=v_{1,\frac{1}{2}}<v_{1,\frac{3}{2}}<
\ldots <v_{1,N_{v_1}+\frac{1}{2}}=V_{1,c},
 \qquad \\-V_{2,c}&=v_{2,\frac{1}{2}}<v_{2,\frac{3}{2}}<
\ldots <v_{2,N_{v_2}+\frac{1}{2}}=V_{2,c},
\end{align*}
  The elements are defined as
\begin{eqnarray*}
\label{2dcell2} & K_{i,j_1,j_2}=[x_{2,i-\frac{1}{2}},x_{2,i+\frac{1}{2}}]
\times [v_{1,j_1-\frac{1}{2}},v_{1,j_1+\frac{1}{2}}] \times [v_{2,j_2-\frac{1}{2}},v_{2,j_2+\frac{1}{2}}] , 
\quad K_{x_2,i}=[x_{2,i-1/2},x_{2,i+1/2}],\\
& K_{v_1,j_1}=[v_{1,j_1-1/2},v_{1,j_1+1/2}]\,,  \quad K_{v_2,j_2}=[v_{2,j_2-1/2},v_{2,j_2+1/2}],  \end{eqnarray*}
for $  i=1,\ldots N_{x_2}, \; j_1=1,\ldots N_{v_1} , \;j_2=1,\ldots N_{v_2} .
$
Let  $\Delta x_{2,i}= x_{2,i+1/2}-x_{2,i-1/2}$, $\Delta v_{1,j_1}=v_{1,j_1+1/2}-v_{1,j_1-1/2}$, $\Delta v_{2,j_2}=v_{2,j_2+1/2}-v_{2,j_2-1/2}$ be the length of each interval. $x_{2,i}^{(l)}$  for $l=1, \ldots, k+1$ be the $(k+1)$ Gauss quadrature points on $K_{x_2,i}$, and $v_{1,j_1}^{(m_1)}$ for $m_1=1, \ldots, k+1$ be the $(k+1)$ Gauss quadrature points on $K_{v_1,j_1}$ and $v_{2,j_2}^{(m_2)}$ for $m_2=1, \ldots, k+1$ be the $(k+1)$ Gauss quadrature points on $K_{v_2,j_2}$. Now we are ready to describe our scheme for each split equation.
\bigskip

\underline{Algorithm {\bf Scheme-a}$(\Delta t)$}
\medskip

To solve equation (a) from $t^n$ to $t^\na$
 \[
  \textrm{ (a)} \left\{
  \begin{array}{l}
\displaystyle\partial_t f  +  v_2 f_{x_2}   = 0~, \\
\partial_t E_1 = 0,\\
\partial_t E_2 = 0,\\
\partial_t B_3 = 0,
  \end{array} \right.
\]
\begin{enumerate}
\item For each $j_1=1,\ldots N_{v_1} , \,j_2=1,\ldots N_{v_2},\, m_1=1, \ldots, k+1,\, m_2=1, \ldots, k+1$, 
we seek $g_{j_1,j_2}^{(m_1,m_2)}(x_2) \in \mW_h^k$, such that
\begin{eqnarray}
 &\displaystyle\int_{K_{x_2,i}} \frac{g_{j_1,j_2}^{(m_1,m_2)}(x_2)-f_h^n(x_2, v_{1,j_1}^{(m_1)}, v_{2,j_2}^{(m_2)} )}{\Delta t} \varphi_h \, dx_2\notag\\
   &\displaystyle- \int_{K_{x_2,i}} v_{2,j_2}^{(m_2)}\frac{g_{j_1,j_2}^{(m_1,m_2)}(x_2)+f_h^n(x_2, v_{1,j_1}^{(m_1)}, v_{2,j_2}^{(m_2)}  )}{2} (\varphi_h)_{x_2} \, dx_2 \notag\\[2mm]
 &\displaystyle+ \reallywidehat{\displaystyle v_{2,j_2}^{(m_2)} \frac{g_{j_1,j_2}^{(m_1,m_2)}(x_{2,i+\frac{1}{2}})+f_h^n(x_{2,i+\frac{1}{2}}, v_{1,j_1}^{(m_1)}, v_{2,j_2}^{(m_2)} )}{2} } (\varphi_h)_{i+\frac{1}{2}}^- \label{schemeas} \\[2mm]&\displaystyle
    -\reallywidehat{\displaystyle v_{2,j_2}^{(m_2)} \frac{g_{j_1,j_2}^{(m_1,m_2)}(x_{2,i-\frac{1}{2}})+f_h^n(x_{2,i-\frac{1}{2}},v_{1,j_1}^{(m_1)}, v_{2,j_2}^{(m_2)}  )}{2} } (\varphi_h)_{i-\frac{1}{2}}^+ =0 \notag
\end{eqnarray}
holds for any test function $\varphi_h(x_2) \in \mW_h^k$.
\item Let $f_h^{n+1}$ be the unique polynomial in $\mS_h^k$, such that 
$$f_h^\na(x_{2,i}^{(l)},v_{1,j_1}^{(m_1)}, v_{2,j_2}^{(m_2)} )=g_{j_1,j_2}^{(m_1,m_2)}(x_{2,i}^{(l)}), \quad \forall i, j_1,j_2,\,l, m_1,m_2.$$
\end{enumerate}
\bigskip

\underline{Algorithm {\bf Scheme-b}$(\Delta t)$}
\medskip

To solve equation (b) from $t^n$ to $t^\na$
\[
  \textrm{(b)} \left\{
  \begin{array}{l}
\partial_t f  +   E_1f_{v_1}+E_2f_{v_2}   = 0~, \\
\partial_t E_1 = -j_1,\\
\partial_t E_2 = -j_2,\\
\partial_t B_3 = 0,
  \end{array} \right.
\]

\begin{enumerate}
\item For each $i=1,\ldots N_{x_2} , l=1, \ldots, k+1$, we seek $g_i^{(l)}(v_1,v_2) \in \mZ_h^k$, $E_{1,i}^{(l)}$, and $E_{2,i}^{(l)}$, such that
 for any test function $\varphi_h(v_1,v_2) \in \mZ_h^k$, we have 
 \begin{align} 
&\hspace{-15mm} \int_{K_{v_2,j_2}}  \int_{K_{v_1,j_1}}\frac{g_i^{(l)}(v_1,v_2)-f_h^n(x_{2,i}^{(l)}, v_1,v_2)}{\Delta t} \varphi_h \,  dv_1dv_2\notag\\[2mm]
& \hspace{-15mm}-\int_{K_{v_2,j_2}}  \int_{K_{v_1,j_1}}\frac{g_i^{(l)}(v_1,v_2)+f_h^n(x_{2,i}^{(l)}, v_1,v_2)}{2}  
    \left(\frac{E_{1,h}^n(x_{2,i}^{(l)})+E_{1,i}^{(l)}}{2} (\varphi_h)_{v_1}+\frac{E_{1,h}^n(x_{2,i}^{(l)})+E_{1,i}^{(l)}}{2} (\varphi_h)_{v_1}\right) dv_1dv_2 \notag\\
  &\hspace{-15mm} + \int_{K_{v_2,j_2}}\left(\frac{E_{1,h}^n(x_{2,i}^{(l)})+E_{1,i}^{(l)}}{2} \right)\displaystyle\reallywidehat{\displaystyle\frac{g_i^{(l)}(v_{1,j_1+\frac{1}{2}},v_2)
  +f_h^n(x_{2,i}^{(l)}, v_{1,j_1+\frac{1}{2}},v_2)}{2}} \varphi_h\left(v^-_{1, j_1+\frac{1}{2}},v_2\right)dv_2 \notag\\
  &\hspace{-15mm} - \int_{K_{v_2,j_2}}\left(\frac{E_{1,h}^n(x_{2,i}^{(l)})+E_{1,i}^{(l)}}{2}\right)\reallywidehat{\displaystyle\frac{g_i^{(l)}(v_{1,j_1-\frac{1}{2}},v_2)
  +f_h^n(x_{2,i}^{(l)}, v_{1,j_1-\frac{1}{2}},v_2)}{2}} \varphi_h\left(v^+_{1, j_1-\frac{1}{2}},v_2\right)dv_2  \label{scheme5b}\\
  &\hspace{-15mm}+ \int_{K_{v_1,j_1}}\left(\frac{E_{2,h}^n(x_{2,i}^{(l)})+E_{2,i}^{(l)}}{2}\right) \reallywidehat{\displaystyle\frac{g_i^{(l)}(v_1,v_{2,j_2+\frac{1}{2}})
  +f_h^n(x_{2,i}^{(l)}, v_1, v_{2,j_2+\frac{1}{2}})}{2}} \varphi_h\left(v_1,v^-_{2, j_2+\frac{1}{2}}\right)dv_1 \notag\\
  &\hspace{-15mm}- \int_{K_{v_1,j_1}}\left(\frac{E_{2,h}^n(x_{2,i}^{(l)})+E_{2,i}^{(l)}}{2}\right) \reallywidehat{\displaystyle\frac{g_i^{(l)}(v_1,v_{2,j_2-\frac{1}{2}})
  +f_h^n(x_{2,i}^{(l)}, v_1, v_{2,j_2-\frac{1}{2}})}{2}} \varphi_h\left(v_1,v^+_{2, j_2-\frac{1}{2}}\right)dv_1=0 \notag
\end{align}
\begin{align*}
&\frac{E_{1,i}^{(l)}-E_{1,h}^n(x_{2,i}^{(l)})}{\Delta t}=  -\frac{1}{2}(J_{1,h}^n(x_{2,i}^{(l)})+J_{1,i}^{(l)}), \; \notag\\
&\frac{E_{2,i}^{(l)}-E_{2,h}^n(x_{2,i}^{(l)})}{\Delta t}=  -\frac{1}{2}(J_{2,h}^n(x_{2,i}^{(l)})+J_{2,i}^{(l)}),
\end{align*}
where
$$ J_{1,h}^n(x_2)=\int_{\Omega_{v_2}} \int_{\Omega_{v_1}} f_h^n(x_2, v_1,v_2) v_1 \,dv_1dv_2, \quad J_{1,i}^{(l)}=\int_{\Omega_{v_2}} \int_{\Omega_{v_1}}g_i^{(l)}(v_1,v_2) v_1 \,dv_1dv_2, $$
$$ J_{2,h}^n(x_2)=\int_{\Omega_{v_2}} \int_{\Omega_{v_1}} f_h^n(x_2, v_1,v_2) v_2 \,dv_1dv_2, \quad J_{2,i}^{(l)}=\int_{\Omega_{v_2}} \int_{\Omega_{v_1}}g_i^{(l)}(v_1,v_2) v_2 \,dv_1dv_2.$$
\item Let $f_h^{n+1}$ be the unique polynomial in $\mS_h^k$, such that 
$$f_h^\na(x_{2,i}^{(l)},v_{1,j_1}^{(m_1)}, v_{2,j_2}^{(m_2)} )=g_i^{(l)}(v_{j_1}^{(m_1)},v_{j_2}^{(m_2)}), \quad \forall i, j_1,j_2,\,l, m_1,m_2.$$
Let $E_{1,h}^{n+1},\,E_{2,h}^{n+1}$ be the unique polynomials in $\mW_h^k$, such that
$$E_{1,h}^\na(x_{2,i}^{(l)})=E_{1,i}^{(l)},\quad\quad E_{2,h}^\na(x_{2,i}^{(l)})=E_{2,i}^{(l)} \quad \forall i, l.$$
\end{enumerate}
\bigskip

\underline{Algorithm {\bf Scheme-c}$(\Delta t)$}
\medskip

To solve equation (c) from $t^n$ to $t^\na$,

\[
  \textrm{(c)}\left\{
  \begin{array}{l}
\partial_t f  +  v_2B _3f_{v_1}-v_1B_3\,f_{v_2}   = 0~, \\
\displaystyle\partial_t E_1 =(B_3)_{x_2}\\
\partial_t B_3 =   (E_1)_{x_2},\\
\partial_t E_2 = 0,
  \end{array} \right.\]
\begin{enumerate}
\item First to  solve the Maxwell's equation, we seek $E_{1,h}^{n+1},\,B_{3,h}^{n+1}\in \mW_h^k$, such that
\begin{equation}
\label{schemec:e}
\begin{split}
& \hspace{-6mm}\int_{K_{x_2,i}} \frac{E_{1,h}^{n+1}-E_{1,h}^n}{\Delta t} \varphi_h \, dx_2 
=- \int_{K_{x_2,i}}\frac{B_{3,h}^{n+1}+B_{3,h}^{n}}{2} (\varphi_h)_{x_2} \, dx_2  \\[3mm]
& \hspace{-6mm}+ \reallywidehat{\displaystyle\frac{B_{3,h}^{n+1}(x_{2,i+\frac{1}{2}})+B_{3,h}^{n}(x_{2,i+\frac{1}{2}})}{2} } (\varphi_h)_{i+\frac{1}{2}}^- 
    - \reallywidehat{\displaystyle\frac{B_{3,h}^{n+1}(x_{2,i-\frac{1}{2}})+B_{3,h}^{n}(x_{2,i-\frac{1}{2}})}{2} }  (\varphi_h)_{i-\frac{1}{2}}^+ =0 ,
\end{split}
\end{equation}
\begin{equation}
\label{schemec:b}
\begin{split}
  &\hspace{-6mm}\int_{K_{x_2,i}} \frac{B_{3,h}^{n+1}-B_{3,h}^n}{\Delta t} \phi_h \, dx_2 
 =- \int_{K_{x_2,i}}\frac{E_{1,h}^{n+1}+E_{1,h}^{n}}{2} (\phi_h)_{x_2} \, dx_2   \\[3mm] 
&\hspace{-6mm}+ \reallywidehat{\displaystyle\frac{E_{1,h}^{n+1}(x_{2,i+\frac{1}{2}})+E_{1,h}^{n}(x_{2,i+\frac{1}{2}})}{2} } (\phi_h)_{i+\frac{1}{2}}^- 
    - \reallywidehat{\displaystyle\frac{E_{1,h}^{n+1}(x_{2,i-\frac{1}{2}})+E_{1,h}^{n}(x_{2,i-\frac{1}{2}})}{2} }  (\phi_h)_{i-\frac{1}{2}}^+ =0,
\end{split}
\end{equation}
holds for any test function $\varphi_h(x_2), \phi_h(x_2) \in \mW_h^k$.
\item For each $i=1,\ldots N_{x_2} , l=1, \ldots, k+1$, denote $$B(x_{i}^{(l)})=\frac{B_{3,h}^{n+1}(x_{2,i}^{(l)})+B_{3,h}^n(x_{2,i}^{(l)})}{2},$$ we seek $g_i^{(l)}(v_1,v_2) \in \mZ_h^k$, such that
 for any test function $\varphi_h(v_1,v_2) \in \mZ_h^k$, we have 
\begin{align}
&  \int_{K_{v_1,j_1}}\int_{K_{v_2,j_2}} \frac{g_i^{(l)}(v_1,v_2)-f_h^n(x_{2,i}^{(l)}, v_1,v_2)}{\Delta t} \varphi_h \,  dv_2dv_1\notag\\[2mm]
&  -\int_{K_{v_1,j_1}}\int_{K_{v_2,j_2}} \frac{g_i^{(l)}(v_1,v_2)+f_h^n(x_{2,i}^{(l)}, v_1,v_2)}{2}  
    \left(v_2B(x_i^{(l)}) (\varphi_h)_{v_1}-v_1B(x_i^{(l)})(\varphi_h)_{v_2}\right) dv_2dv_1 \notag\\
  & + \int_{K_{v_2,j_2}}v_2B(x_i^{(l)}) \reallywidehat{\displaystyle\frac{g_i^{(l)}(v_{1,j_1+\frac{1}{2}},v_2)
  +f_h^n(x_{2,i}^{(l)}, v_{1,j_1+\frac{1}{2}},v_2)}{2}} \varphi_h\left(v^-_{1, j_1+\frac{1}{2}},v_2\right)dv_2 \notag\\
  & - \int_{K_{v_2,j_2}}v_2B(x_i^{(l)})\reallywidehat{\displaystyle\frac{g_i^{(l)}(v_{1,j_1-\frac{1}{2}},v_2)
  +f_h^n(x_{2,i}^{(l)}, v_{1,j_1-\frac{1}{2}},v_2)}{2}} \varphi_h\left(v^+_{1, j_1-\frac{1}{2}},v_2\right)dv_2 \label{scheme5c}\\
  & -\int_{K_{v_1,j_1}}v_1B(x_i^{(l)})\reallywidehat{\displaystyle\frac{g_i^{(l)}(v_1,v_{2,j_2+\frac{1}{2}})
  +f_h^n(x_{2,i}^{(l)}, v_1, v_{2,j_2+\frac{1}{2}})}{2}} \varphi_h\left(v_1,v^-_{2, j_2+\frac{1}{2}}\right)dv_1 \notag\\
  &+\int_{K_{v_1,j_1}}v_1B(x_i^{(l)}) \reallywidehat{\displaystyle\frac{g_i^{(l)}(v_1,v_{2,j_2-\frac{1}{2}})
  +f_h^n(x_{2,i}^{(l)}, v_1, v_{2,j_2-\frac{1}{2}})}{2}} \varphi_h\left(v_1,v^+_{2, j_2-\frac{1}{2}}\right)dv_1 =0.\notag
\end{align}
\item Let $f_h^{n+1}$ be the unique polynomial in $\mS_h^k$, such that 
$$f_h^\na(x_{2,i}^{(l)},v_{1,j_1}^{(m_1)}, v_{2,j_2}^{(m_2)} )=g_i^{(l)}(v_{j_1}^{(m_1)},v_{j_2}^{(m_2)}), \quad \forall i, j_1,j_2,\,l, m_1,m_2.$$
\end{enumerate}

Similar to the discussion in Section \ref{sec:unsplit}, the flux terms in the algorithms above can be taken as either upwind or central flux for Vlasov solver, and either   central or alternating flux for Maxwell solver. Finally, we recall that the method for the full VM system is defined as
$$\hspace{-10mm}\textnormal{\bf Scheme-5}(\Dt)=\textnormal{\bf Scheme-a}(\Dt/2)\textnormal{\bf Scheme-b}(\Dt/2)\textnormal{\bf Scheme-c}(\Dt)\textnormal{\bf Scheme-b}(\Dt/2)\textnormal{\bf Scheme-a}(\Dt/2),$$
and the corresponding  fourth-order-in-time method is
$$\textnormal{\bf Scheme-5F}(\Delta t)=\textnormal{\bf Scheme-5}(\beta_1 \Delta t) \textnormal{\bf Scheme-5}(\beta_2 \Delta t) \textnormal{\bf Scheme-5}(\beta_3 \Delta t).$$

\medskip
Next, we will discuss the conservation properties of the fully discrete schemes with operator splitting.

\begin{thm}[{Total particle number conservation}]
The DG schemes with  time integrators $\textnormal{\bf Scheme-5}(\Delta t)$ and $\textnormal{\bf Scheme-5F}(\Delta t)$ as described in this section
  preserve the total particle number of the system, i.e.
$$\int_{\Omega_{x_2}} \int_{\Omega_{v_2}}\int_{\Omega_{v_1}} f_h^\na dv_1dv_2 dx_2 =\int_{\Omega_{x_2}} \int_{\Omega_{v_2}}\int_{\Omega_{v_1}} f_h^n  dv_1dv_2 dx_2 .$$
\end{thm}
\emph{Proof.}  We only need to prove  conservation for each of the operators $\textnormal{\bf Scheme-a}(\Delta t)$, $\textnormal{\bf Scheme-b}(\Delta t)$ and 
$\textnormal{\bf Scheme-c}(\Delta t)$. For $\textnormal{\bf Scheme-a}(\Delta t)$,
let $\varphi_h=1$ in \eqref{schemeas}, and sum over all elements $K_{x_2,i}$, we get
$$\int_{\Omega_{x_2}}  g_{j_1,j_2}^{(m_1,m_2)}(x_2) dx_2 =\int_{\Omega_{x_2}} f_h^n(x_2, v_{1,j_1}^{(m_1)}, v_{2,j_2}^{(m_2)} ) dx_2 .$$
Therefore for any $j_1,\,j_2,\,m_1,\, m_2$,
 $$\int_{\Omega_{x_2}}  f_h^{n+1}(x_2, v_{1,j_1}^{(m_1)}, v_{2,j_2}^{(m_2)} )   dx_2
   =\int_{\Omega_{x_2}}  f_h^n(x_2, v_{1,j_1}^{(m_1)}, v_{2,j_2}^{(m_2)} )  dx_2 .$$
Since the (k+1)-point Gauss quadrature formula is exact for polynomial with degree less than $2k+2$, we have
\begin{eqnarray*}
\int_{\Omega_{x_2}} \int_{\Omega_{v_2}}\int_{\Omega_{v_1}} f_h^\na dv_1dv_2 dx_2 
=\sum_{j_1} \sum_{j_2}  \sum_{m_1} \sum_{m_2} w_{m_1}w_{m_2} \int_{\Omega_{x_2}}   f_h^{n+1}(x_2, v_{1,j_1}^{(m_1)}, v_{2,j_2}^{(m_2)} ) dx_2 
\,\Delta v_{1,j_1}\Delta v_{2,j_2}\\
=\sum_{j_1} \sum_{j_2}  \sum_{m_1} \sum_{m_2} w_{m_1}w_{m_2} \int_{\Omega_{x_2}}   f_h^{n}(x_2, v_{1,j_1}^{(m_1)}, v_{2,j_2}^{(m_2)} ) dx_2 
\,\Delta v_{1,j_1}\Delta v_{2,j_2}=\int_{\Omega_{x_2}} \int_{\Omega_{v_2}}\int_{\Omega_{v_1}} f_h^n  dv_1dv_2 dx_2 ,
\end{eqnarray*}
 where $w_{m_1}, w_{m_2}$ are the corresponding Gauss quadrature weights. The proof is similar for $\textnormal{\bf Scheme-b}(\Delta t)$ 
 and $\textnormal{\bf Scheme-c}(\Delta t)$ and is omitted.
 $\Box$
 
\begin{thm}[{Total energy conservation}]
If $k \geq 2$, the   DG schemes with time integrators $\textnormal{\bf Scheme-5}(\Delta t)$ and $\textnormal{\bf Scheme-5F}(\Delta t)$  preserve the discrete total energy $TE_n=TE_{n+1}$, where
$$2 \,(TE_n)= \int_{\Omega_{x_2}} \int_{\Omega_{v_2}}\int_{\Omega_{v_1}} f_h^n (v_1^2+v_2^2) dv_1dv_2 dx_2 +\int_{\Omega_{x_2}}  |E_{1,h}^n|^2+|E_{2,h}^n|^2+|B_{3,h}^n|^2 dx_2.$$
\end{thm}
\emph{Proof.}  We  need to show that the discrete total energy conservation for each of the operators $\textnormal{\bf Scheme-a}(\Delta t)$, $\textnormal{\bf Scheme-b} (\Delta t)$ and $\textnormal{\bf Scheme-c}(\Delta t)$. The proof for  $\textnormal{\bf Scheme-a}(\Delta t)$  and $\textnormal{\bf Scheme-b}(\Delta t)$ is similar to
the proof of the split fully discrete schemes in \cite{cheng_va} and is omitted.

As for $\textnormal{\bf Scheme-c}(\Delta t)$, for the simplicity of description, we first introduce some short-hand notations:
\begin{subequations}
\label{notation}
\begin{align}
&\hspace{-6mm}E_{1,i}^-=\frac{1}{2}\left(E_{1,h}^n(x_{2,i+1/2}^-)+E_{1,h}^{n+1}(x_{2,i+1/2}^-)\right), \; 
  E_{1,i}^+=\frac{1}{2}\left(E_{1,h}^n(x_{2,i+1/2}^+)+E_{1,h}^{n+1}(x_{2,i+1/2}^+)\right);\\[2mm]
&\hspace{-6mm}B_{i}^-    =\frac{1}{2}\left(B_{3,h}^n(x_{2,i+1/2}^-)+B_{3,h}^{n+1}(x_{2,i+1/2}^-)\right),\;\;
     B_{i}^+=\frac{1}{2}\left(B_{3,h}^n(x_{2,i+1/2}^+)+B_{3,h}^{n+1}(x_{2,i+1/2}^+)\right);\\[3mm]
 &\hspace{-6mm} \widehat{E_{1,i}}=\reallywidehat{\displaystyle\frac{E_{1,h}^{n+1}(x_{2,i+\frac{1}{2}})+E_{1,h}^{n}(x_{2,i+\frac{1}{2}})}{2} },\quad\quad \quad\quad
 \widehat{B_i}= \reallywidehat{\displaystyle\frac{B_{3,h}^{n+1}(x_{2,i+\frac{1}{2}})+B_{3,h}^{n}(x_{2,i+\frac{1}{2}})}{2} }.
\end{align}
\end{subequations}

Let
$\varphi_h= E_{1,h}^{n+1}+E_{1,h}^{n}\in\mW_h^k$ in \eqref{schemec:e} and $\phi_h = B_{3,h}^{n+1}+B_{3,h}^{n}\in\mW_h^k$ in \eqref{schemec:b},  and
sum up over all element $K_{x_2,i}$, using the notations of \eqref{notation}, we get
\begin{align}
 &\int_{\Omega_{x_2}} \frac{E_{1,h}^{n+1}-E_{1,h}^n}{\Delta t} \left( E_{1,h}^{n+1}+E_{1,h}^{n}\right)\, dx_2 + 
 \int_{\Omega_{x_2}} \frac{B_{3,h}^{n+1}-B_{3,h}^n}{\Delta t} \left( B_{3,h}^{n+1}+B_{3,h}^{n}\right) dx_2\notag \\
&  = 2\sum_i\left(E_{1,i}^+B_{i}^+ - E^-_{1,i}B_{i}^-+\widehat{B_{i}}(E_{1,i}^--E_{1,i}^+)+\widehat{E_{i}}(B_{1,i}^--B_{1,i}^+)\right) \label{schemec:ebp}.
\end{align}
The numerical fluxes for  \eqref{schemec:e}-\eqref{schemec:b} are defined as follows
\begin{subequations} 
\label{flux:schemec}
  \begin{align}
  \textnormal{central:} &  \quad      \widehat{E_{1,i}}=\frac{E_{1,i}^++E_{1,i}^-}{2},\quad\quad
       \widehat{B_i}=\frac{B_i^++B_i^-}{2};\label{flux:schemec:2}\\
\textnormal{alternating:}&   \quad       \widehat{E_{1,i}}=E_{1,i}^+,\;\; \widehat{B_i}= B_i^-; \;\; \textnormal{or}\;\;   \widehat{E_{1,i}}=E_{1,i}^-,\;\; \widehat{B_i}= B_i^+.\label{flux:schemec:3}
    \end{align}
\end{subequations}
Substituting \eqref{flux:schemec} into \eqref{schemec:ebp}, we get
\begin{align}
  & \int_{\Omega_{x_2}}\frac{E_{1,h}^{n+1}-E_{1,h}^n}{\Delta t} \left( E_{1,h}^{n+1}+E_{1,h}^{n}\right)\, dx_2 + 
     \int_{\Omega_{x_2}}\frac{B_{3,h}^{n+1}-B_{3,h}^n}{\Delta t} \left(B_{3,h}^{n+1}+B_{3,h}^{n}\right) dx_2=0
             \label{schemec:eb}
\end{align}

Let $\varphi_h=v^2=v_1^2+v_2^2$ in \eqref{scheme5c} and sum over all element $K_{v_1,j_1}\times K_{v_2,j_2}$, we get
$$\int_{\Omega_{v_2}}\int_{\Omega_{v_1}}  g_i^{(l)}(v_1,v_2) v^2dv_1dv_2 =\int_{\Omega_{v_2}}\int_{\Omega_{v_1}} f_h^n(x_i^{(l)}, v_1,v_2) ) v^2dv_1dv_2 .$$
 Note that $ f_h(x, v)  v^2$  and $g_i^{(l)}(v_1,v_2)v^2$ are polynomials that are at most degree $k+2$  in each variable of $v=(v_1,v_2)$ and $x_2$. Since the (k+1)-point Gauss quadrature formula is exact for polynomial with degree less than $2k+2$ and  $k+2 \leq 2k+1$ when $k \geq 2$,  we have
\begin{align*}
& \int_{\Omega_{x_2}} \int_{\Omega_{v_2}}\int_{\Omega_{v_1}} f_h^\na (v_1^2+v_2^2) dv_1dv_2 dx_2
\\&
 = \sum_i\sum_l\sum_{j_1,j_2}\sum_{m_1,m_2} w_lw_{m_1}w_{m_2}f_h^\na\left(x_{2,i}^{(l)},v_{1,j_1}^{(m_1)},v_{2,j_2}^{(m_2)}\right)
 \left(|v_{1,j_1}^{(m_1)}|^2+|v_{2,j_2}^{(m_2)}|^2\right)
\Delta x_{2,i}\Delta v_{1,j_1}\Delta v_{2,j_2}\\&
 = \sum_i\sum_l\sum_{j_1,j_2}\sum_{m_1,m_2} w_lw_{m_1}w_{m_2}g_i^{(l)}\left(v_{1,j_1}^{(m_1)},v_{2,j_2}^{(m_2)}\right)
 \left(|v_{1,j_1}^{(m_1)}|^2+|v_{2,j_2}^{(m_2)}|^2\right)
\Delta x_{2,i}\Delta v_{1,j_1}\Delta v_{2,j_2}\\
&
 = \sum_i\sum_lw_l\int_{\Omega_{v_2}}\int_{\Omega_{v_1}}   g_i^{(l)}(v_1,v_2) v^2dv_1dv_2\Delta x_{2,i}\\
&
 = \sum_i\sum_lw_l\int_{\Omega_{v_2}}\int_{\Omega_{v_1}}   f_h^n\left(x_{2,i}^{(l)}, v_1,v_2 \right)v^2dv_1dv_2\Delta x_{2,i}\\
 &= \int_{\Omega_{x_2}} \int_{\Omega_{v_2}}\int_{\Omega_{v_1}} f_h^n \left(v_1^2+v_2^2\right) dv_1dv_2 dx_2\end{align*}
Therefore, putting all the results together for $\textnormal{\bf Scheme-a}(\Delta t)$,  $\textnormal{\bf Scheme-b}(\Delta t)$ 
and $\textnormal{\bf Scheme-c}(\Delta t)$,  we are done. \hfill
 $\Box$

\begin{thm}[{$L^2$ stability}]
The DG schemes with  $\textnormal{\bf Scheme-5}(\Delta t)$, $\textnormal{\bf Scheme-5F}(\Delta t)$
satisfy
$$ \int_{\Omega_{x_2}} \int_{\Omega_{v_2}}\int_{\Omega_{v_1}} |f_h^\na|^2 dv_1dv_2 dx_2 =
     \int_{\Omega_{x_2}} \int_{\Omega_{v_2}}\int_{\Omega_{v_1}} |f_h^n|^2  dv_1dv_2 dx_2 $$
for central flux in Vlasov solver and
$$ \int_{\Omega_{x_2}} \int_{\Omega_{v_2}}\int_{\Omega_{v_1}} |f_h^\na|^2 dv_1dv_2 dx_2
 \leq \int_{\Omega_{x_2}} \int_{\Omega_{v_2}}\int_{\Omega_{v_1}}|f_h^n|^2  dv_1dv_2 dx_2 $$
for upwind flux in Vlasov solver (and again we use downwind flux for the negative time steps in  $\textnormal{\bf Scheme-5F}(\Delta t)$).
\end{thm}
\emph{Proof.}   We only need to prove the theorem for each of the operators $\textnormal{\bf Scheme-a}(\Delta t)$, $\textnormal{\bf Scheme-b}(\Delta t)$  and $\textnormal{\bf Scheme-c}(\Delta t)$. For  $\textnormal{\bf Scheme-a}(\Delta t)$,
let $\varphi_h=g_{j_1,j_2}^{(m_1,m_2)}(x_2)+ f_h^n(x_2,v_{1,j_1}^{(m_1)}, v_{2,j_2}^{(m_2)}) $ in \eqref{schemeas}, and sum over all element $K_{x_2,i}$, we get
$$ \int_{\Omega_{x_2}} \left| g_{j_1,j_2}^{(m_1,m_2)}(x_2)\right|^2 dx_2= \int_{\Omega_{x_2}} \left|f_h^n(x_2, v_{1,j_1}^{(m_1)}, v_{2,j_2}^{(m_2)})\right|^2 dx_2 $$
for central flux and
$$ \int_{\Omega_{x_2}} \left | g_{j_1,j_2}^{(m_1,m_2)}(x_2)\right|^2 dx_2 \leq  \int_{\Omega_{x_2}}\left |f_h^n(x_2, v_{1,j_1}^{(m_1)}, v_{2,j_2}^{(m_2)})\right|^2 dx_2 $$
for upwind flux. Therefore for any $j, m$,
 $$ \int_{\Omega_{x_2}} \left|f_h^\na(x_2, v_{1,j_1}^{(m_1)}, v_{2,j_2}^{(m_2)})\right|^2  dx_2 = \int_{\Omega_{x_2}} \left|f_h^n(x_2,v_{1,j_1}^{(m_1)}, v_{2,j_2}^{(m_2)})\right|^2 dx_2 ,$$
 for central flux and
  $$ \int_{\Omega_{x_2}}  \left|f_h^\na(x_2, v_{1,j_1}^{(m_1)}, v_{2,j_2}^{(m_2)})\right|^2  dx_2 \leq  \int_{\Omega_{x_2}} \left|f_h^n(x_2, v_{1,j_1}^{(m_1)}, v_{2,j_2}^{(m_2)})\right|^2 dx_2 ,$$
  for upwind flux.
 Since the (k+1) Gauss quadrature formula is exact for polynomial of  degree less than $2k+2$, we have
\begin{eqnarray*}
&\hspace{-8mm}\displaystyle\int_{\Omega_{x_2}} \int_{\Omega_{v_2}}\int_{\Omega_{v_1}} |f_h^\na|^2 dv_1dv_2 dx_2 =\sum_{j_1,j_2}\sum_{m_1,m_2} w_{m_1}w_{m_2}  \int_{\Omega_{x_2}} \left |f_h^\na(x_2, v_{1,j_1}^{(m_1)}, v_{2,j_2}^{(m_2)})\right|^2  dx_2 \,\Delta v_{1,j_1}\Delta v_{2,j_2}\\
&\hspace{-8mm}\displaystyle= (\textrm{or} \leq) \sum_{j_1,j_2} \sum_{m_1,m_2} w_{m_1}w_{m_2} \int_{\Omega_{x_2}}  \left|f_h^n(x_2,v_{1,j_1}^{(m_1)},v_{2,j_2}^{(m_2)})\right|^2  dx_2 \,\Delta v_{1,j_1}\Delta v_{2,j_2}\displaystyle= \int_{\Omega_{x_2}} \int_{\Omega_{v_2}}\int_{\Omega_{v_1}} |f_h^n|^2  dv_1dv_2 dx_2 ,
\end{eqnarray*}
  where $w_{m_1}, w_{m_2}$ are the corresponding Gauss quadrature weights.  The proof is similar for $\textnormal{\bf Scheme-b}(\Delta t)$
 and  $\textnormal{\bf Scheme-c}(\Delta t)$. \hfill$\Box$


In summary, the schemes with the operator splitting is fully implicit, energy conservative, and $L^2$ stable. Each of the split equation is only in $\bx$ or $\bv$ space, and can be computed efficiently.
Similar to \cite{cheng_va}, we need a
Jacobian-free Newton-Krylov solver \cite{knoll2004jacobian}   to compute the nonlinear systems resulting from $\textnormal{\bf Scheme-b}(\Delta t)$.

\section{Numerical Results}
\setcounter{equation}{0}
\setcounter{figure}{0}
\setcounter{table}{0}
\label{sec:numerical}

In this section, we show the numerical results of the proposed methods $\textnormal{\bf Scheme-1}(\Dt)$,  $\textnormal{\bf Scheme-2}(\Dt)$, $\textnormal{\bf Scheme-5}(\Dt)$ for the streaming Weibel instability \cite{califano1998ksw},  which is a reduced version of the VM system with the simple form \eqref{example}. 

The initial conditions   are given by 
\beq 
f(x_2, v_1, v_2, 0)=\frac{1}{\pi \beta}
e^{- v_2^2 /\beta} [\delta e^{- (v_1-v_{0,1})^2 /\beta}
+(1-\delta) e^{- (v_1+v_{0,2})^2 /\beta} ],
\eeq

\beq 
E_1(x_2,  0)=E_2(x_2,  0)=0, \qquad B_3(x_2, 0)=b \sin(k_0 x_2),\eeq
where $\beta^{1/2}$ is the thermal velocity and $\delta$ is a parameter measuring the symmetry of the electron beams ($\delta = 0.5$ in the symmetric case)
and $b$  is the amplitude of the initial perturbation to the magnetic field. $b=0$ is an equilibrium state composed of counter-streaming beams 
propagating perpendicular to the direction of inhomogeneity. As in \cite{califano1998ksw}, we trigger the instability by taking $\beta = 0.01, b = 0.001$.
In the numerical runs, we consider two cases: 
 \begin{align*}
 \textnormal{Run 1}:&\quad \delta = 0.5,\;v_{0,1}=v_{0,2}=0.3,\; \quad \;\;\; k_0=0.2; \quad \textrm{(initially symmetric beams)}\\
 \textnormal{Run 2}:&\quad \delta = \frac{1}{6},\;\;\;v_{0,1}=0.5,\;v_{0,2}=0.1,\; k_0=0.2. \quad \textrm{(initially nonsymmetric beams)}
 \end{align*}
 
In this problem, $f(x_2, v_1, v_2, 0)$  does not depend on $x_2$ and the initial particle
density is uniform and equals to  a constant, i.e.  $\rho(x_2,v_1,v_2, 0) =1$. We compute the solution on
the domain of $0\leq x_2 \leq L$, where $L=2 \pi/ k_0$, $k_0$
denotes the wave number. Periodic boundary conditions are assumed in
the $x_2$ direction.  The domain for $(v_1, v_2)$ is chosen such
that $f\simeq0$ on the boundaries. For the accuracy test, we set $\Ov=[-1.2,1.2]^2$. For other numerical results, we set $\Ov=[-1.5,1.5]^2$ to eliminate the
boundary effects and to accurately reflect the conservation properties of our methods.

 $\textnormal{\bf Scheme-1}(\Dt)$,  $\textnormal{\bf Scheme-2}(\Dt)$ are subject to CFL conditions.  While for the fully implicit method $\textnormal{\bf Scheme-5}(\Dt)$ , to save computational time, we use a fixed time step $\Delta t$. 
For  $\textnormal{\bf Scheme-5}$, we use KINSOL from SUNDIALS \cite{hindmarsh2005sundials} to solve the nonlinear algebraic systems resulting from the discretization of equation (b).

In all the runs below, we use the upwind flux for Vlasov solver. As discussed in \cite{cheng_va}, the central flux for Vlasov equation does not build 
any numerical dissipation into the scheme and this is not desired when filamentation occurs. 
For Maxwell solver,  as demonstrated in \cite{cheng_vm}, the upwind flux can cause energy dissipation. Therefore, in the scope of the current paper, we only consider the central flux \eqref{central:2} and  alternating flux \eqref{alternating:1} for the Maxwell solver. Those two flux choices will be extensively studied in the accuracy and conservation test with various time discretizations, and due to the superior performance of the alternating flux in terms of accuracy, we provide more simulation results with alternating fluxes in Figures \ref{figure_kineticeb}-\ref{figure_semiimplicitaltx2}. For simplicity, we use uniform meshes in $x_2$ ,$v_1$ and $v_2$ directions, while we note that nonuniform mesh can also be easily adapted under this DG framework.

\subsection{Accuracy tests}
In this subsection, we   test the orders of accuracy of the proposed schemes. The VM system is
time reversible, which provides a way to measure the errors of our schemes.  Let
$f(\bx,\bv,0),\;\bE(\bx,0),\;\bB(\bx,0)$
be the the initial conditions of the VM system and 
 $f(\bx,\bv,T)$, $\bE(\bx,T),\;\bB(\bx,T)$ be the solutions at $t=T$.
If we enforce
 $f(\bx,-\bv,T),\;\bE(\bx,T),\;-\bB(\bx,T)$ be the initial conditions for the VM system at $t=0$, then at $t=T$, we will
 recover
$f(\bx,-\bv,0),\;\bE(\bx,0),\;-\bB(\bx,0)$.  In Tables \ref{table_scheme1c} to \ref{table_scheme5a},  we run the VM system to $T=5$ and then back to $T=10$, and compare the numerical solution  with the exact initial conditions. The mesh is taken to be uniform with $N_{x_2}=N_{v_1}=N_{v_2}$.

%
%
%
 
 Tables 
 \ref{table_scheme1c} and \ref{table_scheme1a} list the $L^2$ errors and
orders for time integrator $\textnormal{\bf Scheme-1}$ with two flux choices for the Maxwell's equations:  the central flux and 
 the alternating flux.
  The parameters are those of Run 1 with symmetric counter-streaming.
 To match the accuracy of the temporal and spatial discretizations, we take   $\Delta t\sim \min(\Delta x, \Delta v)$ for space $\mG_h^{1}$, and $\Delta t\sim \min(\Delta x, \Delta v)^{3/2}$ for $\mG_h^{2}$, and $\Delta t\sim \min(\Delta x, \Delta v)^2$ for $\mG_h^{3}$.  Because of the stability restriction of the explicit scheme, we take $CFL = 0.3$ for $\mG_h^{1}$,  and the coefficient to be $0.15, 0.08$ for $\mG_h^{2}$ and $\mG_h^{3}$, respectively. From these tables, we can see that for all three polynomial spaces, we obtain the optimal $(k+1)$-th order for $f$, while the convergence order of $E_2$ is higher.
We also observe that  schemes with the upwind and alternating fluxes achieve optimal  order of $k+1$ for $B_3$ and $E_1$, while for odd $k$, the central flux gives suboptimal
order of accuracy for $B_3$ and $E_1$.

Tables 
  \ref{table_scheme2c} and \ref{table_scheme2a} list the $L^2$ errors and
orders for time integrator $\textnormal{\bf Scheme-2}$ with the central and alternating flux choices for the Maxwell's equations. We use the same parameter and CFL conditions as in $\textnormal{\bf Scheme-1}$. The conclusions for these tables are similar to $\textnormal{\bf Scheme-1}$.

 Tables \ref{table_scheme5c} and \ref{table_scheme5a}  list the $L^2$ errors and
orders for time integrator $\textnormal{\bf Scheme-5}$ with the central and alternating flux choices for the Maxwell's equation.   The parameters are those of Run 1 with symmetric counter-streaming.  The tolerance parameter in KINSOL solver is set to be $\epsilon_{tol}=10^{-12}$. $\Delta t$  is fixed to save computational time, and their values are listed in  Tables \ref{table_scheme5c} and \ref{table_scheme5a}. We observe the optimal $(k+1)$-th order for $f$, except for $\mS_h^{2}, \mW_h^2$. We believe this is because the mesh is still under-resolved to observe optimal order of convergence. This phenomenon is also present in \cite{cheng_va} for a similar type of mesh and space. For the  $E_1, E_2, B_3$ components, the order is sub-optimal for some mesh.
That's because the error in these components is so small, so that the tolerance parameter $\epsilon_{tol}$ for the Newton-Krylov solver has polluted the error in the calculation. 

Tables \ref{table_scheme5fc} and \ref{table_scheme5fa} list the $L^2$ errors and
orders for time integrator $\textnormal{\bf Scheme-5F}$ with the central and alternating flux choices for the Maxwell's equations. We have fixed $\Delta t$ to save computational time and use the same parameter setting as for $\textnormal{\bf Scheme-5}$.  Optimal convergence rate of fourth order in both space and time for f is observed, while the errors for $E_1,E_2,B_3$ has been polluted by the tolerance parameter $\epsilon_{tol}$ in the Newton-Krylov solver.
 \begin{table}[!htbp]
 \centering
\caption{Time discretization $\textnormal{\bf Scheme-1}$, with DG methods using the indicated polynomial space. Central flux for Maxwell's equation. $L^2$ error and order.}
\label{table_scheme1c}
\smallskip
\begin{tabular}{|c|c|c|c|c|c|c|}
\hline
\multicolumn{2}{|c|}{\multirow{2}*{\backslashbox{Space}{Mesh }}}&$20^3$&\multicolumn{2}{|c|}{$40^3$}&\multicolumn{2}{|c|}{$80^3$}\\
\cline{3-7}
 \multicolumn{2}{|c|}{} & Error& Error & Order& Error & Order\\
\hline
\multirow{4}*{$\mG_h^{1}, \mU_h^1$}
& $f$	&1.78E-01&	5.04E-02&	1.82&	1.30E-02&	1.95\\  \cline{2-7}
&$B_3$	&1.33E-05&	8.49E-06&	0.65&	5.04E-06&	0.75\\  \cline{2-7}
&$E_1$	&1.87E-06&	1.32E-06&	0.50&	5.85E-07&	1.17\\  \cline{2-7}
&$E_2$	&9.28E-07&	1.93E-07&	2.27&	2.05E-08&         3.23
\\
\hline
\multirow{4}*{$\mG_h^{2}, \mU_h^2$}
&$f$	&5.62E-02&	7.72E-03	&2.86	&1.02E-03&	2.92\\  \cline{2-7}
&$B_3$	&2.10E-07&	1.47E-08	&3.84	&1.14E-09&	3.69\\  \cline{2-7}
&$E_1$	&3.04E-08&	4.30E-09	&2.82	&1.98E-10&	4.44\\  \cline{2-7}
&$E_2$	&1.67E-07&	2.20E-08	&2.92	&1.47E-09&	3.90
\\
\hline
\multirow{4}*{$\mG_h^{3}, \mU_h^3$}
&$f$	&1.23E-02&	1.04E-03&	3.56&	7.01E-05&	3.89\\  \cline{2-7}
&$B_3$	&1.01E-07&	4.03E-09&	4.65&	1.12E-10&	5.17\\  \cline{2-7}
&$E_1$	&4.91E-08&	2.48E-10&	7.63&	3.04E-11&	3.03\\  \cline{2-7}
&$E_2$	&1.38E-08&	7.93E-10&	4.12&	1.58E-11&	5.65
\\ \hline
\end{tabular}
\end{table}
 \begin{table}[!htbp]
 \centering
\caption{Time discretization $\textnormal{\bf Scheme-1}$, with DG methods using the indicated polynomial space. Alternating flux for Maxwell's equation. $L^2$ error and order.}
\label{table_scheme1a}
\smallskip
\begin{tabular}{|c|c|c|c|c|c|c|}
\hline
\multicolumn{2}{|c|}{\multirow{2}*{\backslashbox{Space}{Mesh }}}&$20^3$&\multicolumn{2}{|c|}{$40^3$}&\multicolumn{2}{|c|}{$80^3$}\\
\cline{3-7}
 \multicolumn{2}{|c|}{} & Error& Error & Order& Error & Order\\
\hline
\multirow{4}*{$\mG_h^{1}, \mU_h^1$}
&$f$	&1.78E-01&	5.04E-02&	1.82&	1.30E-02&	1.95\\  \cline{2-7}
&$B_3$	&2.89E-06&	6.60E-07&	2.13&	1.66E-07&	1.99\\  \cline{2-7}
&$E_1$	&3.81E-07&	1.39E-07&	1.45&	3.46E-08&	2.01\\  \cline{2-7}
&$E_2$	&1.02E-06&	2.16E-07&	2.24&	2.21E-08&	3.29
\\
\hline
\multirow{4}*{$\mG_h^{2}, \mU_h^2$}
&$f$	&5.62E-02&	7.72E-03	&2.86	&1.02E-03&	2.92\\  \cline{2-7}
&$B_3$	&2.76E-07&	1.89E-08	&3.87	&1.19E-09&	3.99\\  \cline{2-7}
&$E_1$	&4.13E-08&	4.22E-09	&3.29	&1.97E-10&	4.42\\  \cline{2-7}
&$E_2$	&1.67E-07&	2.20E-08	&2.92	&1.47E-09&	3.90
\\
\hline
\multirow{4}*{$\mG_h^{3}, \mU_h^3$}
&$f$	&1.23E-02&	1.04E-03&	3.56&	7.01E-05&	3.89\\  \cline{2-7}
&$B_3$ &9.79E-08&	1.63E-09&	5.96&	4.04E-11&	5.33\\  \cline{2-7}
&$E_1$ &2.03E-08&	2.54E-10&	7.60&	4.28E-12&	5.89\\  \cline{2-7}
&$E_2$	&1.38E-08&	7.93E-10&	4.12&	1.59E-11&	5.64
\\ \hline
\end{tabular}
\end{table}

 \begin{table}[!htbp]
 \centering
\caption{Time discretization $\textnormal{\bf Scheme-2}$, with DG methods using the indicated polynomial space. Central flux for Maxwell's equation. $L^2$ error and order.}
\label{table_scheme2c}
\smallskip
\begin{tabular}{|c|c|c|c|c|c|c|}
\hline
\multicolumn{2}{|c|}{\multirow{2}*{\backslashbox{Space}{Mesh }}}&$20^3$&\multicolumn{2}{|c|}{$40^3$}&\multicolumn{2}{|c|}{$80^3$}\\
\cline{3-7}
 \multicolumn{2}{|c|}{} & Error& Error & Order& Error & Order\\
\hline
\multirow{4}*{$\mG_h^{1}, \mU_h^1$}
&f	&1.78E-01&	5.04E-02&	1.82&	1.30E-02&	1.95\\  \cline{2-7}
&B3	&1.34E-05&	8.50E-06&	0.66&	5.04E-06&	0.75\\  \cline{2-7}
&E1	&1.85E-06&	1.31E-06&	0.50&	5.84E-07&	1.17\\  \cline{2-7}
&E2	&9.28E-07&	1.93E-07&	2.27&	2.05E-08&	3.23
\\
\hline
\multirow{4}*{$\mG_h^{2}, \mU_h^2$}
&f	&5.62E-02&	7.72E-03	&2.86	&1.02E-03&	2.92\\  \cline{2-7}
&B3	&2.21E-07&	1.56E-08	&3.82	&1.15E-09&	3.76\\  \cline{2-7}
&E1	&2.30E-08&	4.03E-09	&2.51	&1.90E-10&	4.41\\  \cline{2-7}
&E2	&1.67E-07&	2.20E-08	&2.92	&1.47E-09&	3.90
\\
\hline
\multirow{4}*{$\mG_h^{3}, \mU_h^3$}
&f	&1.23E-02&	1.04E-03&	3.56&	7.01E-05&	3.89\\  \cline{2-7}
&B3	&1.01E-07&	4.03E-09&	4.65&	1.12E-10&	5.17\\  \cline{2-7}
&E1	&4.92E-08&	5.33E-10&	6.53&	3.40E-11&	3.97\\  \cline{2-7}
&E2	&1.38E-08&	7.93E-10&	4.12&	1.58E-11&	5.65
\\ \hline
\end{tabular}
\end{table}
 \begin{table}[!htbp]
 \centering
\caption{Time discretization $\textnormal{\bf Scheme-2}$, with DG methods using the indicated polynomial space. Alternating flux for Maxwell's equation. $L^2$ error and order.}
\label{table_scheme2a}
\smallskip
\begin{tabular}{|c|c|c|c|c|c|c|}
\hline
\multicolumn{2}{|c|}{\multirow{2}*{\backslashbox{Space}{Mesh }}}&$20^3$&\multicolumn{2}{|c|}{$40^3$}&\multicolumn{2}{|c|}{$80^3$}\\
\cline{3-7}
 \multicolumn{2}{|c|}{} & Error& Error & Order& Error & Order\\
\hline
\multirow{4}*{$\mG_h^{1}, \mU_h^1$}
&f	&1.78E-01&	5.04E-02&	1.82&	1.30E-02&	1.95\\  \cline{2-7}
&B3	&2.99E-06&	6.55E-07&	2.19&	1.64E-07&	2.00\\  \cline{2-7}
&E1	&1.34E-07&	3.64E-08&	1.88&	1.06E-08&	1.78\\  \cline{2-7}
&E2	&1.02E-06&	2.16E-07&	2.24&	2.21E-08&	3.29
\\
\hline
\multirow{4}*{$\mG_h^{2}, \mU_h^2$}
&f	&5.62E-02&	7.72E-03	&2.86	&1.02E-03&	2.92\\  \cline{2-7}
&B3	&2.71E-07&	1.90E-08	&3.83	&1.19E-09&	4.00\\  \cline{2-7}
&E1	&3.67E-08&	4.06E-09	&3.18	&1.90E-10&	4.42\\  \cline{2-7}
&E2	&1.67E-07&	2.20E-08	&2.92	&1.47E-09&	3.90
\\
\hline
\multirow{4}*{$\mG_h^{3}, \mU_h^3$}
&f	&1.23E-02	&1.04E-03	&3.56&	7.01E-05&	3.89\\  \cline{2-7}
&B3	&9.77E-08	&1.63E-09	&5.91&	4.04E-11&	5.33\\  \cline{2-7}
&E1	&1.99E-08	&5.37E-10	&5.21&	1.58E-11&	5.09\\  \cline{2-7}
&E2	&1.38E-08	&7.93E-10	&4.12&	1.59E-11&	5.64
\\ \hline
\end{tabular}
\end{table}
 \begin{table}[!htbp]
 \centering
\caption{Time discretization $\textnormal{\bf Scheme-5}$, with DG methods using the indicated polynomial space. Central flux for Maxwell's equation. $L^2$ error and order.}
\label{table_scheme5c}
\smallskip
\begin{tabular}{|c|c|c|c|c|c|c|}
\hline
\multicolumn{1}{|c|}{Space}& & Error& Error & Order& Error & Order\\
\hline
\multirow{3}*{$\mS_h^{1}, \mW_h^1$}& Mesh&$40^3$&\multicolumn{2}{|c|}{$60^3$}&\multicolumn{2}{|c|}{$80^3$}\\
\cline{2-7}
&$f$	&2.74E-03&	1.36E-03	&1.73	&8.16E-04&	1.78\\ \cline{2-7}
&$B_3$	&7.83E-06&	5.93E-06	&0.69	&4.77E-06&	0.76\\ \cline{2-7}
\multirow{2}*{$\Delta t=0.1 \frac{40}{N_x}$}&$E_1$	&1.38E-06&	8.88E-07	&1.09	&6.15E-07&	1.28\\ \cline{2-7}
&$E_2$	&1.85E-08&	4.17E-09	&3.67	&2.86E-09&	1.31
\\
\hline
\multirow{3}*{$\mS_h^{2}, \mW_h^2$}& Mesh&$80^3$&\multicolumn{2}{|c|}{$100^3$}&\multicolumn{2}{|c|}{$120^3$}\\
\cline{2-7}
&$f$	         &4.27E-04	&2.30E-04	&2.15	&1.41E-04&	2.19\\ \cline{2-7}
&$B_3$	&4.33E-08	&2.14E-08	&2.45	&1.31E-08&	2.20\\ \cline{2-7}
\multirow{2}*{$\Delta t=0.2\left( \frac{40}{N_x}\right)^{3/2}$}
&$E_1$	&3.38E-08	&1.64E-08	&2.51	&9.69E-09&	2.36\\ \cline{2-7}
&$E_2$	&2.68E-10	&8.63E-11	&3.94	&4.53E-11&	2.89\\
\hline

\multirow{2}*{$\mS_h^{3}, \mW_h^3$}& Mesh&$40^3$&\multicolumn{2}{|c|}{$60^3$}&\multicolumn{2}{|c|}{$80^3$}\\
\cline{2-7}

&$f$	&1.00E-04	&2.18E-05&	3.76& 7.65E-06&	3.64\\ \cline{2-7}
&$B_3$&	1.54E-07	&6.08E-08&	2.29& 2.99E-08&	2.47 \\ \cline{2-7}
\multirow{2}*{{$\Delta t=0.2 (\frac{40}{N_x})^{2}$}}&$E_1$&	2.75E-08	&2.93E-08&	-0.16& 2.11E-08&	1.14\\ \cline{2-7}
&$E_2$&	5.71E-10	&2.38E-10&	2.16&1.43E-10	&1.77 
\\ \hline
\end{tabular}
\end{table}

 \begin{table}[!htbp]
 \centering
\caption{Time discretization $\textnormal{\bf Scheme-5}$, with DG methods using the indicated polynomial space. Alternating flux for Maxwell's equation. $L^2$ error and order.}
\label{table_scheme5a}
\smallskip
\begin{tabular}{|c|c|c|c|c|c|c|}
\hline
\multicolumn{1}{|c|}{Space} &  &Error& Error & Order& Error & Order\\
\hline
\multirow{3}*{$\mS_h^{1}, \mW_h^1$}&Mesh &$40^3$&\multicolumn{2}{|c|}{$60^3$}&\multicolumn{2}{|c|}{$80^3$}\\
\cline{2-7}

&$f$	&2.74E-03	&1.36E-03	&1.73	&8.15E-04&	1.78 \\ \cline{2-7}
&$B_3$&	1.90E-07	&6.13E-08	&2.79	&2.33E-08&	3.36 \\ \cline{2-7}
\multirow{2}*{$\Delta t=0.1 \frac{40}{N_x}$}&$E_1$	&5.73E-08	&4.62E-09	&6.21	&4.60E-09	&0.02 \\ \cline{2-7}
&$E_2$	&1.98E-08	&2.83E-09	&4.80	&8.40E-10	&4.22\\
\hline
\multirow{3}*{$\mS_h^{2}, \mW_h^2$}&Mesh &$60^3$&\multicolumn{2}{|c|}{$80^3$}&\multicolumn{2}{|c|}{$100^3$}\\
\cline{2-7}
&$f$	&4.27E-04&	2.30E-04&	2.15	&1.41E-04	&2.19\\ \cline{2-7}
&$B_3$&	3.83E-08&	1.58E-08	&3.08&	7.88E-09&	3.12\\ \cline{2-7}
\multirow{2}*{$\Delta t=0.2\left( \frac{40}{N_x}\right)^{3/2}$}
&$E_1$&	1.59E-09&	6.27E-10&	3.23&	2.58E-107&	3.98\\ \cline{2-7}
&$E_2$&	2.59E-10	&8.15E-11	&4.02	&4.29E-11&	2.88
\\
\hline
\multirow{2}*{$\mS_h^{3}, \mW_h^3$}&Mesh &$40^3$&\multicolumn{2}{|c|}{$60^3$}&\multicolumn{2}{|c|}{$80^3$}\\
\cline{2-7}

&$f$	&1.00E-04&	2.18E-05	&3.76&	7.62E-06	&3.65\\ \cline{2-7}
&$B_3$	&1.22E-07&	4.41E-08&	2.51	&2.02E-08&	2.71\\ \cline{2-7}
\multirow{2}*{$\Delta t=0.2 (\frac{40}{N_x})^2$}&$E_1$	&1.58E-08&	3.97E-09	&3.41&	1.30E-09	&3.88\\ \cline{2-7}
&$E_2$	&5.35E-10&	2.31E-10&	2.07	&1.41E-10&	1.72
\\ \hline
\end{tabular}
\end{table}

 \begin{table}[!htbp]
 \centering
\caption{Time discretization $\textnormal{\bf Scheme-5F}$, with DG methods using the indicated polynomial space. 
Central flux for Maxwell's equation. $L^2$ error and order. $\Delta t=0.2 \frac{20}{N_x}$.}
\label{table_scheme5fc}
\smallskip
\begin{tabular}{|c|c|c|c|c|c|c|}
\hline
\multicolumn{2}{|c|}{\multirow{2}*{\backslashbox{Space}{Mesh }}}&$20^3$&\multicolumn{2}{|c|}{$40^3$}&\multicolumn{2}{|c|}{$60^3$}\\
\cline{3-7}
 \multicolumn{2}{|c|}{} & Error& Error & Order& Error & Order\\
\hline
\multirow{4}*{$\mS_h^{3}, \mW_h^3$}
&$f$	&3.05E-03&	2.30E-04&	3.73	&3.67E-05&	4.53\\ \cline{2-7}
&$B_3$	&4.65E-07&	1.15E-07&	2.02&	4.13E-08&	2.53\\ \cline{2-7}
&$E_1$	&4.49E-07 &	2.87E-08 &	3.97	&3.17E-08&	-0.25\\ \cline{2-7}
&$E_2$	&3.12E-09 &	5.95E-10	&2.39&	2.18E-10	&2.48\\ \hline
\end{tabular}
\end{table}

 \begin{table}[!htbp]
 \centering
\caption{Time discretization $\textnormal{\bf Scheme-5F}$, with DG methods using the indicated polynomial space. 
Alternating flux for Maxwell's equation. $L^2$ error and order. $\Delta t=0.2 \frac{20}{N_x}.$}
\label{table_scheme5fa}
\smallskip
\begin{tabular}{|c|c|c|c|c|c|c|}
\hline
\multicolumn{2}{|c|}{\multirow{2}*{\backslashbox{Space}{Mesh }}}&$20^3$&\multicolumn{2}{|c|}{$40^3$}&\multicolumn{2}{|c|}{$60^3$}\\
\cline{3-7}
 \multicolumn{2}{|c|}{} & Error& Error & Order& Error & Order\\
\hline
\multirow{4}*{$\mS_h^{3}, \mW_h^3$}
&$f$	&3.05E-03&	2.30E-04&	3.73&	3.67E-05&	4.53\\ \cline{2-7}
&$B_3$	&4.95E-07&	1.27E-07&	1.96&	4.45E-08&	2.59\\ \cline{2-7}
&$E_1$	&1.44E-07&	1.36E-08&	3.40&	4.98E-09&	2.48\\ \cline{2-7}
&$E_2$	&3.18E-09&	5.86E-10&	2.44&	2.15E-10&	2.47\\ \hline
\end{tabular}
\end{table}

\subsection{Conservation properties}
In this subsection, we will verify the conservation properties of the proposed methods $\textnormal{\bf Scheme-1}$, $\textnormal{\bf Scheme-2}$ and $\textnormal{\bf Scheme-5}$. In particular, we test $\textnormal{\bf Scheme-1}$ and $\textnormal{\bf Scheme-2}$ with space $\mG_h^2$, and denote them by ``$\textnormal{\bf Scheme-1}\textnormal{ and } P^2$" and ``$\textnormal{\bf Scheme-2}\textnormal{ and } P^2$" in the figures, respectively. We run these two schemes on $100^3$ mesh with $CFL = 0.15$. We test $\textnormal{\bf Scheme-5}$ with space $\mS_h^2$ and denote it by ``$\textnormal{\bf Scheme-5}\textnormal{ and } Q^2$''. To save computational time for this fully implicit scheme, we set $\epsilon_{tol}=10^{-8}$ for the Newton-Krylov solver and run this scheme on $80^3$ mesh with fixed time step $\Delta t = 0.2$. 

\begin{figure}[!htbp]
\centering
\subfigure[Total particle number. Run1 ]{\includegraphics[width=0.45\textwidth]{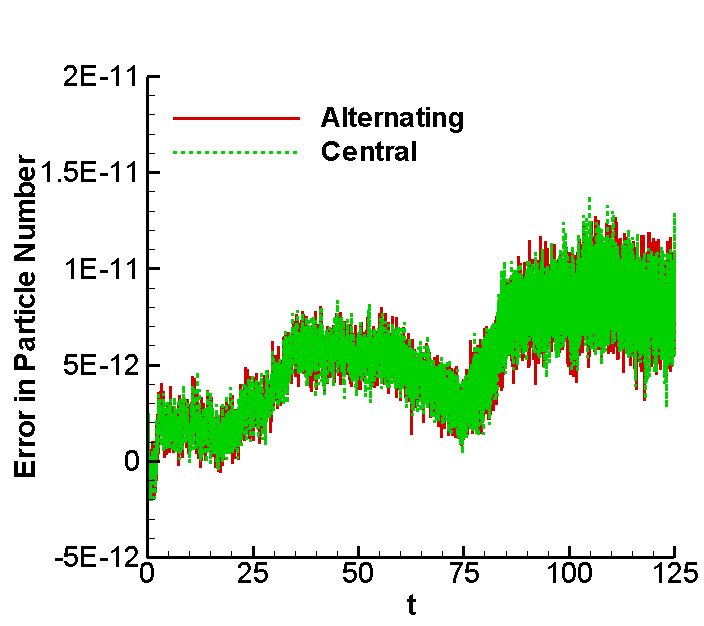}}
\subfigure[Total particle number. Run2 ]{\includegraphics[width=0.45\textwidth]{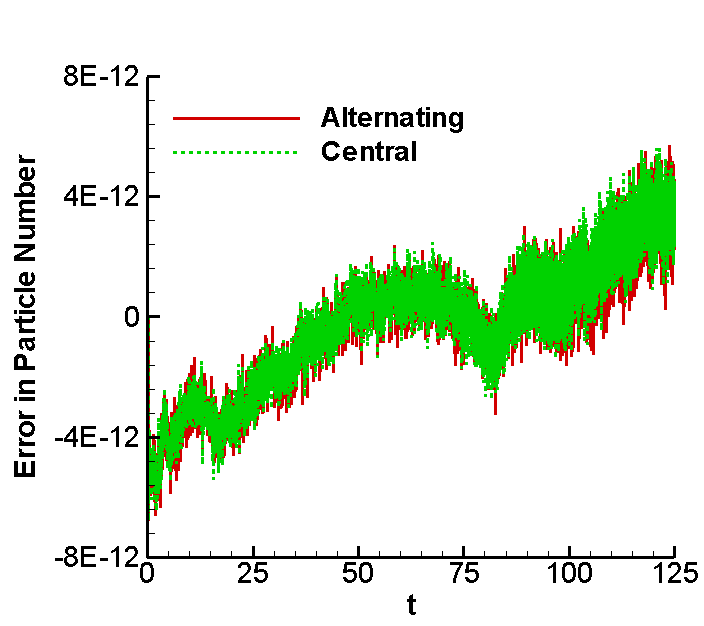}}\\
\subfigure[Total energy. Run1 ]{\includegraphics[width=0.45\textwidth]{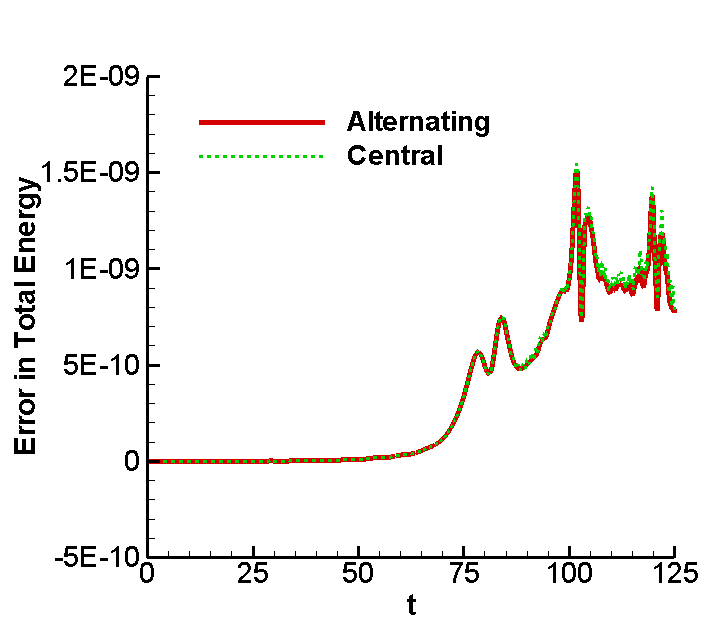}}
\subfigure[Total energy. Run2 ]{\includegraphics[width=0.45\textwidth]{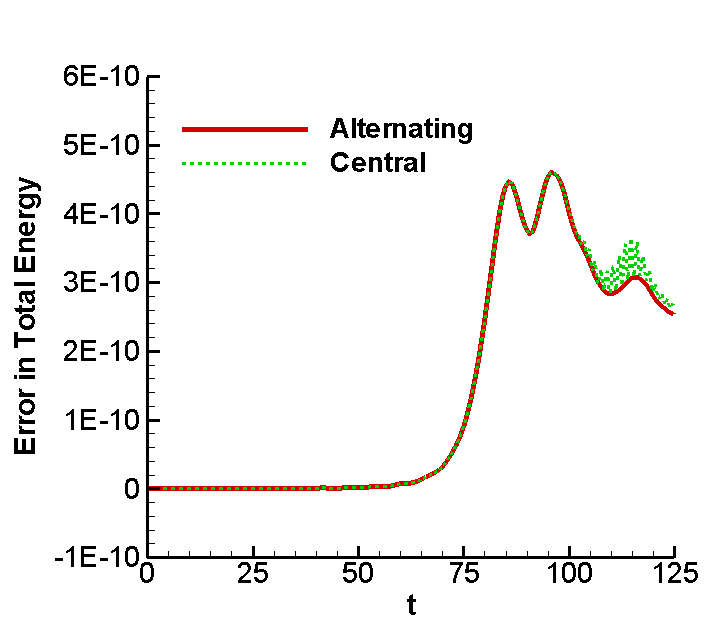}}
\caption{Evolution of the error in total particle number and total energy computed by $\textnormal{\bf Scheme-1}\textnormal{ and } P^2$  with indicated fluxes.  $100^3$ mesh.  $CFL = 0.15$.}
\label{figure_explicit}
\end{figure}

\begin{figure}[!htbp]
\centering
\subfigure[Total particle number. Run1 ]{\includegraphics[width=0.45\textwidth]{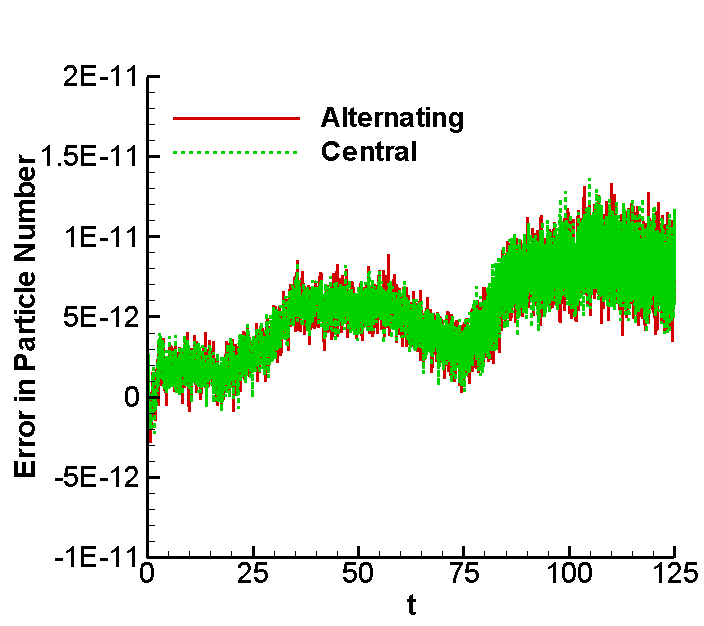}}
\subfigure[Total particle number. Run2 ]{\includegraphics[width=0.45\textwidth]{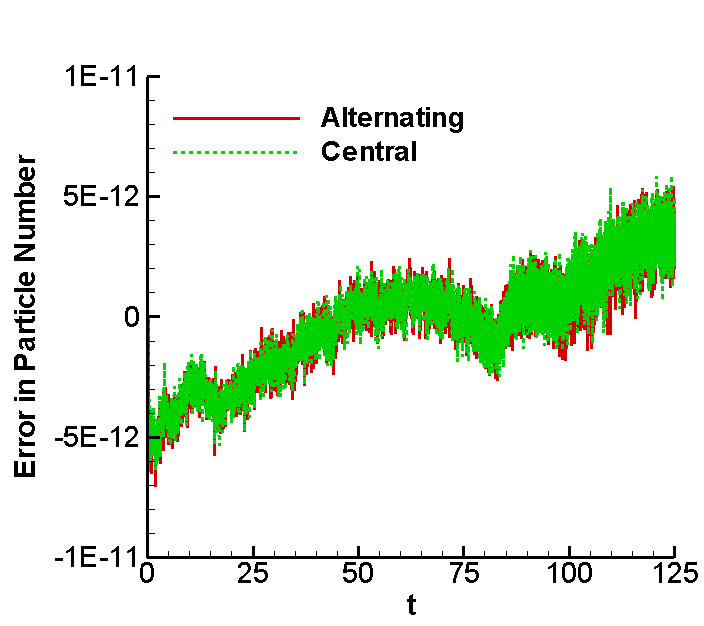}}\\
\subfigure[Total energy. Run1 ]{\includegraphics[width=0.45\textwidth]{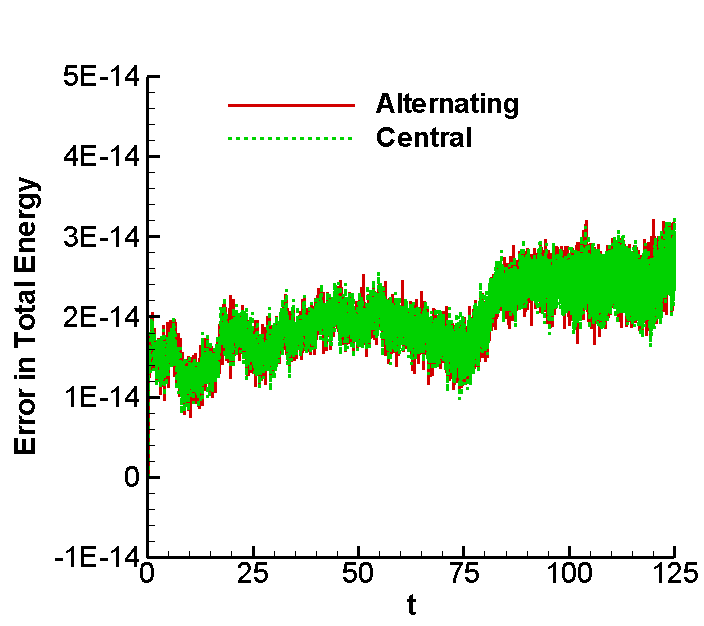}}
\subfigure[Total energy. Run2 ]{\includegraphics[width=0.45\textwidth]{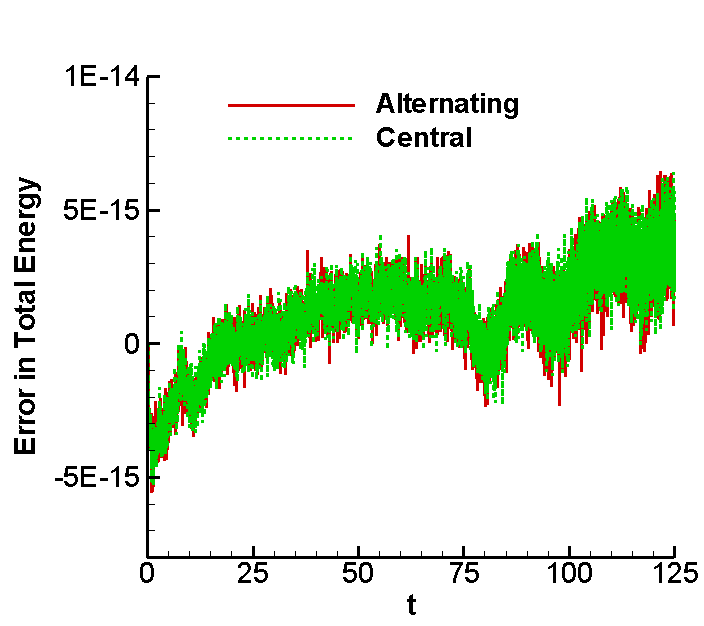}}\\
\caption{Evolution of the error in total particle number and total energy computed by $\textnormal{\bf Scheme-2}\textnormal{ and } P^2$  with indicated fluxes.  $100^3$ mesh.  $CFL = 0.15$.}
\label{figure_semiimplicit}
\end{figure}

In Figure \ref{figure_explicit}, we plot the error of the total particle number and total energy of parameter choice of Run 1 and Run 2 for $\textnormal{\bf Scheme-1}\textnormal{ and } P^2$ with the central and alternating fluxes for the Maxwell solver.  We  observe that the errors   stay small, below $10^{-11}$ for  total particle number, $10^{-9}$ for total energy with parameter Run 1, and $10^{-10}$ for total energy with parameter Run 2.
 The conservation is especially good with  $\textnormal{\bf Scheme-2}$ as in Figure \ref{figure_semiimplicit}.  The errors  are below $10^{-11}$ for the total particle number and below $10^{-14}$ for the total energy.  The difference of  the energy conservation in $\textnormal{\bf Scheme-1}$ and $\textnormal{\bf Scheme-2}$ reflects the exact conservation of   $\textnormal{\bf Scheme-2}$ and near conservation of $\textnormal{\bf Scheme-1}$  as illustrated in  Theorem \ref{thm:energy}.
  
\begin{figure}[!htbp]
\centering
\subfigure[Total particle number. Run1 ]{\includegraphics[width=0.45\textwidth]{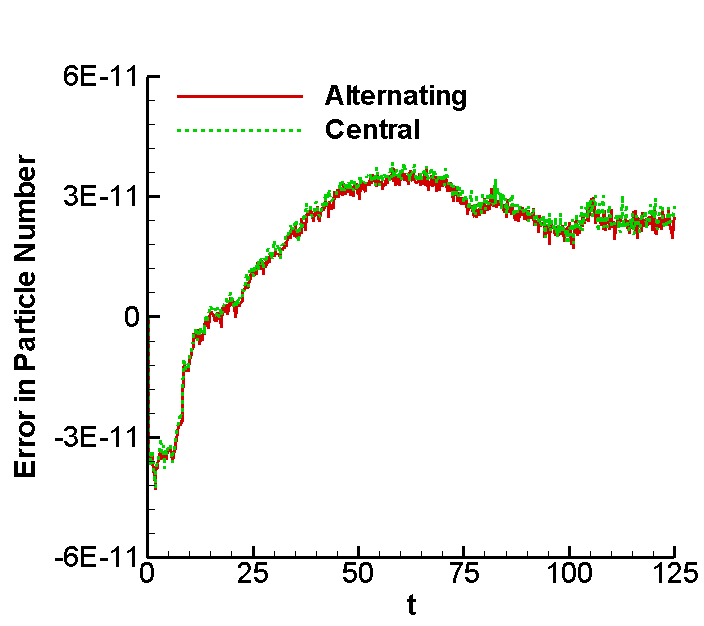}}
\subfigure[Total particle number. Run2 ]{\includegraphics[width=0.45\textwidth]{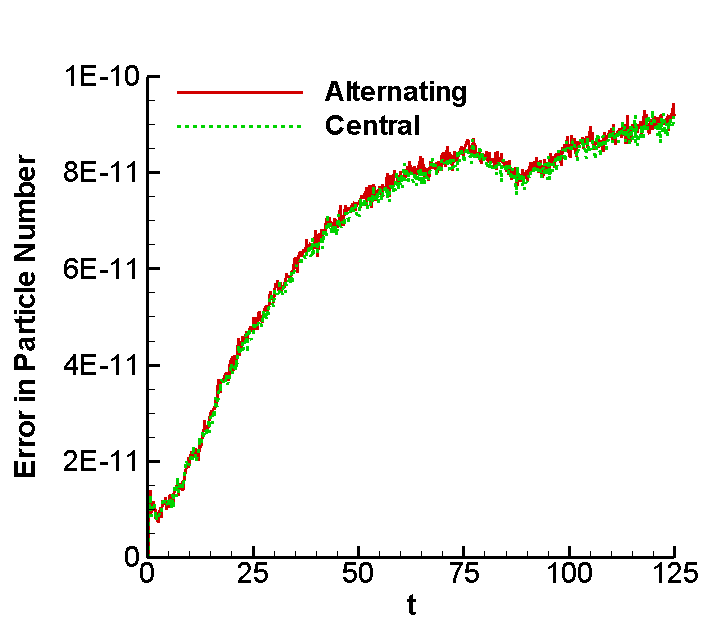}}\\
\subfigure[Total energy. Run1 ]{\includegraphics[width=0.45\textwidth]{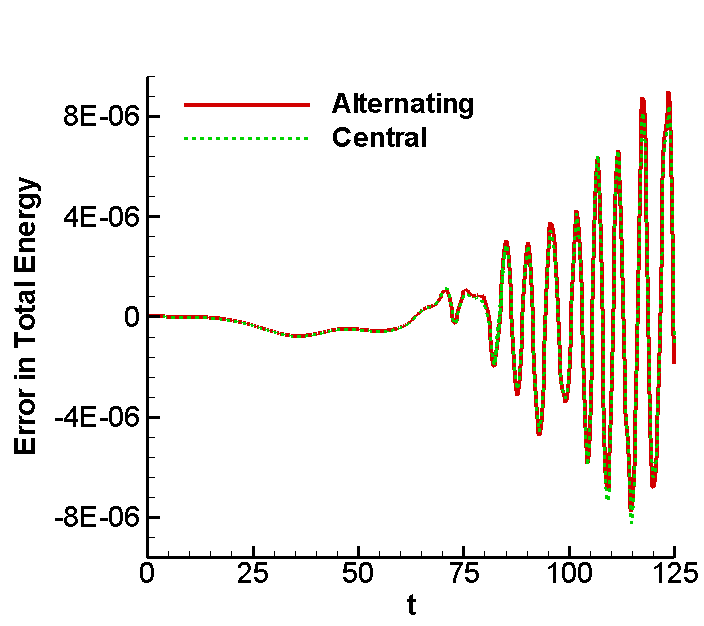}}
\subfigure[Total energy. Run2 ]{\includegraphics[width=0.45\textwidth]{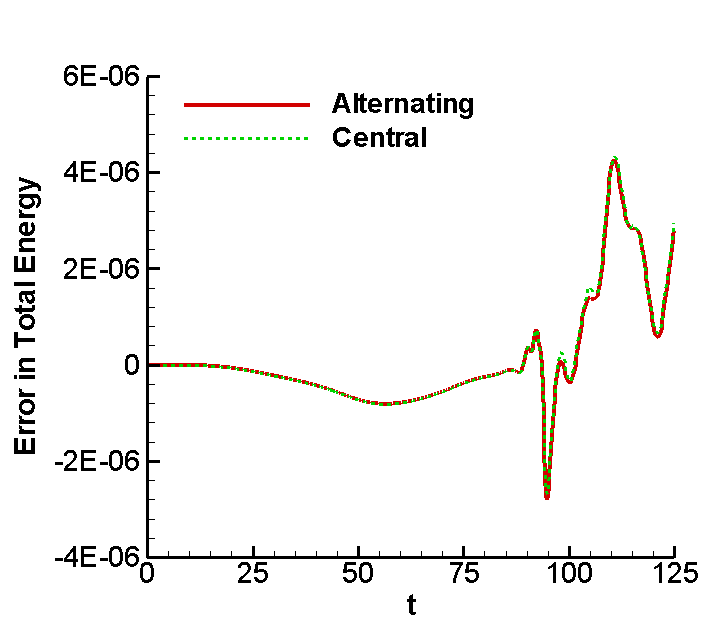}}\\
\caption{Evolution of the error in total particle number and total energy computed by $\textnormal{\bf Scheme-5}\textnormal{ and } Q^2$  with indicated fluxes.  $80^3$ mesh.  $\Delta t = 0.2$.}
\label{figure_implicit}
\end{figure}

 Figure \ref{figure_implicit} shows the error of the total particle number and total energy for   $\textnormal{\bf Scheme-5}$. The errors for mass are below $10^{-11}$. The errors for the total energy are below $10^{-6}$ and they are larger mainly due to the error caused by  the Newton-Krylov solver, and is related to $\epsilon_{tol}=10^{-8}$. All these results agree well with the theorems in the previous section.
 
\begin{figure}[!htbp]
\centering
\subfigure[Total particle number. Run1 ]{\includegraphics[width=0.45\textwidth]{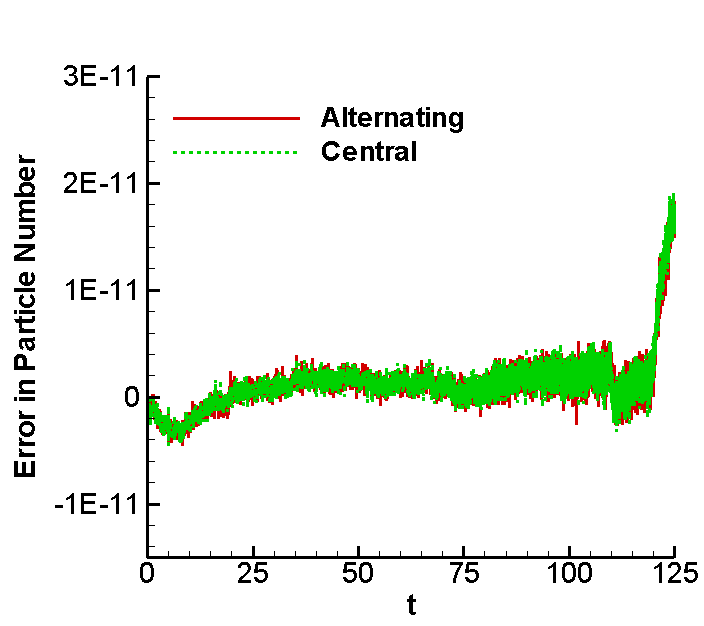}}
\subfigure[Total particle number. Run2 ]{\includegraphics[width=0.45\textwidth]{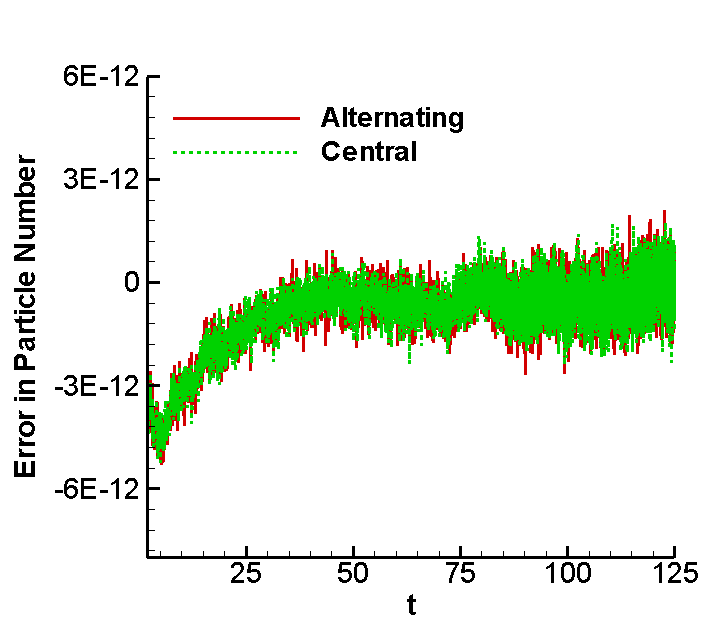}}\\
\subfigure[Total energy. Run1 ]{\includegraphics[width=0.45\textwidth]{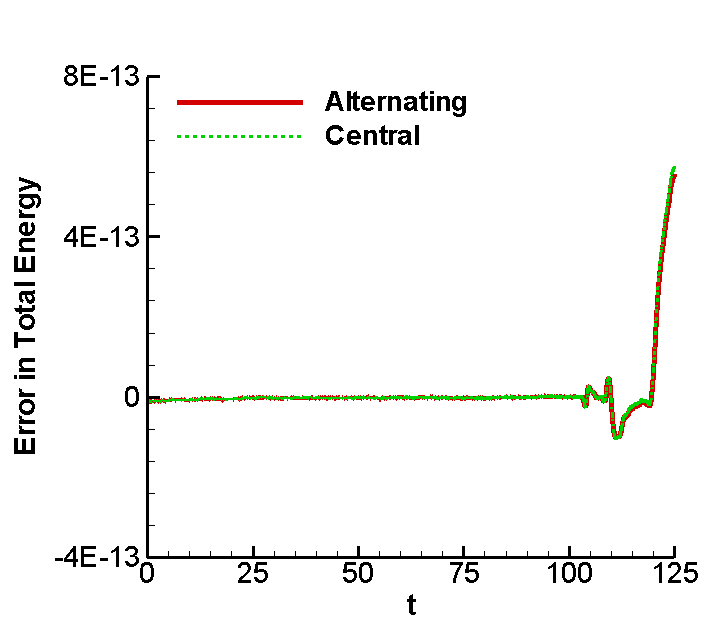}}
\subfigure[Total energy. Run2 ]{\includegraphics[width=0.45\textwidth]{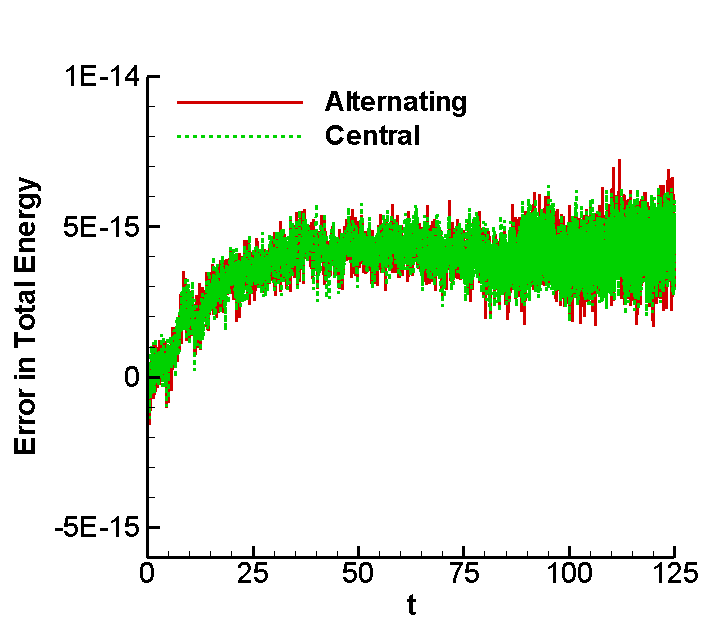}}\\
\caption{Evolution of the error in total particle number and total energy computed by $\textnormal{\bf Scheme-2}\textnormal{ and } P^2$  with indicated fluxes.  $50^3$ mesh.  $CFL = 0.15$.}
\label{figure_semiimplicitmesh50}
\end{figure}

 In Figure \ref{figure_semiimplicitmesh50}, we use a coarse mesh ($50^3$) to plot the errors in the conserved quantities to demonstrate that the conservation properties of our schemes are mesh independent. We use  $\textnormal{\bf Scheme-2} \textnormal{ and } P^2$ to demonstrate the behavior. Upon comparison with the results from finer mesh  in Figures \ref{figure_semiimplicit}, we conclude that the mesh size has no impact on the conservation of total particle number and total energy as predicted by  Theorems \ref{thm:mass} and \ref{thm:fullyenergy2}. This demonstrates the distinctive feature of our scheme: the total particle number and energy can be well preserved even with an under-resolved mesh.

\subsection{Collections of numerical data}
In this subsection, we collect some sample numerical data to benchmark our schemes. The results are computed by  $\textnormal{\bf Scheme-2}\textnormal{ and } P^2$ on a $80^3$ mesh with alternating flux for the Maxwell solver.

\begin{figure}[!htbp]
\centering
\subfigure[Kinetic energy. Run1]{\includegraphics[width=0.45\textwidth]{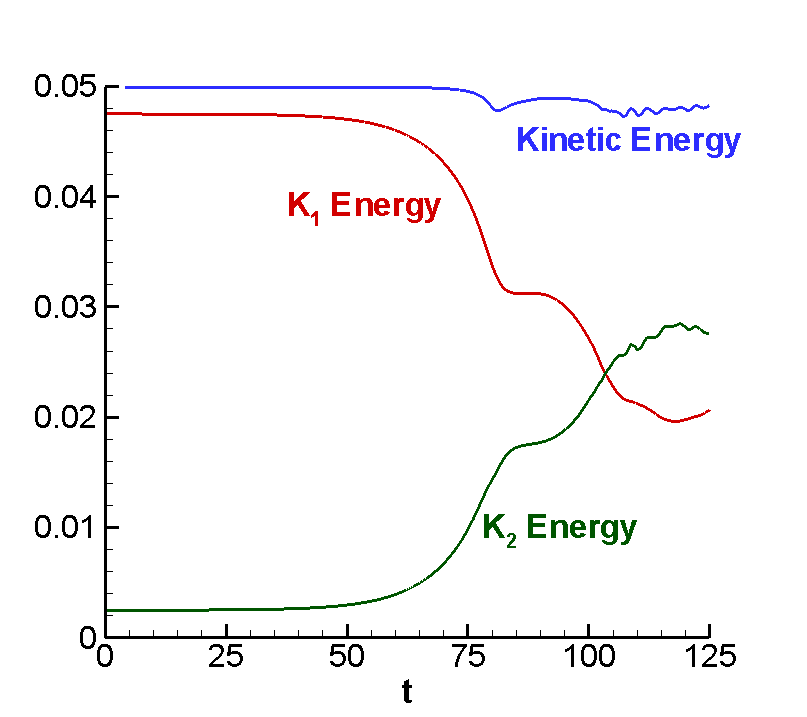}}
\subfigure[Kinetic energy. Run2]{\includegraphics[width=0.45\textwidth]{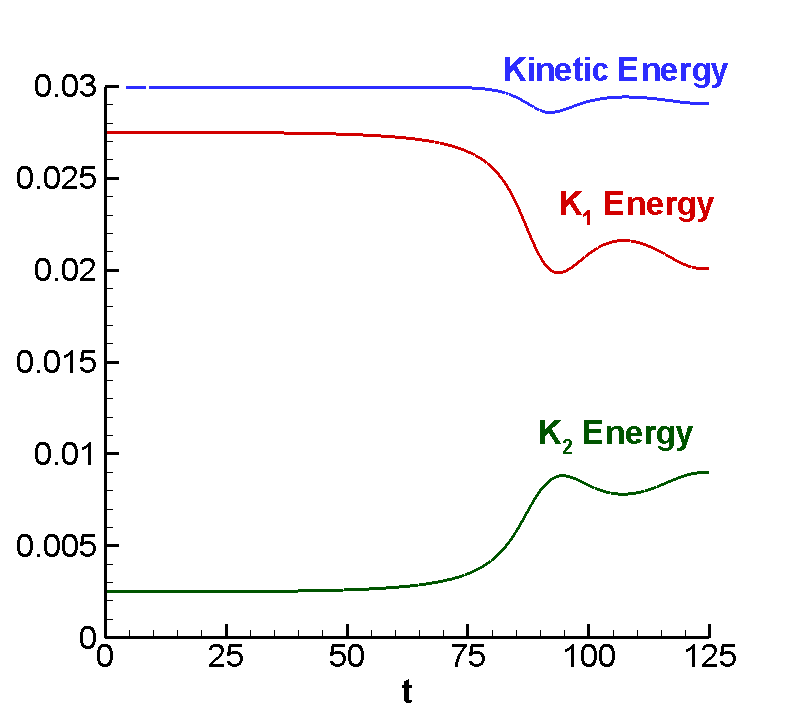}}\\
\subfigure[Electric and magnetic energy. Run1]{\includegraphics[width=0.45\textwidth]{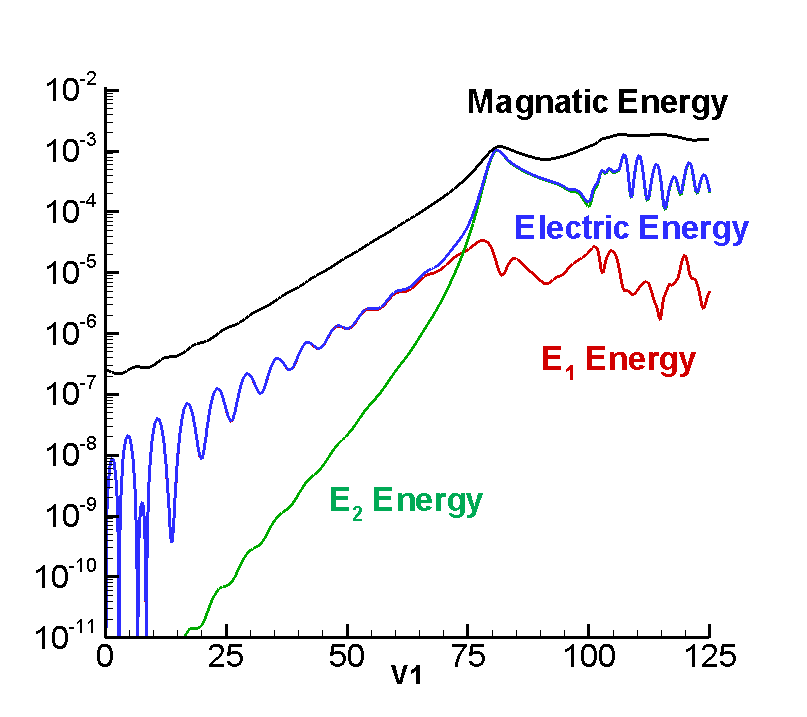}}
\subfigure[Electric and magnetic energy. Run2]{\includegraphics[width=0.45\textwidth]{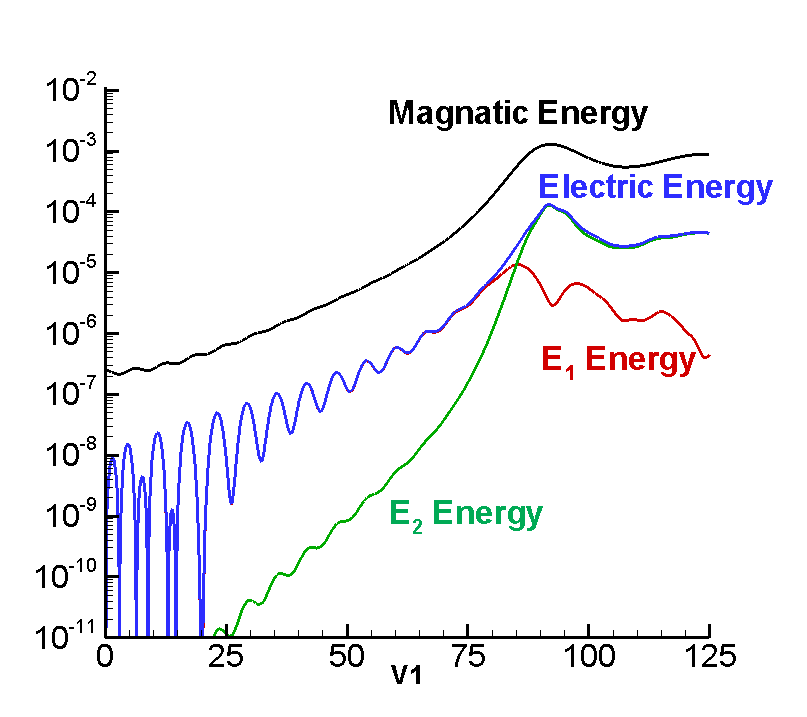}}
\caption{Evolution of the kinetic, electric and magnetic energies for the streaming Weibel instability. $\textnormal{\bf Scheme-2}\textnormal{ and } P^2$.  
$80^3$ mesh.
 $CFL = 0.15$. Alternating flux for Maxwell solver. }
 \label{figure_kineticeb}
\end{figure}

\begin{figure}[!htbp]
\centering
\subfigure[$t=55$. Run1]{\includegraphics[width=0.45\textwidth]{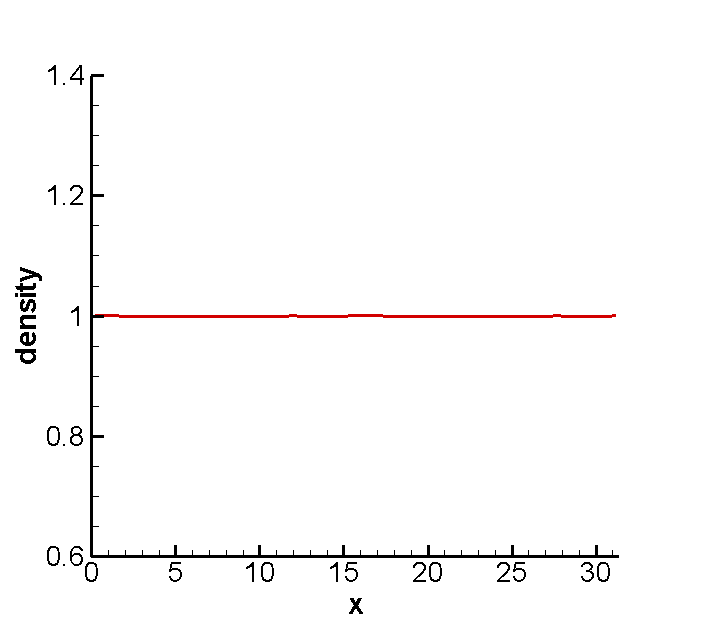}}
\subfigure[$t=55$. Run2]{\includegraphics[width=0.45\textwidth]{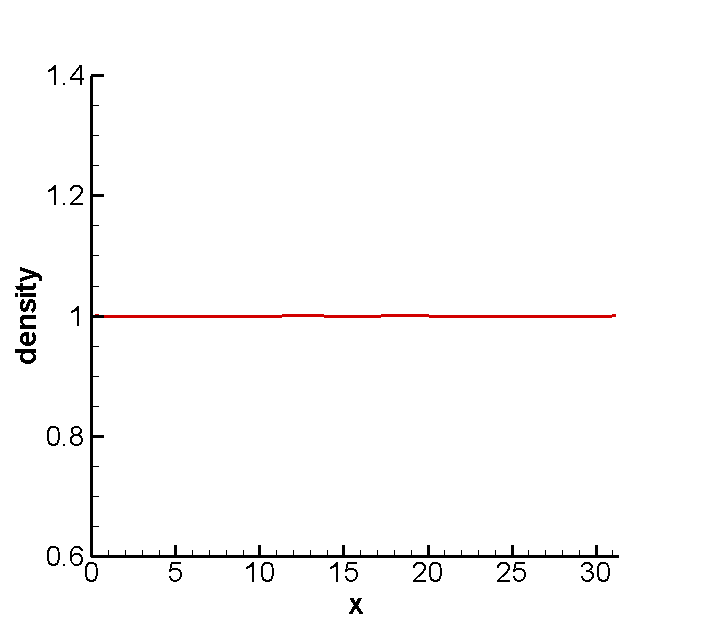}}\\
\subfigure[$t=82$. Run 1 ]{\includegraphics[width=0.45\textwidth]{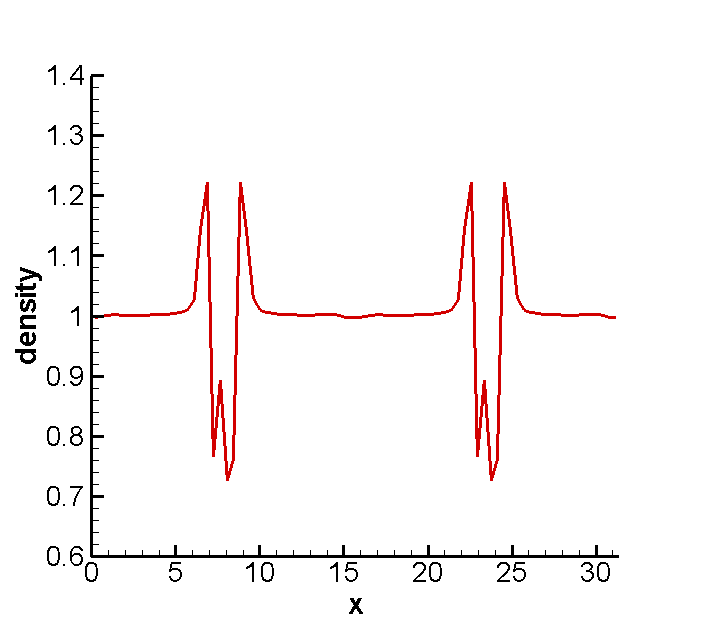}}
\subfigure[$t=82$. Run 2 ]{\includegraphics[width=0.45\textwidth]{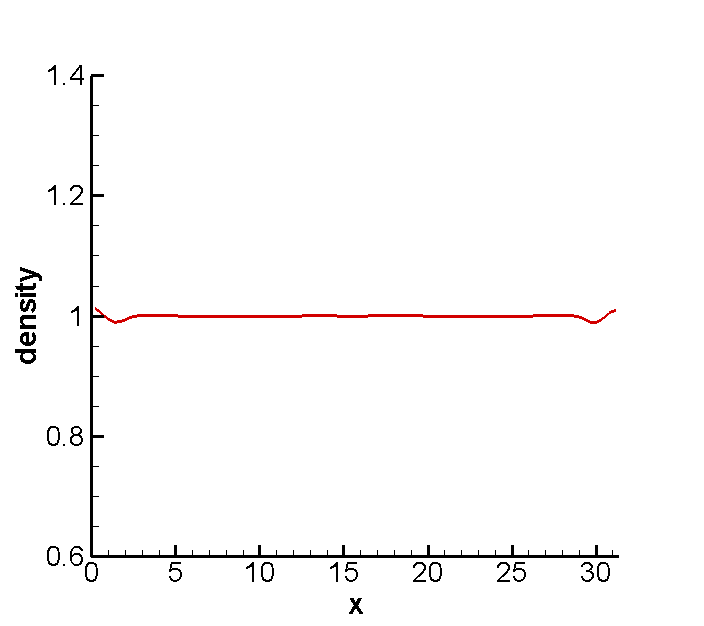}}\\
\subfigure[$ t=125$. Run 1]{\includegraphics[width=0.45\textwidth]{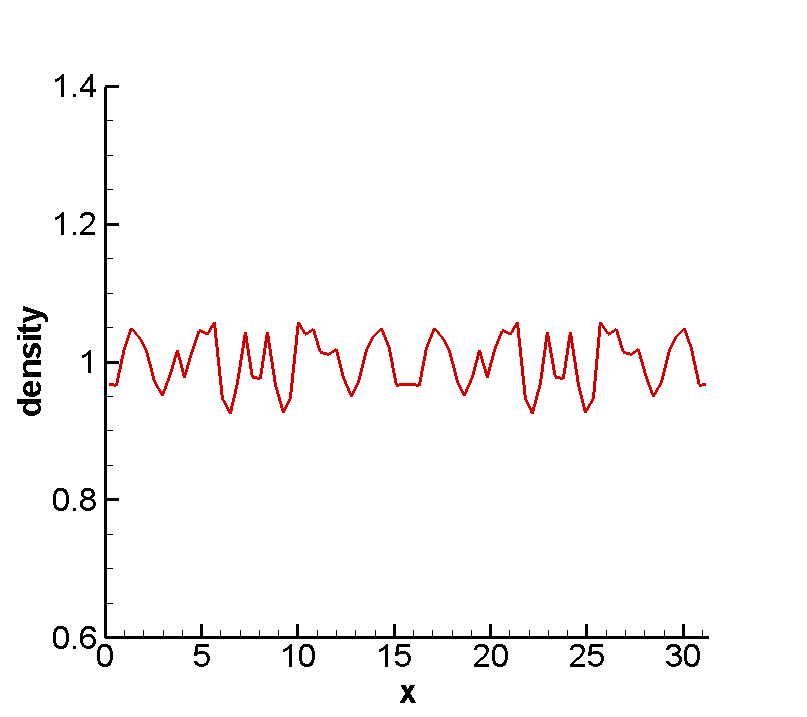}}
\subfigure[$ t=125$. Run 2]{\includegraphics[width=0.45\textwidth]{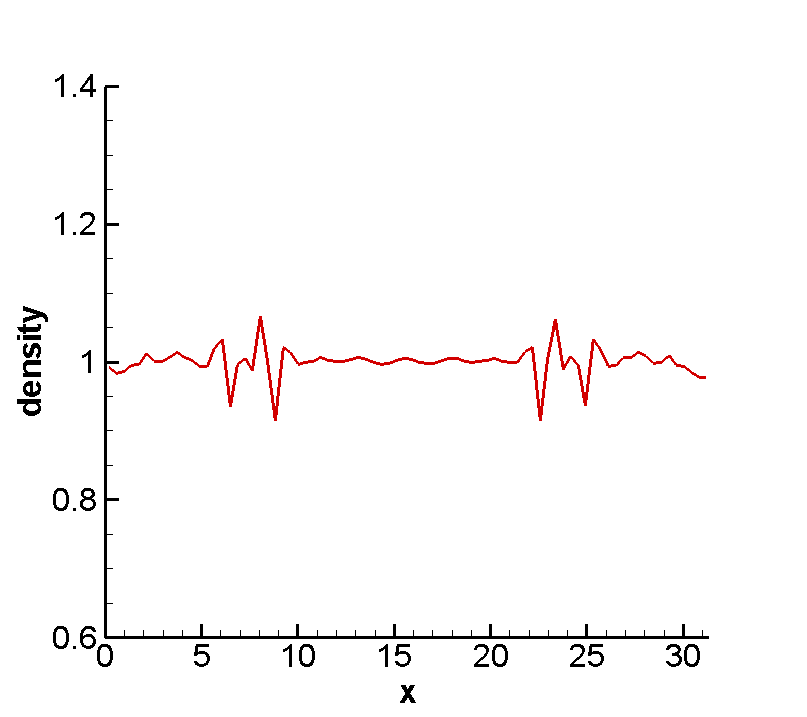}}
\caption{Plots of the density function for the streaming Weibel instability at selected  time $t$. $\textnormal{\bf Scheme-2}\textnormal{ and } P^2$. 
 $80^3$ mesh.
 $CFL = 0.15$. Alternating flux for Maxwell solver. }
\label{figure_semiimplicitaltdens}
\end{figure}

Figure \ref{figure_kineticeb} plots  the time evolution of the kinetic, electric and magnetic energies with parameter choice of Run 1 and Run 2. In particular, we plot the separate components of the kinetic energy, which are defined by
$
K_1=\frac{1}{2L}\int_{0}^{L}\int_{\Omega_v}fv_1^2dv_1dv_2dx_2, \; K_2=\frac{1}{2L}\int_{0}^{L}\int_{\Omega_v}fv_2^2dv_1dv_2dx_2,$
and the separate components of electric energy with $E_1$ energy and $E_2$ energy defined by
$\frac{1}{2L}\int_{0}^{L}E_1^2dx_2$ and $\frac{1}{2L}\int_{0}^{L}E_2^2dx_2,
$ respectively.  Figure (a) and (b) show the transference
of kinetic energy from one component to the other with a deficit converted into field energy, which is consistent with the total energy conservation, as shown in Figure \ref{figure_semiimplicit}.  After a rapid transient, the magnetic and inductive electric fields grow initially at a linear growth rate. For $t\sim60$, nonlinear effects become important and thus the instability speeds up. For longer times, $t\geq 70$ kinetic effects come into play and the instability saturates. The magnetic energy becomes statistically constant, while the electric energy reaches its maximum value at saturation and then starts to decrease. This is in agreement with the fact that as soon as the instability saturates and  the growth rate decreases, the wave becomes dominated by the magnetic field.
Here we  also observe that the growth rate of $E_2$ energy is about twice of the growth rate of the magnetic energy. This behavior was anticipated 
in \cite{califano1998ksw} in the context of a two-fluid model and also agrees with \cite{cheng_vm}. It is due to wave coupling and a modulation of the electron
density induced by the spatial modulation of $B_3^2$. The density modulation, including the expected spikes, is seen in Figure \ref{figure_semiimplicitaltdens}.

\begin{figure}[!htbp]
\centering
\subfigure[Log Fourier Modes of $E_1$. Run1]{\includegraphics[width=0.45\textwidth]{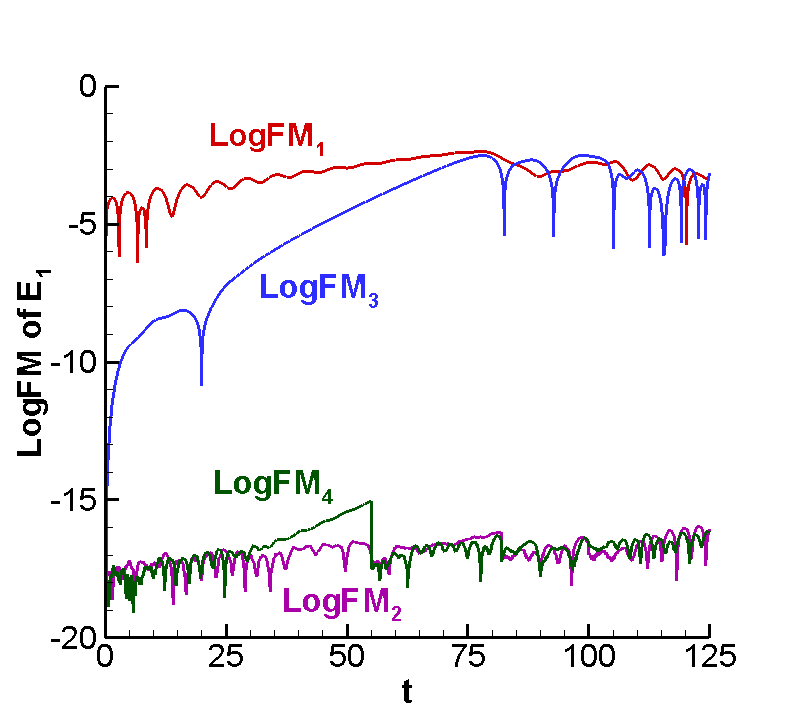}}
\subfigure[Log Fourier Modes of $E_1$. Run2]{\includegraphics[width=0.45\textwidth]{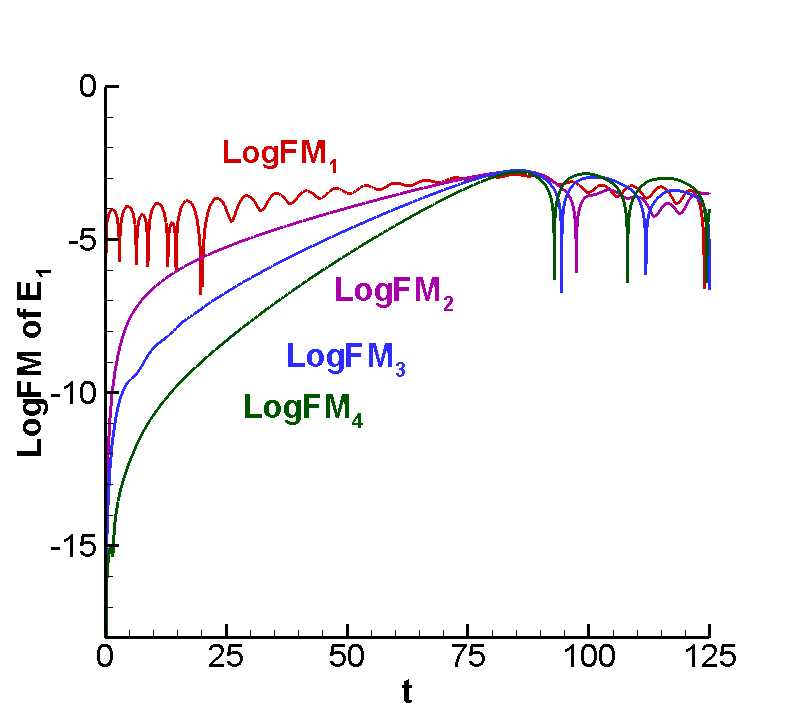}}\\
\subfigure[Log Fourier Modes of $E_2$. Run1]{\includegraphics[width=0.45\textwidth]{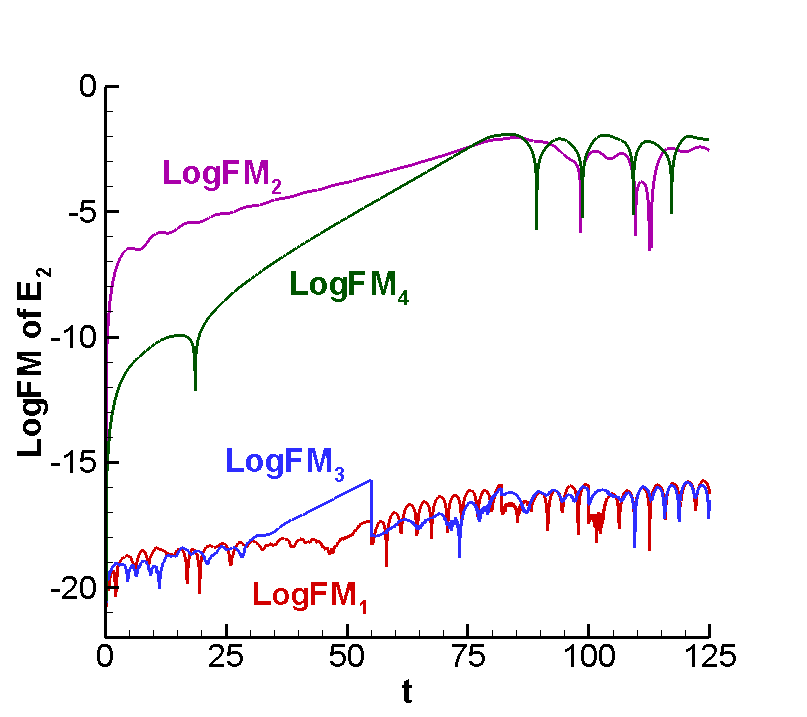}}
\subfigure[Log Fourier Modes of $E_2$. Run2]{\includegraphics[width=0.45\textwidth]{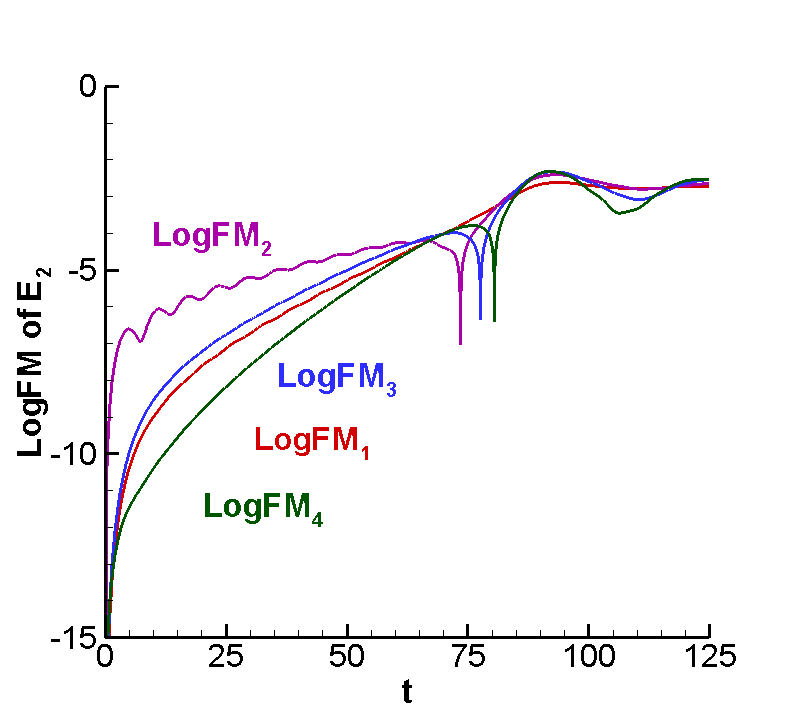}}\\
\subfigure[Log Fourier Modes of $B_3$. Run1]{\includegraphics[width=0.45\textwidth]{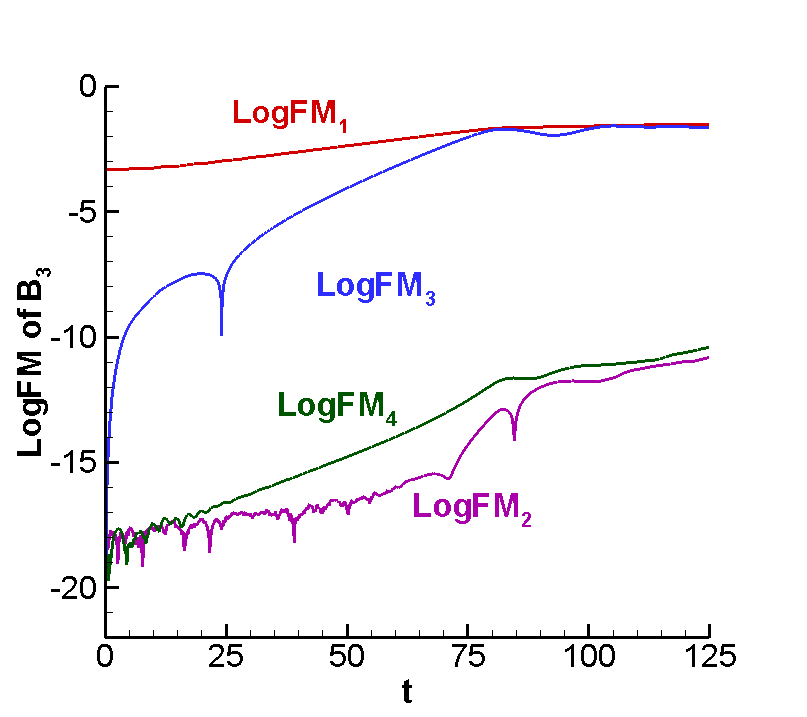}}
\subfigure[Log Fourier Modes of $B_3$. Run2]{\includegraphics[width=0.45\textwidth]{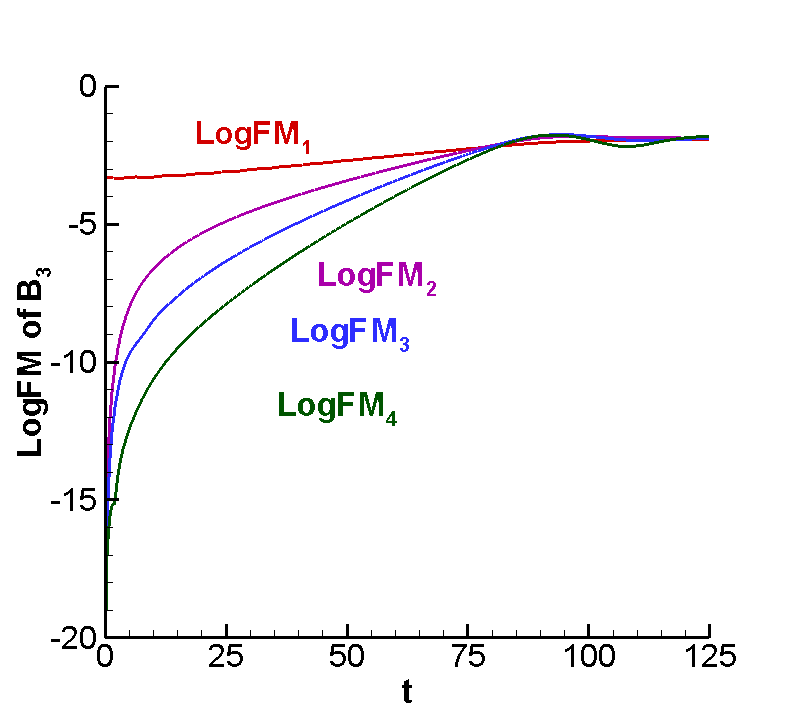}}
\caption{Log Fourier Modes. $\textnormal{\bf Scheme-2}\textnormal{ and } P^2$.  
$80^3$ mesh.
 $CFL = 0.15$. Alternating flux for Maxwell solver. }
 \label{figure_logfm}
\end{figure}

\begin{figure}[!htbp]
\centering
\subfigure[$t=55$. Run 1. ]{\includegraphics[width=0.45\textwidth]{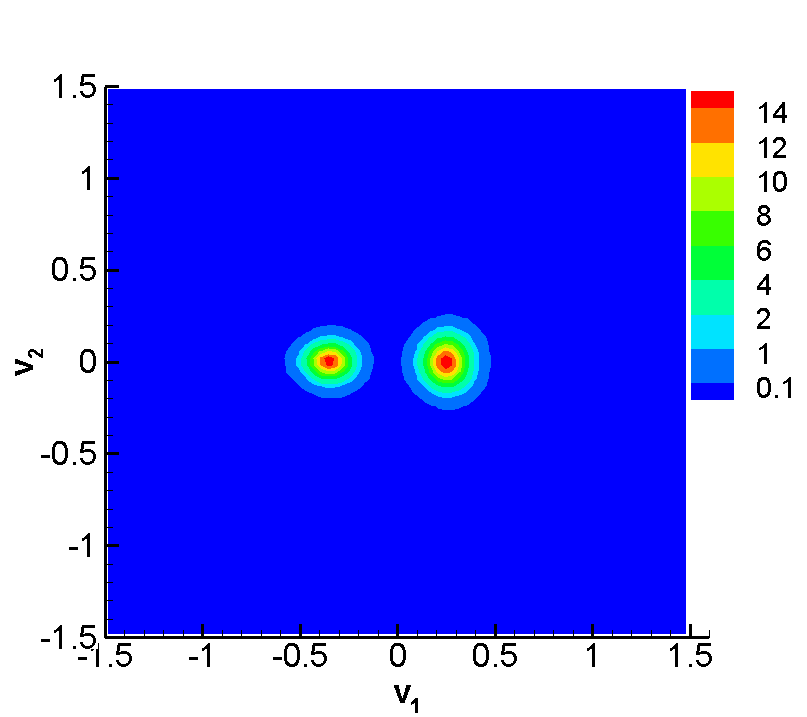}}
\subfigure[$t=55$. Run 2. ]{\includegraphics[width=0.45\textwidth]{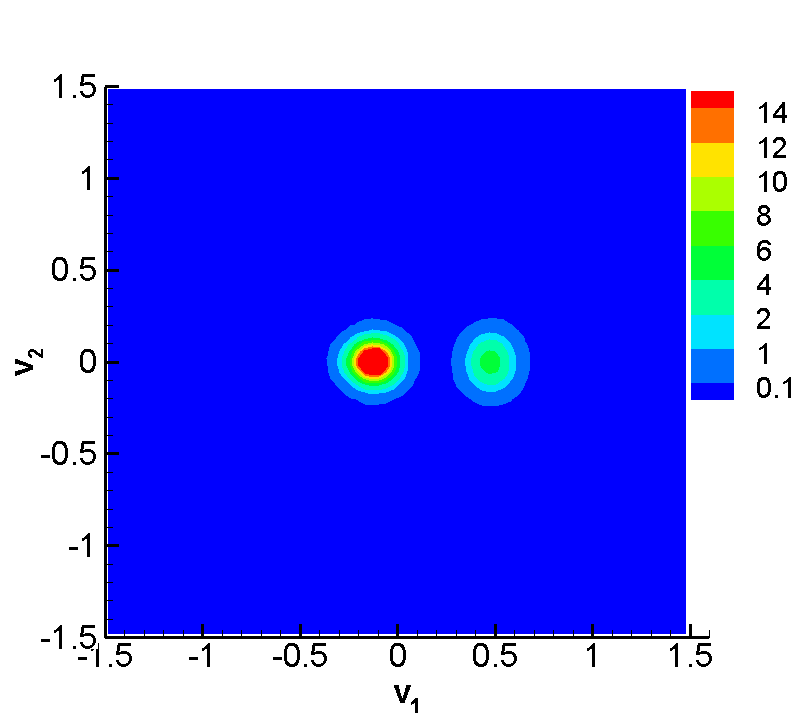}}\\
\subfigure[$t=82$. Run 1. ]{\includegraphics[width=0.45\textwidth]{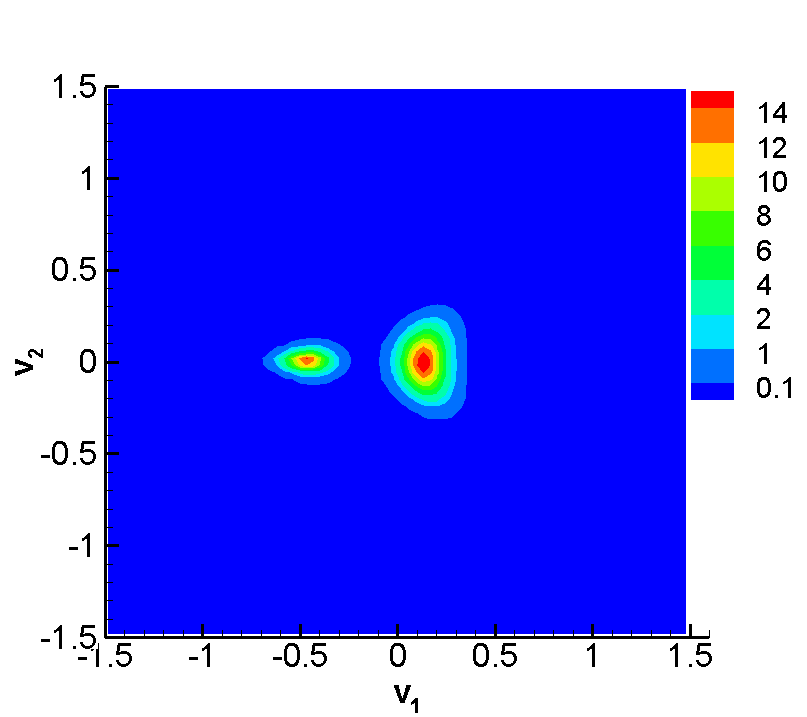}}
\subfigure[$t=82$. Run 2. ]{\includegraphics[width=0.45\textwidth]{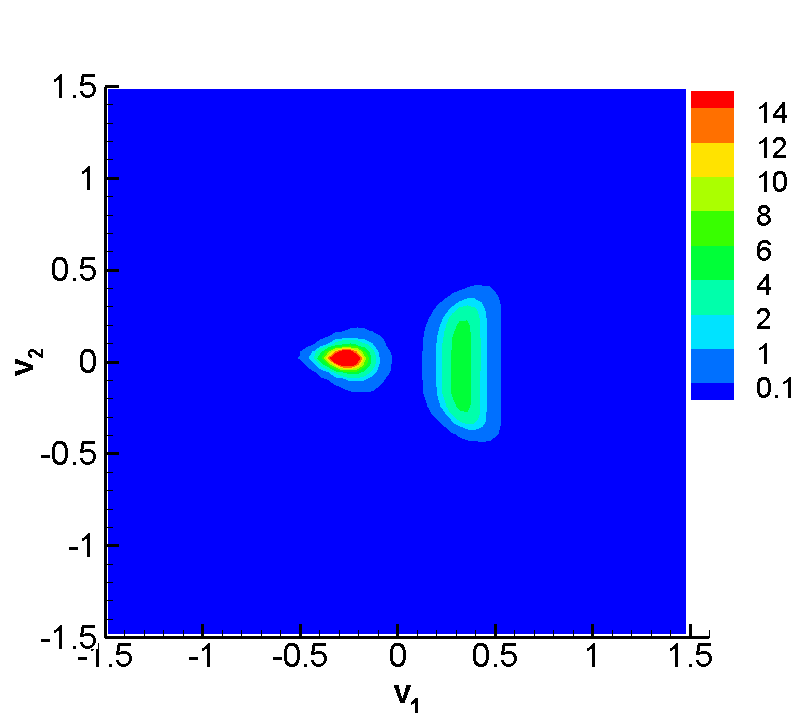}}\\
\subfigure[$t=125$. Run 1. ]{\includegraphics[width=0.45\textwidth]{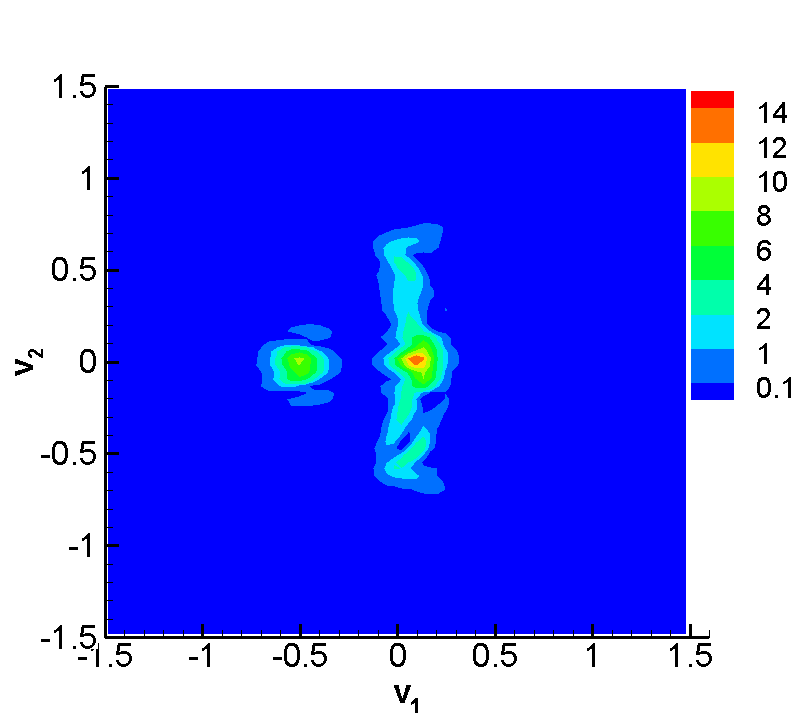}}
\subfigure[$t=125$. Run 1. ][$t=82$. Run 2. ]{\includegraphics[width=0.45\textwidth]{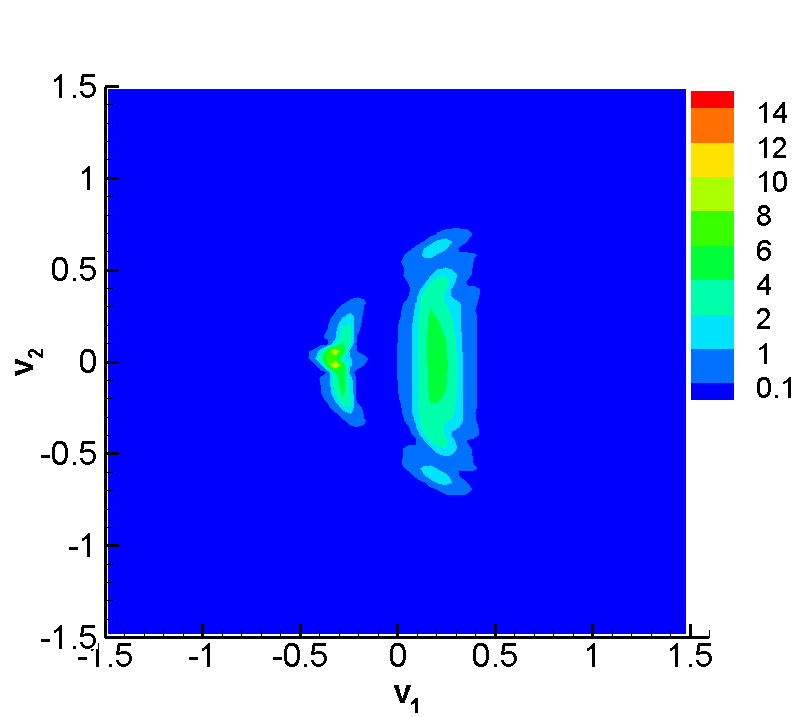}}
\caption{2D contour plots of the distribution function  at selected location $x_2=0.0625\pi$ and time $t$. $\textnormal{\bf Scheme-2}\textnormal{ and } P^2$.  $80^3$ mesh. $CFL = 0.15$. Alternating flux for Maxwell solver. }
\label{figure_semiimplicitaltx1}
\end{figure}

\begin{figure}[!htbp]
\centering
\subfigure[$t=55$. Run 1. ]{\includegraphics[width=0.45\textwidth]{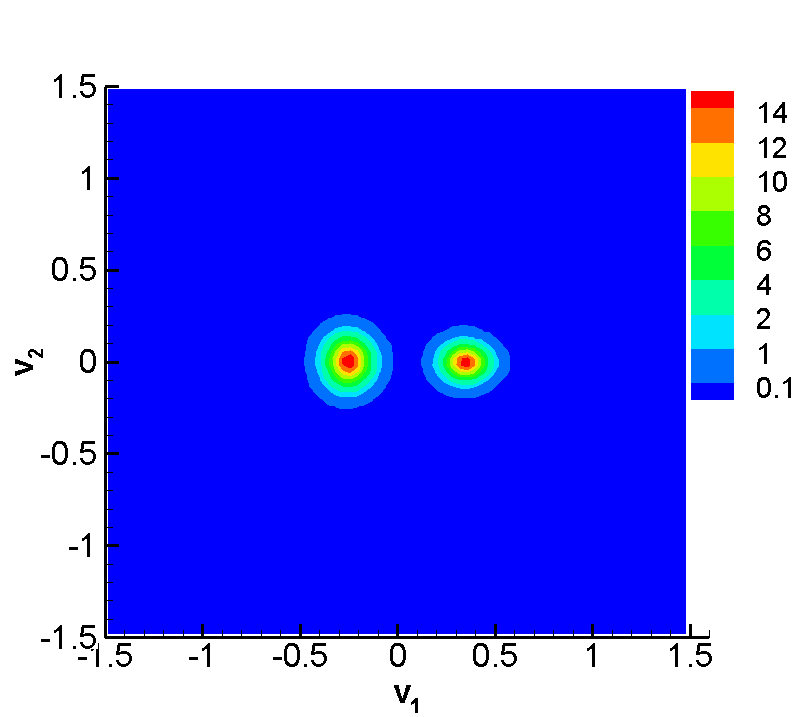}}
\subfigure[$t=55$. Run 2. ]{\includegraphics[width=0.45\textwidth]{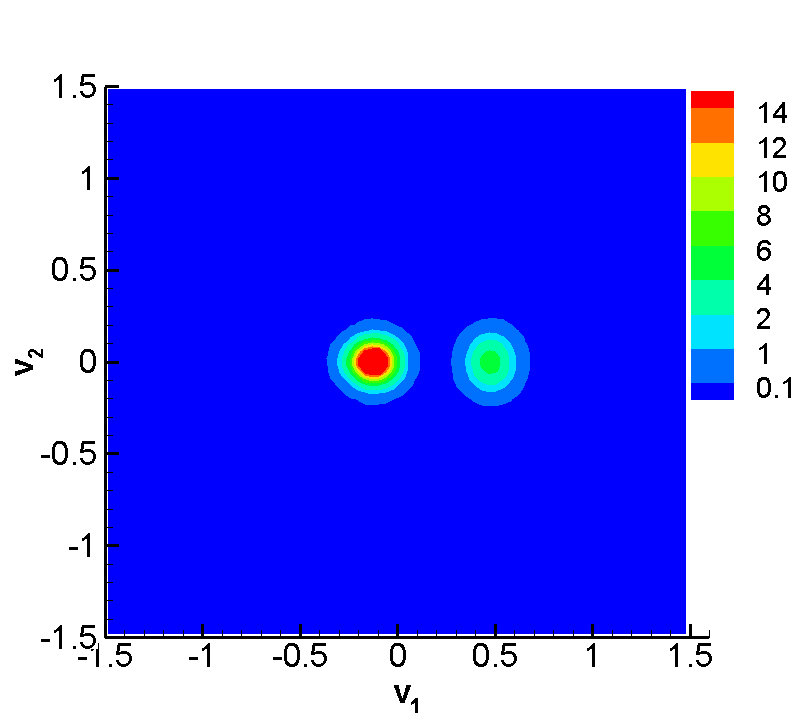}}\\
\subfigure[$t=82$. Run 1. ]{\includegraphics[width=0.45\textwidth]{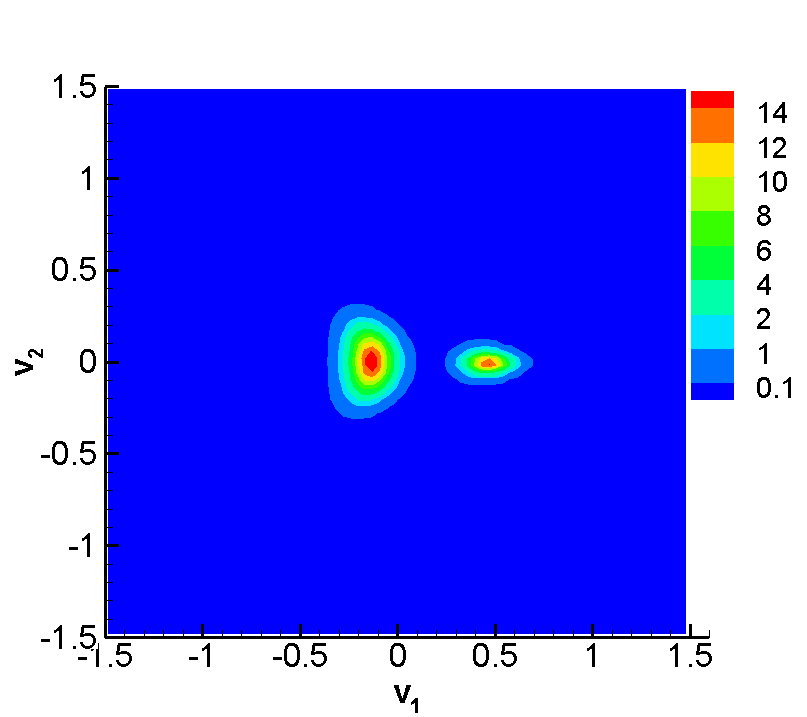}}
\subfigure[$t=82$. Run 2. ]{\includegraphics[width=0.45\textwidth]{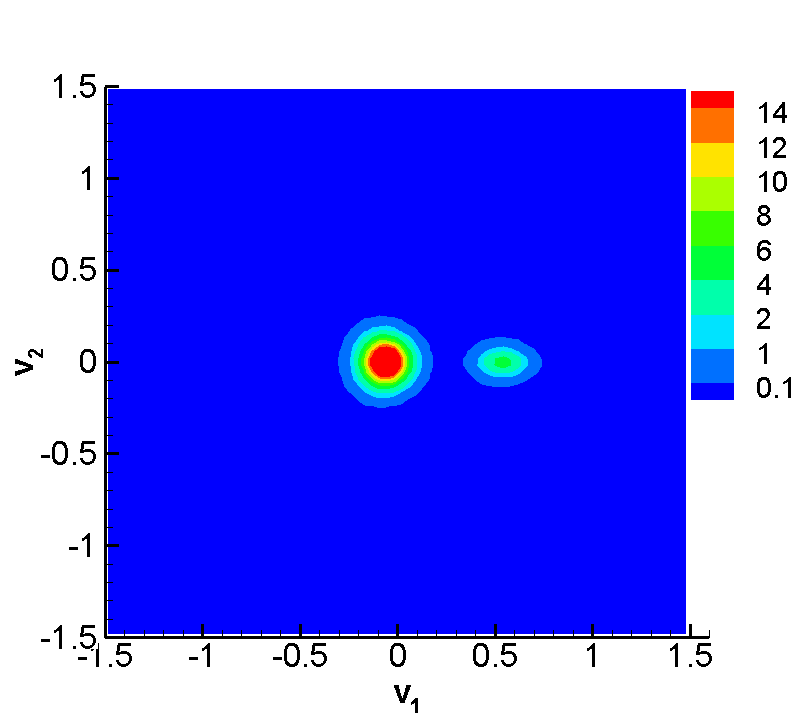}}\\
\subfigure[$t=125$. Run 1. ]{\includegraphics[width=0.45\textwidth]{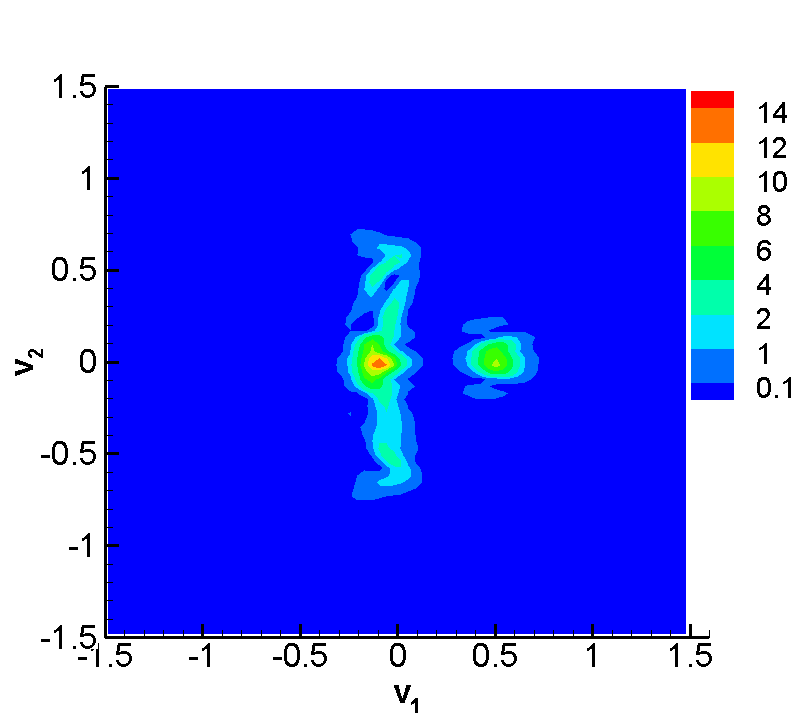}}
\subfigure[$t=125$. Run 2. ]{\includegraphics[width=0.45\textwidth]{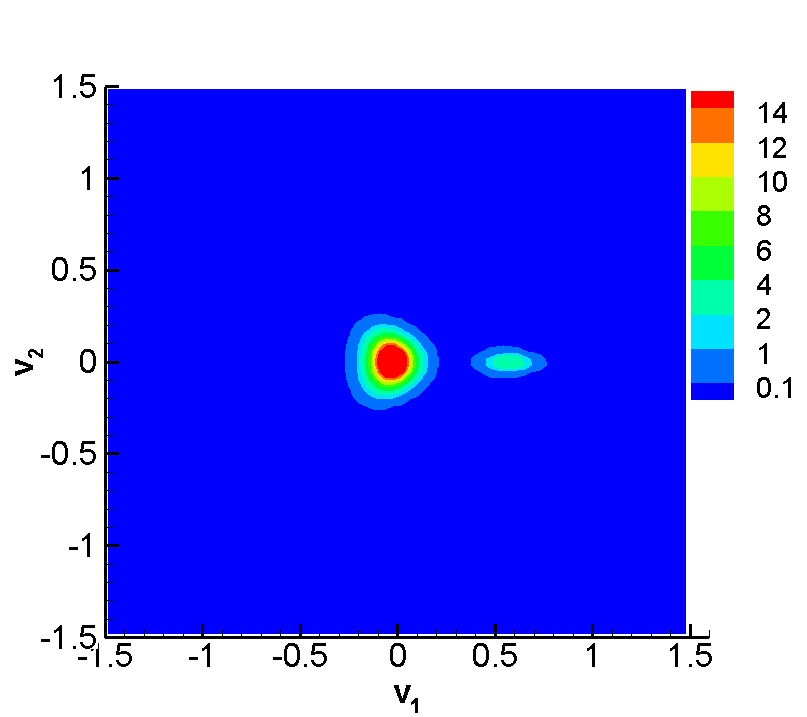}}
\caption{2D contour plots of the distribution function at selected location $x_2=4.9375\pi$ and time $t$. $\textnormal{\bf Scheme-2}\textnormal{ and } P^2$.  $80^3$ mesh. $CFL = 0.15$. Alternating flux for Maxwell solver. }
\label{figure_semiimplicitaltx2}
\end{figure}

\begin{figure}[!htbp]
\centering
\subfigure[Alternating. Run 1. ]{\includegraphics[width=0.45\textwidth]{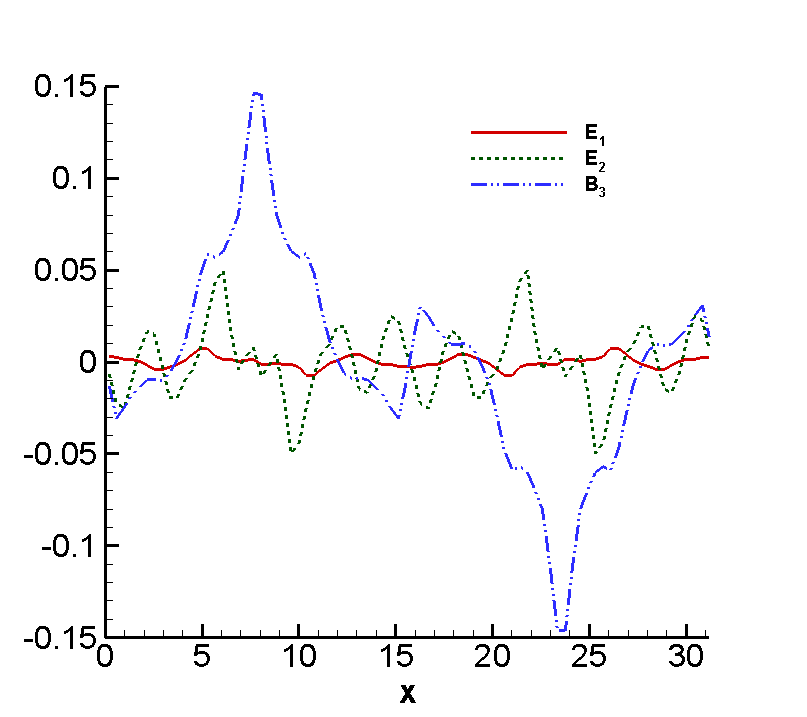}}
\subfigure[Alternating. Run 2. ]{\includegraphics[width=0.45\textwidth]{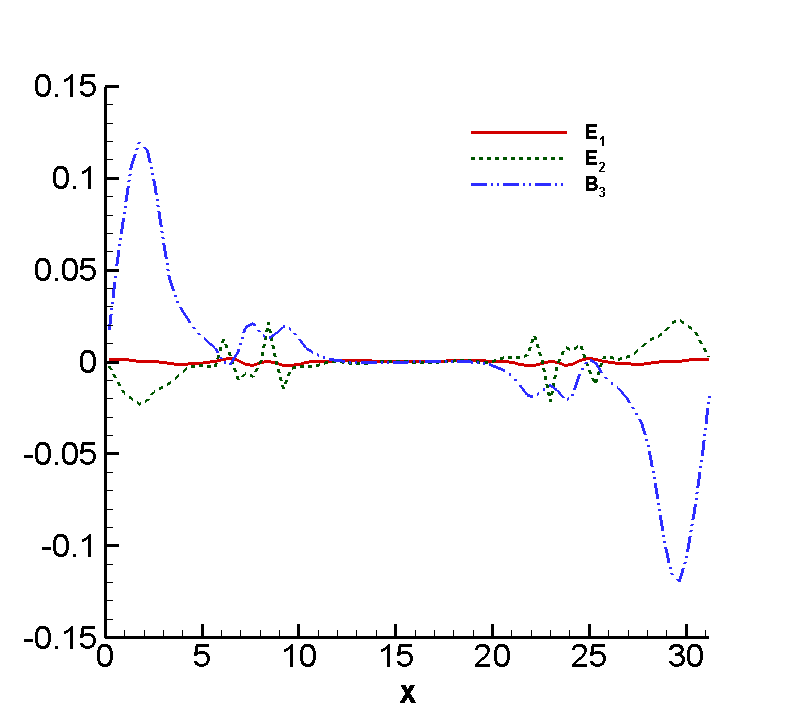}}
\caption{Plots of the electric and magnetic fields at time $t=125$.   $\textnormal{\bf Scheme-2}\textnormal{ and } P^2$.  $80^3$ mesh. $CFL = 0.15$. Alternating flux for Maxwell solver.}
\label{figure_semiimplicitebplot}
\end{figure}

In Figures \ref{figure_logfm},  we plot the first four  Log Fourier modes for the fields $E_1,\, E_2,\, B_3$  with parameter choice of Run 1 and Run 2, where
the $n$-th Log Fourier mode for a function $W(x,t)$ \cite{Heath} is defined as
$$
logF\!M_n(t)=\log_{10} \left(\frac{1}{L} \sqrt{\left|\int_0^L W(x, t) \sin(\kappa nx) \, dx \right|^2 +
\left|\int_0^L W(x, t) \cos(\kappa nx)\,  dx \right|^2} \right).
$$
Here $\kappa=\kappa_0=0.2$.
The Log Fourier modes generated by our methods agree well with \cite{cheng_vm}. 
Figures \ref{figure_semiimplicitaltx1} and \ref{figure_semiimplicitaltx2} plot the 2D contours of $f$ at  selected  time $t$ at the position $x_2=0.0625 \pi$ (near left endpoint of the domain) and the position $x_2=4.9375$ (near the middle of the domain), respectively. The times chosen correspond to those for the density of Figure \ref{figure_semiimplicitaltdens}.  
 In Figure \ref{figure_semiimplicitebplot}, We also plot  the electric and magnetic fields at the final time $t=125$ for completeness. 
All of our results are in reasonable agreement comparing with  \cite{califano1998ksw, cheng_vm}.

\section{Concluding Remarks}
\label{sec:conclusion}
In this paper, we  generalize the idea in our previous work for the VA system \cite{cheng_va} and propose energy-conserving solvers for  the VM system.
 The conservation in the total particle number and total energy is achieved on the fully discrete level in our schemes after taking additional care of both the temporal and spatial discretizations.
The main components of our methods include second order and above, explicit or implicit energy-conserving temporal discretizations, and  DG methods for Vlasov and Maxwell's equations with carefully chosen numerical fluxes.   In particular, energy-conserving operator splitting is proposed for the fully implicit schemes to treat the issue of high dimensionality.  Numerical tests such as the streaming Weibel instability are provided to demonstrate the accuracy
 and conservation of the schemes. Our next goal is to generalize the methods to multi-species systems in higher dimensions.

\section*{Acknowledgements}
 YC is supported by grants NSF DMS-1217563,  DMS-1318186, AFOSR FA9550-12-1-0343 and the startup fund from Michigan State University. AJC is supported by AFOSR grants FA9550-11-1-0281, FA9550-12-1-0343 and FA9550-12-1-0455,  NSF grant DMS-1115709
and MSU foundation SPG grant RG100059.  We gratefully acknowledge the support from Michigan Center for Industrial and Applied Mathematics.

\bibliographystyle{abbrv}
\bibliography{yingda,ref_cheng,refer,ref_cheng_plasma_2}

\begin{thebibliography}{10}

\bibitem{Ayuso2010}
B.~Ayuso, J.~A. Carrillo, and C.-W. Shu.
\newblock Discontinuous {Galerkin} methods for the multi-dimensional
  {Vlasov-Poisson} problems.
\newblock {\em Math. Models Methods Appl. Sci}.

\bibitem{Ayuso2009}
B.~Ayuso, J.~A. Carrillo, and C.-W. Shu.
\newblock Discontinuous {Galerkin} methods for the one-dimensional
  {Vlasov-Poisson} system.
\newblock {\em Kinet. Relat. Models}, 4:955--989, 2011.

\bibitem{ayuso2012high}
B.~Ayuso and S.~Hajian.
\newblock {High order and energy preserving discontinuous Galerkin methods for
  the Vlasov-Poisson system}.
\newblock 2012.
\newblock preprint.

\bibitem{barth2006role}
T.~Barth.
\newblock On the role of involutions in the discontinuous {Galerkin}
  discretization of {Maxwell} and magnetohydrodynamic systems.
\newblock In {\em IMA Volume on Compatible spatial discretizations}, pages
  69--88. Springer, 2006.

\bibitem{brackbill1982implicit}
J.~Brackbill and D.~Forslund.
\newblock An implicit method for electromagnetic plasma simulation in two
  dimensions.
\newblock {\em J. Comput. Phys.}, 46(2):271--308, 1982.

\bibitem{califano2001ffm}
F.~Califano, N.~Attico, F.~Pegoraro, G.~Bertin, and S.~Bulanov.
\newblock {Fast formation of magnetic islands in a plasma in the presence of
  counterstreaming electrons}.
\newblock {\em Phys. Rev. Lett.}, 86(23):5293--5296, 2001.

\bibitem{califano1965ikp}
F.~Califano, F.~Pegoraro, and S.~Bulanov.
\newblock {Impact of kinetic processes on the macroscopic nonlinear evolution
  of the electromagnetic-beam-plasma instability}.
\newblock {\em Phys. Rev. Lett.}, 84:3602, 1965.

\bibitem{califano1998ksw}
F.~Califano, F.~Pegoraro, S.~Bulanov, and A.~Mangeney.
\newblock {Kinetic saturation of the Weibel instability in a collisionless
  plasma}.
\newblock {\em Phys. Rev. E}, 57(6):7048--7059, 1998.

\bibitem{chen2011energy}
G.~Chen, L.~Chac{\'o}n, and D.~Barnes.
\newblock An energy-and charge-conserving, implicit, electrostatic
  particle-in-cell algorithm.
\newblock {\em J. Comput. Phys.}, 230(18):7018--7036, 2011.

\bibitem{cheng_va}
Y.~Cheng, A.~J. Christlieb, and X.~Zhong.
\newblock Energy conserving schemes for {Vlasov-Amp\`{e}re} systems.
\newblock {\em J. Comput. Phys.}, 256:630--655, 2014.

\bibitem{cheng_vm}
Y.~Cheng, I.~M. Gamba, F.~Li, and P.~J. Morrison.
\newblock {Discontinuous Galerkin schemes for Vlasov-Maxwell systems}.
\newblock {\em SIAM Journal on Numerical Analysis}, 52:1017--1049, 2014.

\bibitem{cheng_vp}
Y.~Cheng, I.~M. Gamba, and P.~J. Morrison.
\newblock Study of conservation and recurrence of {Runge--Kutta} discontinuous
  {Galerkin} schemes for {Vlasov-Poisson} systems.
\newblock {\em J. Sci. Comput.}, 56:319--349, 2013.

\bibitem{Chung2012}
E.~T. Chung, P.~Ciarlet, and T.~F. Yu.
\newblock {Convergence and superconvergence of staggered discontinuous Galerkin
  methods for the three-dimensional Maxwell's equations on Cartesian grids}.
\newblock {\em J. Comput. Phys.}, 235:14--31, 2013.

\bibitem{Cockburn_2000_history}
B.~Cockburn, G.~Karniadakis, and C.-W. Shu.
\newblock The development of discontinuous {Galerkin} methods.
\newblock In B.~Cockburn, G.~Karniadakis, and C.-W. Shu, editors, {\em
  Discontinuous {Galerkin} methods: theory, computation and applications},
  volume~11, pages 3--50. Springer, 2000.

\bibitem{Cockburn_2001_RK_DG}
B.~Cockburn and C.-W. Shu.
\newblock {Runge-Kutta} discontinuous {Galerkin} methods for
  convection-dominated problems.
\newblock {\em J. Sci. Comput.}, 16:173--261, 2001.

\bibitem{cohen1989performance}
B.~Cohen, A.~Langdon, D.~Hewett, and R.~Procassini.
\newblock Performance and optimization of direct implicit particle simulation.
\newblock {\em J. Comput. Phys.}, 81(1):151--168, 1989.

\bibitem{de1992easily}
J.~De~Frutos and J.~Sanz-Serna.
\newblock An easily implementable fourth-order method for the time integration
  of wave problems.
\newblock {\em J. Comput. Phys.}, 103(1):160--168, 1992.

\bibitem{DipernaCPAM89}
R.~DiPerna and P.-L. Lions.
\newblock {Global weak solutions of Vlasov-Maxwell systems}.
\newblock {\em Commun. Pur. Appl. Math}, 42:729--757, 1989.

\bibitem{eliasson2003numerical}
B.~Eliasson.
\newblock {Numerical modelling of the two-dimensional Fourier transformed
  Vlasov-Maxwell system}.
\newblock {\em J. Comput. Phys.}, 190(2):501--522, 2003.

\bibitem{fezoui2005convergence}
L.~Fezoui, S.~Lanteri, S.~Lohrengel, and S.~Piperno.
\newblock Convergence and stability of a discontinuous {Galerkin} time-domain
  method for the 3d heterogeneous {Maxwell} equations on unstructured meshes.
\newblock {\em ESAIM: Mathematical Modelling and Numerical Analysis},
  39(06):1149--1176, 2005.

\bibitem{Filbet_compare_2003}
F.~Filbet and E.~Sonnendr{\"u}cker.
\newblock Comparison of {Eulerian} {Vlasov} solvers.
\newblock {\em Computer Physics Communications}, 150:247--266, 2003.

\bibitem{forest1990fourth}
E.~Forest and R.~Ruth.
\newblock Fourth-order symplectic integration.
\newblock {\em Physica D: Nonlinear Phenomena}, 43(1):105--117, 1990.

\bibitem{GlasseyCMP88}
R.~Glassey and J.~Schaeffer.
\newblock {Global existence for the relativistic Vlasov-Maxwell system with
  nearly neutral initial data}.
\newblock {\em Comm. Math. Phys.}, 119:353--384, 1988.

\bibitem{GlasseyCMP97}
R.~Glassey and J.~Schaeffer.
\newblock {The ``two and one-half-dimensional" relativistic Vlasov Maxwell
  system}.
\newblock {\em Commun. Math. Phys.}, 185:257--284, 1997.

\bibitem{GlasseyARMA98I}
R.~Glassey and J.~Schaeffer.
\newblock {The relativistic Vlasov-Maxwell system in two space dimensions. I}.
\newblock {\em Arch. Ration. Mech. Anal.}, 141:331--354, 1998.

\bibitem{GlasseyARMA98II}
R.~Glassey and J.~Schaeffer.
\newblock {The relativistic Vlasov-Maxwell system in two space dimensions. II}.
\newblock {\em Arch. Ration. Mech. Anal.}, 141:355--374, 1998.

\bibitem{GlasseyARMA86}
R.~T. Glassey and W.~A. Strauss.
\newblock {Singularity formation in a collisionless plasma could occur only at
  high velocityes}.
\newblock {\em Arch. Ration. Mech. Anal.}, 92:59--90, 1986.

\bibitem{GlasseyCMP87}
R.~T. Glassey and W.~A. Strauss.
\newblock {Absence of shocks in an initially dilute collisionless plasma}.
\newblock {\em Comm. Math. Phys.}, 113:191--208, 1987.

\bibitem{hairer2006geometric}
E.~Hairer, C.~Lubich, and G.~Wanner.
\newblock {\em Geometric numerical integration: structure-preserving algorithms
  for ordinary differential equations}, volume~31.
\newblock Springer, 2006.

\bibitem{heath2012discontinuous}
R.~Heath, I.~Gamba, P.~Morrison, and C.~Michler.
\newblock {A discontinuous Galerkin method for the Vlasov-Poisson system}.
\newblock {\em J. Comput. Phys.}, 231(4):1140--1174, 2012.

\bibitem{Heath_thesis}
R.~E. Heath.
\newblock {Numerical analysis of the discontinuous Galerkin method applied to
  plasma physics}.
\newblock 2007.
\newblock Ph. D. dissertation, the University of Texas at Austin.

\bibitem{Heath}
R.~E. Heath, I.~M. Gamba, P.~J. Morrison, and C.~Michler.
\newblock A discontinuous {Galerkin} method for the {Vlasov-Poisson} system.
\newblock {\em J. Comp. Phys.}, 231:1140--1174, 2012.

\bibitem{hesthaven2007nodal}
J.~Hesthaven and T.~Warburton.
\newblock {\em Nodal discontinuous {Galerkin} methods: algorithms, analysis,
  and applications}, volume~54.
\newblock Springer, 2007.

\bibitem{hindmarsh2005sundials}
A.~C. Hindmarsh, P.~N. Brown, K.~E. Grant, S.~L. Lee, R.~Serban, D.~E.
  Shumaker, and C.~S. Woodward.
\newblock Sundials: Suite of nonlinear and differential/algebraic equation
  solvers.
\newblock {\em ACM T. Math. Software}, 31(3):363--396, 2005.

\bibitem{Hesthaven_09_divclean}
G.~Jacobs and J.~Hesthaven.
\newblock Implicit-explicit time integration of a high-order particle-in-cell
  method with hyperbolic divergence cleaning.
\newblock {\em Comput. Phys. Comm.}, 180:1760--1767, 2009.

\bibitem{Jacobs_Hesthaven_2006}
G.~B. Jacobs and J.~S. Hesthaven.
\newblock High-order nodal discontinuous galerkin particle-in-cell method on
  unstructured grids.
\newblock {\em J. Comput. Phys.}, 214:96--121, May 2006.

\bibitem{knoll2004jacobian}
D.~A. Knoll and D.~E. Keyes.
\newblock {Jacobian-free Newton-Krylov} methods: a survey of approaches and
  applications.
\newblock {\em J. Comput. Phys}, 193(2):357--397, 2004.

\bibitem{leimkuhler2005simulating}
B.~Leimkuhler and S.~Reich.
\newblock {\em {Simulating Hamiltonian dynamics}}, volume~14.
\newblock Cambridge University Press, 2005.

\bibitem{mangeney2002nsi}
A.~Mangeney, F.~Califano, C.~Cavazzoni, and P.~Travnicek.
\newblock {A numerical scheme for the integration of the Vlasov-Maxwell system
  of equations}.
\newblock {\em J. Comput. Phys.}, 179(2):495--538, 2002.

\bibitem{Markidis20117037}
S.~Markidis and G.~Lapenta.
\newblock The energy conserving particle-in-cell method.
\newblock {\em J. Comput. Phys.}, 230(18):7037 -- 7052, 2011.

\bibitem{mclachlan2002splitting}
R.~McLachlan and G.~Quispel.
\newblock Splitting methods.
\newblock {\em Acta Numerica}, 11(0):341--434, 2002.

\bibitem{Munz:2000}
C.-D. Munz, P.~Omnes, R.~Schneider, E.~Sonnendr\"{u}cker, and U.~Vo$\beta$.
\newblock {Divergence Correction Techniques for Maxwell Solvers Based on a
  Hyperbolic Model}.
\newblock {\em J. Comput. Phys.}, 161:484--511, 2000.

\bibitem{piperno_2006}
S.~Piperno.
\newblock Symplectic local time-stepping in non-dissipative {DGTD} methods
  applied to wave propagation problems.
\newblock {\em ESAIM: Mathematical Modelling and Numerical Analysis},
  40(05):815--841, 2006.

\bibitem{piperno2002nondiffusive}
S.~Piperno, M.~Remaki, and L.~Fezoui.
\newblock A nondiffusive finite volume scheme for the three-dimensional
  {Maxwell's} equations on unstructured meshes.
\newblock {\em SIAM J Numer Anal}, 39(6):2089--2108, 2002.

\bibitem{qiu2011positivity}
J.~Qiu and C.~Shu.
\newblock {Positivity preserving semi-Lagrangian discontinuous Galerkin
  formulation: Theoretical analysis and application to the Vlasov-Poisson
  system}.
\newblock {\em J. Comput. Phys.}, 230(23):8386--8409, 2011.

\bibitem{rodrigue2001vector}
G.~Rodrigue and D.~White.
\newblock A vector finite element time-domain method for solving {Maxwell's}
  equations on unstructured hexahedral grids.
\newblock {\em SIAM J. Sci. Comput.}, 23(3):683--706, 2001.

\bibitem{rossmanith2011positivity}
J.~Rossmanith and D.~Seal.
\newblock {A positivity-preserving high-order semi-Lagrangian discontinuous
  Galerkin scheme for the Vlasov-Poisson equations}.
\newblock {\em J. Comput. Phys.}, 230(16):6203--6232, 2011.

\bibitem{sanz1991order}
J.~Sanz-Serna and L.~Abia.
\newblock Order conditions for canonical {Runge-Kutta} schemes.
\newblock {\em SIAM J Numer Anal}, 28(4):1081--1096, 1991.

\bibitem{strang1968construction}
G.~Strang.
\newblock On the construction and comparison of difference schemes.
\newblock {\em SIAM J Numer Anal}, 5(3):506--517, 1968.

\bibitem{taflove2000computational}
A.~Taflove and S.~Hagness.
\newblock Computational electrodynamics: the {FDTD} method.
\newblock {\em Artech House Boston, London}, 2000.

\bibitem{umeda2009two}
T.~Umeda, K.~Togano, and T.~Ogino.
\newblock Two-dimensional full-electromagnetic {Vlasov} code with conservative
  scheme and its application to magnetic reconnection.
\newblock {\em Comput. Phys. Commun.}, 180(3):365--374, 2009.

\bibitem{yee1966numerical}
K.~Yee.
\newblock Numerical solution of initial boundary value problems involving
  {Maxwell's} equations in isotropic media.
\newblock {\em IEEE T. Antenn. Propag.}, 14(3):302--307, 1966.

\bibitem{yoshida1990construction}
H.~Yoshida.
\newblock Construction of higher order symplectic integrators.
\newblock {\em Phys. Lett. A}, 150(5):262--268, 1990.

\end{thebibliography}

\end{document}